\documentclass[10pt]{article}

\usepackage{bm}
\usepackage{fullpage}  

\usepackage{url}
\usepackage{algorithm, algorithmicx,algpseudocode}
\usepackage{multirow} 
\usepackage{hhline}  
\usepackage{verbatim}
\usepackage{bbm}
\usepackage{graphicx}
\usepackage{amsmath,amssymb,amsfonts,graphicx,amsthm,mathtools,nicefrac}
\usepackage{lscape}
\usepackage{color}
\usepackage{authblk}
\usepackage{enumerate}
\usepackage{subcaption}
\usepackage{booktabs}

\allowdisplaybreaks
\newcommand\numberthis{\addtocounter{equation}{1}\tag{\theequation}}

\newtheorem{theorem}{Theorem}[section]
\newtheorem{assumption}{Assumption}[section]
\newtheorem{remark}{Remark}[section]

\numberwithin{equation}{section}

\newtheorem{lemma}[theorem]{Lemma}
\newtheorem{corollary}[theorem]{Corollary}

\DeclareMathOperator{\diag}{\mathrm{diag}}

\DeclareMathOperator{\out}{out} 
\DeclareMathOperator{\Var}{Var} 
 
\DeclareMathOperator{\KL}{KL} 
\DeclareMathOperator{\Geo}{Geo}

\def\la {\left\langle}
\def\ra {\right\rangle} 

\newcommand{\matsnorm}[2]{\left\| #1\right\|_{{#2}}}

\newcommand{\opnorm}[1]{\ensuremath{\matsnorm{#1}{}}}

\newcommand{\twonorm}[1]{\ensuremath{\matsnorm{#1}{\footnotesize{2}}}}

\newcommand{\bfm}[1]{\bm{#1}}

\newcommand{\E}[2][]{\mathbb{E}_{#1} \left\{ #2 \rule{0mm}{3mm}\right\}}
 
\def\va{\bfm a}   \def\mA{\bfm A}  
\def\vb{\bfm b}   \def\mB{\bfm B}  
     
\def\vd{\bfm d}     
     
  \def\mF{\bfm F}  
\def\vg{\bfm g}     
   \def\mH{\bfm H}  
   \def\mI{\bfm I}  
   \def\mJ{\bfm J}

     \def\R{\mathbb{R}}
\def\vs{\bfm s}

\def\vv{\bfm v}     
   \def\mW{\bfm W}  
\def\vx{\bfm x}     
\def\vy{\bfm y}

\def\calA{{\cal  A}}

\def\calD{{\cal  D}} 
\def\calE{{\cal  E}} 
 
\def\calG{{\cal  G}}

\def\calN{{\cal  N}} 
\def\calO{{\cal  O}}

\def\calS{{\cal  S}}

\newcommand{\bfsym}[1]{\bm{#1}}

\def\btheta{\bfsym {\theta}}

\def \by{\bar{\vy}}

\def \bv{\bar{\vv}}

\def \tran {\mathsf{T}}

\def \bzero{\bm 0}
\def \bone{\bm 1}

\usepackage{bbm}
\def \one {\mathbbm{1}}

\def \bnablaf {\overline{\nabla V}}

\def \pW {\left(\mW\otimes \mI_d \right)}
\def \bbtheta{\bar{\btheta}}
\def \bx{\bar{\btheta}}
\def \bd{\bar{\vd}}
\def\vx{\btheta}

\def \pJ { \frac{1}{n}\left(\mJ_n \otimes \mI_d \right)}
\def \bnablaV {\overline{\nabla V}}
\def \tnablaV {\widetilde{\nabla V}}

\def \vvx {\bfm x}

\begin{document}
	
\title{Decentralized Natural Policy Gradient with Variance Reduction for Collaborative Multi-Agent Reinforcement Learning\footnotetext{Jinchi Chen and Jie Feng contribute equally to this work.}}
\author[1]{Jinchi Chen}
\author[1]{Jie Feng}
\author[1,2]{Weiguo Gao}
\author[1]{Ke Wei}
\affil[1]{School of Data Science, Fudan University, Shanghai, China.\vspace{.15cm}}
\affil[2]{School of Mathematical Sciences, Fudan University, Shanghai, China.}


\date{\today}

\maketitle
\begin{abstract}
This paper studies a policy optimization problem arising from collaborative multi-agent reinforcement learning in a decentralized setting where agents communicate with their neighbors over an undirected graph to maximize the sum of their cumulative rewards. A novel decentralized natural policy gradient method, dubbed Momentum-based Decentralized Natural Policy Gradient (MDNPG), is proposed, which incorporates natural gradient, momentum-based variance reduction, and gradient tracking into the decentralized stochastic gradient ascent framework. The  $\calO(n^{-1}\epsilon^{-3})$ sample complexity for  MDNPG to converge to an $\epsilon$-stationary point has been established under standard assumptions, where $n$ is the number of agents. It indicates that  MDNPG can achieve the optimal convergence rate for decentralized policy gradient methods and possesses a linear speedup in contrast to  centralized optimization methods. Moreover, superior empirical performance of MDNPG over other state-of-the-art algorithms has been demonstrated by extensive numerical experiments.
\end{abstract}

\section{Introduction}

Reinforcement learning (RL) is a sequential decision-making task in which an agent seeks a strategy that maximizes the long-term return received from the environment via interaction with the system. Recent years have witnessed considerable theoretical and empirical advances in RL, see for example \cite{kaelbling1996reinforcement, arulkumaran2017deep, rajeswaran2020game, agarwal2021theory} and references therein. In particular, when combined with deep learning, RL has achieved the most recent state of the art in various data-driven applications, including robotics~\cite{kober2013reinforcement}, finance~\cite{liu2020finrl} and game playing~\cite{mnih2013playing}.

Markov decision processes (MDPs) are  widely used to model how agents interact with an environment. An MDP can be defined as a tuple $ \la\calS,\calA, P, r,\gamma\ra$, where $\calS$ is a finite state space,  $\calA$ is a finite action space,  $P: \calS \times \calA \rightarrow \Delta(\calS)$ is the state transition model which determines the probability from $(\vs, \va)$ to state $\vs'$, $r$: $\calS\times\calA\times\calS\rightarrow [-1,1]$ is  the immediate reward function associated with the transition from $(\vs,\va)$ to $\vs'$, and $\gamma\in[0,1)$ is the discount factor. Moreover, a  policy, denoted $\pi:\calS\rightarrow\Delta(\calA)$, specifies a decision-making strategy, that is, $\pi(\va|\vs)$ is the probability of executing action $\va$ at state $\vs$. Given an initial state distribution $\rho(\vs_0)$,  let $\tau = (\vs^0, \va^0, r^0, \vs^1, \va^1, r^1, \cdots, \vs^{H-1}, \va^{H-1}, r^{H-1}, \vs^H)$ be a trajectory of time horizon $H$ induced by a policy $\pi$, where $r^h=r(\vs^h,\va^h,\vs^{h+1})$. 
The overall goal of RL is to find a policy that maximizes the expected  discounted cumulative rewards, which can be formulated as the following optimization problem:
\begin{align}
\label{eq: rleq1}
    \max_{\pi\in\Pi} \left\{ V(\pi) := \E[\tau \sim p(\cdot|\pi)]{R(\tau) }\right\},
\end{align}
where $R(\tau)=\sum_{h=0}^{H-1}\gamma^h r^h$ is the discounted return obtained from the trajectory $\tau$,  $p(\cdot|\pi)$ is the distribution of the trajectories, and $\Pi$ represents the policy space.


There are several classical categories of RL algorithms. Model-based approaches, such as policy iteration and value iteration  (see e.g., \cite{puterman2014markov}) find the optimal policy based on the ideas of of fixed point iteration. Whereas in model-free settings, value-based methods, like temporal difference learning and Q-learning (see e.g., \cite{sutton2018reinforcement,bertsekas2019reinforcement}) solely use reward obtained from the environment to seek the optimal strategy. These methods can be roughly thought of as approximate dynamic programming with Monte Carlo learning. In contrast, policy gradient methods~\cite{williams1992simple,sutton1999policy,konda1999actor} maximize the objective function  in \eqref{eq: rleq1} by gradient ascent with a differentiable parameterized policy in the model-free manner. Gradient-based approaches have a few advantages. For example, they can generate stochastic policies, which are more exploratory and are easily extended to continuous control problems. Coupled with neural networks, they have gained tremendous success in many applications due to their flexibility and adaptability. Moreover, the theoretical guarantees for gradient-based methods are relatively more complete, even in conjunction with simple function approximations~\cite{sutton1999policy, xu2020improving, agarwal2021theory}.  

In this paper, we restrict our attention to policy optimization based methods. Using a parameterized policy $\pi_{\btheta}$ where $\btheta\in \mathbb{R}^d$,~\eqref{eq: rleq1} can be expressed as a finite dimensional optimization problem:
\begin{align}
\label{eq: rleq2}
\max_{\btheta\in\R^d} \left\{ V(\btheta) := \E[\tau\sim p(\cdot | \btheta)]{R(\tau)} \right\}.
\end{align}
After parameterization, the distribution of the trajectories, denoted $p(\tau|\btheta)$, is given by
\begin{align}
\label{dist of tau}
	p(\tau| \btheta) := \rho(\vs^0) \prod_{h=0}^{H-1} \pi_{\btheta}(\va^h| \vs^h) P(\vs^{h+1} | \vs^h, \va^h),
\end{align}
where we recall that $\rho(\vs^0)$ is the initial state distribution.

A direct method for solving problem~\eqref{eq: rleq2} is  policy gradient (PG). Despite its simplicity, PG is not invariant to reparameterization. As an alternative, natural policy gradient (NPG) methods~\cite{kakade2001natural,bagnell2003covariant,peters2008natural,bhatnagar2007incremental} utilize the intrinsic distance between policies, i.e., the Kullback-Leibler (KL) divergence, to modify the search direction so that parameterization invariance can be preserved. As two variants of NPG methods, trust region policy optimization (TRPO)~\cite{schulman2015trust} combines NPG with a line search procedure to guarantee improvement, whereas proximal policy optimization (PPO)~\cite{schulman2017proximal} uses a simplified objective with a penalty term or a clipped ratio rather than the KL constraint.

\subsection{Collaborative multi-agent reinforcement learning} \label{section: marl}
More recently, there has been a growing interest in multi-agent reinforcement learning (MARL) which allows agents to address problems simultaneously in more complicated settings, such as fully cooperative, fully competitive, and mixed of the two~\cite{busoniu2008comprehensive,nowe2012game,zhang2021multi}. MARL arises in many applications, including autonomous driving~\cite{shalev2016safe}, game playing~\cite{vinyals2019grandmaster}, and wireless networks~\cite{yao2019collaborative}. In this paper, we study an $n$-agent fully cooperative setting in which the goal of agents is to cooperatively maximize the global value function defined as follows:
\begin{align*}
	\max_{\btheta\in\mathbb{R}^d} \left\{ V(\btheta):=\frac{1}{n}\sum_{i=1}^{n}V_i(\btheta)\right\},
	\numberthis\label{eq: marl1}
\end{align*}
where $\btheta\in\mathbb{R}^d$ is the parameter of policy and $V_i(\btheta)$ is the value function of the $i$-th agent. 
Let $\tau_i = (\vs^0, \va^0, r_i^0, \vs^1, \va^1, r_i^1, \cdots, \vs^{H-1}, \va^{H-1}, r_i^{H-1}, \vs^H)$ be the trajectory induced by the policy $\pi_{\btheta}$ for the $i$-th agent$\footnote{For ease of notation, we drop the subscript $i$ for $\vs_i^h$ and $\va_i^h$ but only keep the subscript for $r_i^h$.}$. 
The  discounted return $R(\tau_i)$ of the $i$-th agent over trajectory $\tau_i$ is given by \begin{align*}
R(\tau_i) = \sum_{h=0}^{H-1} \gamma_i^h r_i^h.
\end{align*}
Therefore, $V_i(\btheta)$ in~\eqref{eq: marl1} has the following expression:
\begin{align*}
    V_i(\btheta):=\E[\tau_i\sim p(\cdot | \btheta)]{R(\tau_i)}. 
\end{align*}

Problem~\eqref{eq: marl1} can be used to model different cooperative MARL settings. Here we give two examples.
\paragraph{Collaborative reinforcement learning.}
In the collaborative RL setting, agents aim to maximize the sum of their cumulative rewards in a global environment~\cite{zhang2018fully, jiang2021mdpgt}. Consider an $n$-agent MDP denoted by a tuple $\la \calS, \{\calA_i\}_{i=1}^n, P, \{r_i\}_{i=1}^n, \gamma_i\ra$, where 
\begin{itemize}
    \item $\calS$ is the global state space shared by all agents, 
    \item $\calA:=\calA_1 \times \cdots \times \calA_n$ is the joint action space of all agents,
    \item $P:\calS\times \calA \rightarrow \Delta(\calS)$ is the state transition model,
    \item $r_i:\calS \times \calA \times \calS \rightarrow[-1,1]$ is the reward function of agent $i$,
    \item $\gamma_i$ is the discount factor for the $i$-th agent.
\end{itemize}
Let $\vs\in\calS$ be the global state, $\va = (\va_1,\cdots, \va_n)\in \calA$ be the joint action, and $\va_i\in\calA_i$ be the local action executed by the $i$-th agent. We assume that the state $\vs$ and the action $\va$ are observed globally whereas the reward $r_i$ is locally observable.  Define $\pi:\calS\rightarrow\Delta(\calA)$ as a joint policy, where $\pi(\va|\vs)$ specifies the probability that the agents select action $\va$ at state $\vs$. Since each agent makes decisions independently, we have that $\pi(\va|\vs) = \prod_{i=1}^{n} \pi_i(\va_i| \vs)$. Further, suppose that  the $i$-th policy is parameterized by $\btheta_{[i]}\in\R^{d_i}$, denoted $\pi_{\btheta_{[i]}}$. The probability of executing $\va$ at sate $\vs$ can be rewritten as
\begin{align}
\label{eq: pi1}
	\pi_{\btheta}(\va | \vs) := \prod_{i=1}^{n} \pi_{\btheta_{[i]}}(\va_i | \vs).
\end{align}
In this scenario, $\btheta$ in \eqref{eq: marl1}  is given by $\btheta = \begin{bmatrix}
\btheta_{[1]}^\tran &\cdots & \btheta_{[n]}^\tran
\end{bmatrix}^\tran \in\R^{d}$ and $d = \sum_{i=1}^{n}d_i$. Notice that  each  reward $r_i$ not only relies on the local parameter but also relies on the parameters of other agents.

\paragraph{Multi-task reinforcement learning.}
Multi-task reinforcement learning (MTRL) refers to the problem in which different agents learn a shared policy in different but similar environments~\cite{zeng2021decentralized}. MTRL can utilize similarities across different environments to enhance learning efficiency and generalization. Such an approach has received a lot of attention in recent years~\cite{wilson2007multi,teh2017distral, crawshaw2020multi}. In the MTRL setting, the MDP for the $i$-th agent is expressed as $\la\calS_i, \calA, P_i, r_i, \gamma_i\ra$. The setup for different agents can differ in terms of:
\begin{itemize}
    \item $\calS_i$, the state space for the $i$-th environment (similar or overlapping), 
    \item $P_i: \calS_i\times\calA\rightarrow\Delta(\calS_i)$, the transition model for the $i$-th environment,
    \item $r_i: \calS_i\times\calA \times \calS_i \rightarrow[-1,1]$, the reward function for the $i$-th agent,
    \item $\gamma_i$, the discount factor for the $i$-th agent.
\end{itemize}  
Note that in order for the agents to share a common policy $\pi$, they must have the same action space $\calA$ and the states in each state space $\calS_i$ must have the same format. In MTRL, though agents independently select actions in their own environments, the policy should be parameterized using a single parameter $\btheta\in\R^d$, yielding the parameterized policy $\pi_{\btheta}$ . 
 \\

Roughly speaking, collaborative RL shares a global state space, while MTRL shares a common action space. For conciseness, we refer to both of  the aforementioned  settings as collaborative multi-agent reinforcement learning. It should be easy to see whether ``collaborative'' refers to two tasks or the particular one task from the context.


\subsection{Decentralized optimization setup}
Since each agent only has access to local information, solving problem~\eqref{eq: marl1} needs to aggregate all local computations  to update the learning parameter. The centralized optimization method uses a central coordinator for data collection and information transmission, inevitably leading to high communication costs. Moreover, the central coordinator does not exist or may be too expensive to deploy in real applications. By contrast, in a decentralized framework that is considered in this paper each agent only communicates with its neighbors through a communication network. Let $\calG=(\calN, \calE)$ be the communication network which is indeed an undirected graph, where $\calN=\{1,\cdots,n\}$ is the set of agents, and $\calE\subseteq \calN\times\calN$ is the collection of edges. Note that a pair $(i,j)\in \calE$ represents that $i$ can communicate with $j$. For the $i$-th agent, define the set of its neighbors  as $\calN(i)=\{j\in \calN |(i,j)\in\calE \text{ or } i=j\}$. In addition, we can associate  a weight matrix $\mW=[W_{ij}]\in\mathbb{R}^{n\times n}$ with the graph $\calG$, where $W_{ij}>0$ if  $(i,j)\in \calE$, and $W_{ij}=0$ otherwise.  
Assuming $\calG$ is a connected graph, it is not hard to see that problem~\eqref{eq: marl1} is equivalent to
\begin{align*}
	\max_{\btheta_1,\cdots, \btheta_n\in\mathbb{R}^d} &  \quad  \frac{1}{n}\sum_{i=1}^{n}V_i(\btheta_i),  \\
	\text{subject to } & \quad \btheta_i=\btheta_j, \quad 
	\text{for all } (i,j) \in \calE.
	\numberthis\label{eq: marl3}
\end{align*}


\subsection{Main contributions and outline of this paper}

The main contributions of this work are summarized as follows.
\begin{itemize}
    \item We develop a Momentum-based Decentralized Natural Policy Gradient (MDNPG) method for the collaborative MARL problem. MDNPG combines natural gradient with momentum-based variance reduction and gradient tracking to solve the decentralized optimization problem. In a nutshell, natural gradient is a gradient method with suitably chosen preconditioning. Extensive numerical experiments show that introducing this preconditioning in the decentralized setting is also able to improve the empirical performance for collaborative MARL.
    \item
   Theoretical guarantees for  MDNPG have been obtained, showing that MDNPG is able to converge to an $\epsilon$-stationary point in $\calO(n^{-1}\epsilon^{-3})$ iterations provided a mini-batch initialization. Even though the variance reduced decentralized policy gradient has been studied in the collaborative MARL scenario, the existing analysis does not apply directly to MDNPG due to the requirement for the consensus of the precondition matrices. To overcome this difficulty, a novel stochastic ascent inequality  (see Lemma \ref{lemma ascent}) has been established to handle the preconditioning for the non-convex objective in the decentralized setting. This intermediate technical result is of independent interest and may be applied to the analysis of other preconditioned stochastic first order methods in decentralized non-convex optimization.
\end{itemize}

The rest of this paper is organized as follows.
In Section~\ref{section algorithm}, we present a complete description of  MDNPG  and provide theoretical guarantees for it. In addition, more closely related works are reviewed. 
In Section~\ref{section experiments}, we compare  MDNPG  with other state-of-the-art algorithms in single-agent and multi-agent experiments, which demonstrate the efficiency of the proposed method.
The proofs of the main results and key lemmas are presented in Section~\ref{section proofs}.
Finally, in Section~\ref{section conclusions}, we conclude this paper with future research directions.

Throughout this paper,  we refer to $\mA \otimes \mB$ as the Kronecker product. We denote by $\bone_n \in\R^n$ the all-one vector (i.e., all entries of $\bone_n$ are $1$) and by $\mJ_n\in\R^{n\times n}$ the all-one matrix. The $d\times d$ identity matrix is denoted by $\mI_d$. Additionally, we denote by  $\Delta(\calS)$ (or $\Delta(\calA)$)  the probability simplex over the state (or action) space.

\section{MDNPG and  convergence results} \label{section algorithm}

The Momentum-based Decentralized Natural Policy Gradient (MDNPG) algorithm is summarized in Algorithm~\ref{alg: moentum pg}. In the algorithm, agents perform the following steps at each iteration $t$:  gradient estimator calculation, gradient tracking, and parameter update. Notice that each step is simultaneously executed by all agents but is only presented from agent $i$' s view for simplicity. Overall, there are three pillars in MDNPG, which will be detailed next. Compared with  policy gradient based decentralized optimization algorithms~\cite{jiang2021mdpgt, zeng2021decentralized, zhao2021distributed} for collaborative MARL, the key difference is in the parameter update step where a natural gradient direction is used for each agent. 

\begin{algorithm}[ht!]
	\caption{Momentum-based Decentralized Natural Policy Gradient (MDNPG)}
	\label{alg: moentum pg}
	\begin{algorithmic}
		\State Input: number of iterations $T$, horizon $H$, batch size $B$, learning rate $\eta$, momentum parameter $\beta$, initial parameter $\bbtheta^0 \in\R^d$, initial estimator $\vv_i^{-1} = \bzero \in\R^d$, initial tracker $\vy_i^0 = \bzero \in\R^d$.
		\State Initialization: $\btheta_i^0= \bbtheta^0$, $\vv_i^0= \frac{1}{B}\sum_{b=1}^{B} \vg_i( \tau_{i,b}^0 | \btheta_i^0 )$ and $\vy_i^1 = \sum_{j\in\calN(i)} W_{ij}\vv_j^0$ for $i=1,\cdots, n$, where $\{\tau_{i,b}^0\}_{b=1}^B$ represents the $B$ trajectories  i.i.d sampled from $p(\cdot | \btheta_i^0)$.
		
		\For{$t=1,2,\ldots, T$}
		\State Generate an estimator $\vv_i^t$ of $\nabla V_i(\btheta^t)$:
		\begin{align*}
			\vv_i^t = \beta \vg_i(\tau_i^t|  \btheta_i^t) + (1-\beta) \left( \vv_i^{t-1} + \vg_i(\tau_i^t | \btheta_i^t) - \omega(\tau_i^t | \btheta_i^{t-1}, \btheta_i^t)\cdot \vg_i(\tau_i^t | \btheta_i^{t-1})\right).
		\end{align*}
		\State Gradient Tracking: 
		\begin{align*}
			\vy_{i}^{t+1} = \sum_{j\in\calN(i)} W_{ij} \left(\vy^{t}_j + \vv^t_j - \vv^{t-1}_j\right).
		\end{align*}
		\State Parameter Update:
		\begin{align*}
			\btheta_i^{t+1} = \sum_{j\in\calN(i)}W_{ij} \left( \btheta_j^t + \eta \mH_j^t \vy_j^{t+1} \right).
		\end{align*}
		\EndFor
		\State Output: $\btheta_{\out}\in\R^d$ chooses randomly from $\{\vx_i^t\}_{i=1,\ldots,n, t=0,\ldots, T}$.
		
	\end{algorithmic}
\end{algorithm}

\subsection{Three pillars in MDNPG algorithm}

\subsubsection{Pillar I: Decentralized optimization}
To solve problem \eqref{eq: marl3}, each agent can first perform a local gradient update and then seek consensus with its neighbors in order to fulfill the equality constraint. This is the basic idea behind the decentralized gradient ascent method which can be expressed as
\begin{align*}
	\btheta_i^{t+1} = \sum_{j\in\calN(i)} W_{ij}(\btheta_j^t +\eta  \nabla V_j(\btheta_j^t)),
	\numberthis\label{eq: pg1}
\end{align*}
where $\eta$ represents the learning rate,  $\nabla V_j(\btheta_j^t)$ represents the gradient of $V_j(\btheta_j)$ with respect to $\btheta_j^t$, and $W_{ij}$ are the elements of the weight matrix $\mW$ associated with the communication network $\calG$. 
Note that $\mW$ here plays a role of weighted average for consensus which should satisfy certain properties (see Assumption~\ref{assumption 6}).
Despite its simplicity, the original decentralized gradient ascent suffers from slow convergence. To address this issue,  a gradient tracking technique has been developed in~\cite{li2020communication, pu2021distributed}, of which the central idea is to correct biases between local copies of $\btheta$ via tracking the average gradient, i.e., $\frac{1}{n} \sum_{i=1}^n \nabla V_i(\btheta_i)$. The modified version of update~\eqref{eq: pg1} with gradient tracking  consists of the following two steps:
\begin{align*} 
\begin{split}
\vy_{i}^{t+1} &= \sum_{j\in\calN(i)} W_{ij} \left (\vy^{t}_j + \nabla V_j(\btheta^t_j) - \nabla V_j(\btheta^{t-1}_j)\right),\\
	\btheta_i^{t+1} &= \sum_{j\in\calN(i)}W_{ij} \left( \btheta_j^t + \eta \vy_j^{t+1} \right),
\end{split}
	\numberthis\label{eq: gt}
\end{align*}
where $\vy_i$ denotes the gradient tracker for agent $i$. Simple calculation shows that ~\eqref{eq: gt} satisfies the dynamic average consensus property:
\begin{align*}
\frac{1}{n} \sum_{i=1}^n \vy_i^{t+1} = \frac{1}{n} \sum_{i=1}^n \nabla V_i(\btheta^{t}_i), \quad t\geq 1,
\end{align*}
which implies that the average of $\nabla V_i(\btheta^{t}_i)$ is dynamically tracked by the average of $\vy_i^{t+1}$. Further, it can be proved that decentralized optimization methods equipped with gradient tracking can achieve better convergence rate~\cite{li2020communication}.

\subsubsection{Pillar II: Variance reduction}
 Consider the optimization problem 
 \begin{align*} \max_{\btheta\in\R^d}\E{f(\btheta; \xi)},
 \end{align*}
 where $\xi$ represents a random variable drawn from an unknown distribution $\calD$. The stochastic gradient ascent at the $t$-th iteration is given as
\begin{align*}
    \btheta^{t+1} = \btheta^{t} + \eta \cdot \vg(\btheta^{t};\xi^{t}),
\end{align*}
where $\eta$ is the learning rate and $\vg(\btheta^{t};\xi^{t})=\nabla f(\btheta^{t};\xi^{t})$ is the gradient estimator with $\xi^{t}$ being independently sampled from $\calD$. Due to the high variance incurred by the stochastic evaluation of the gradient,  vanilla stochastic gradient methods suffer from slow convergence. Thus, in order to accelerate the methods, various variance reduction methods have been proposed and studied in the past decades, such as SVRG~\cite{johnson2013accelerating}, SAGA~\cite{defazio2014saga}, SARAH~\cite{nguyen2017sarah}, and SPIDER~\cite{fang2018spider}. More recently, a momentum-based variance reduction method \cite{cutkosky2019momentum,tran2019hybrid} is proposed, in which the gradient estimator is given by 
\begin{align*}
	\vv^t = \beta \underbrace{\vg(\btheta^t;\xi^t)}_{\text{SGD}}  + (1-\beta)  \underbrace{\left( \vv^{t-1} + \vg(\btheta^t;\xi^t) - \vg(\btheta^{t-1};\xi^t)\right)}_{\text{SARAH}},
	\numberthis\label{eq: momentum}
\end{align*}
where $\beta\in(0,1]$ is the momentum parameter. A key feature of the momentum-based method is that it is a single-loop algorithm which leverages the benefits of both the unbiased SGD estimator~\cite{bottou2012stochastic} and the novel SARAH  estimator~\cite{nguyen2017sarah}. Thus it can avoid  the high computational cost of batch gradients to reduce variance. 

In this work, we will adopt the momentum-based variance reduction method for the policy gradient estimation. The gradient of $V_i(\btheta_i)$ in~\eqref{eq: marl3} with respect to $\btheta_i$ can be computed as follows
\begin{align*}
	\nabla V_i(\btheta_i) 
	&= \nabla_{\btheta_i} \E[\tau_i\sim p(\cdot |\btheta_i)]{R(\tau_i)}\\
	&=\int_{\tau_i} \nabla_{\btheta_i} p(\tau_i|\btheta_i) R(\tau_i) d\tau_i\\
	&=\int_{\tau_i} p(\tau_i|\btheta_i)  \frac{\nabla_{\btheta_i} p(\tau_i|\btheta_i) }{p(\tau_i|\btheta_i)} R(\tau_i) d\tau_i\\
	&=\E[\tau_i\sim p(\cdot |\btheta_i)]{\nabla_{\btheta_i} \log p(\tau_i|\btheta_i) R(\tau_i) }.
	\numberthis\label{eq: grad1}
\end{align*}
 The commonly used gradient estimators of policy gradient include REINFORCE \cite{williams1992simple} or GPOMDP \cite{baxter2001infinite}. For the $i$-th agent, we adopt REINFORCE with a  baseline $b_i$ as the policy gradient estimator: 
\begin{align*}
	\vg_i(\tau_i|\btheta_i) = \left[\sum_{h=0}^{H-1} \nabla_{\btheta_i} \log \pi_{\btheta_i}(\va^h| \vs^h)  \right] \cdot \left[ \sum_{h=0}^{H-1} \gamma^h r^h_i - b_i\right],
	\numberthis\label{eq: grad2}
\end{align*}
where $\tau_i$ denotes a trajectory generated under policy $\pi_{\btheta_i}$. {Without loss of generality, we use the same $\gamma$ for all agents.}

Note that in the momentum-based gradient estimator~\eqref{eq: momentum} for the ordinary stochastic optimization, the $\xi^t$ sampled from $\mathcal{D}$ is independent of $\btheta^t$. However, in \eqref{eq: grad2}, the sampled trajectory $\tau^t_i$ is determined by the distribution $p(\cdot|\btheta^{t}_i)$. It is easily seen that $\mathbf{g}_i(\tau^t_i|\btheta^{t-1}_i)$ is a biased estimator for $\nabla V(\btheta^{t-1}_i)$.  To ensure the unbiased property, we can utilize the importance sampling technique,
\begin{align*}
	\E[\tau_i^t\sim p(\cdot|\btheta^{t}_i)]{\omega(\tau_i^t |\btheta^{t-1}_i, \btheta^t_i)  \mathbf{g}_i(\tau_i^t|\btheta^{t-1}_i ) } =\nabla V_i (\btheta^{t-1}_i),
\end{align*}
 where $\omega(\tau^t_i |\btheta^{t-1}_i, \btheta^t_i)$ represents the importance weight defined as
\begin{align}
	\label{eq important sampling}
\omega(\tau^t_i |\btheta^{t-1}_i,\btheta^t_i)= \frac{p(\tau^t_i|\btheta^{t-1}_i)}{p(\tau^t_i|\btheta^t_i)} = \prod_{h=0}^{H-1}\frac{\pi_{\btheta^{t-1}_i}(\va^h|\vs^h)}{\pi_{\btheta^t_i}(\va^h|\vs^h)}.
\end{align}
Therefore, the momentum-based variance reduction~\eqref{eq: momentum} for the policy gradient of the $i$-th agent, denoted $\vv^{t}_{i}$, is given by~\cite{huang2020momentum}:
\begin{align*}
	\vv_i^t = \beta \vg_i(\tau_i^t|  \btheta_i^t) + (1-\beta)  \left( \vv_i^{t-1} + \vg_i(\tau_i^t | \btheta_i^t) - \omega(\tau_i^t | \btheta_i^{t-1}, \btheta_i^t)\cdot \vg_i(\tau_i^t | \btheta_i^{t-1})\right).
\end{align*}

\subsubsection{Pillar III: Natural policy gradient}


In order to introduce the natural policy gradient method, we consider the policy optimization problem \eqref{eq: rleq2}:
\begin{align*}
\max_{\btheta\in\R^d} \left\{ V(\btheta) = \E[\tau\sim p(\cdot | \btheta)]{R(\tau)} \right\},
\end{align*}
where $\tau$ is the $H$-horizon trajectory and $p(\tau|\btheta)$  given by \eqref{dist of tau}.  The PG update $\btheta^{t+1} = \btheta^t + \eta \nabla V(\btheta^t)$ is a gradient ascent method over the parameter space,  which is also the minimizer of the following problem 
\begin{align}
\label{opti 1}
    \min_{\btheta} \la -\nabla V(\btheta^t), \btheta - \btheta^t \ra + \frac{1}{2\eta}\twonorm{\btheta - \btheta^t}^2.
\end{align}
However, since the objective function essentially relies on the distributions of $\tau$, it is more natural to conduct a search over distribution space, leading to the following sub-problem for updating $\btheta^t$:
\begin{align}
\label{penalized KL}
    \min_{\btheta} \la -\nabla V(\btheta^t), \btheta - \btheta^t \ra + \frac{1}{2\eta}\KL( p(\tau|\btheta^t); p(\tau| \btheta) ),
\end{align}
where the KL divergence is used to enable the search  around $p(\tau|\btheta^t)$ over the distribution space.

Since $\KL( p(\tau|\btheta^t); p(\tau| \btheta^t) ) = 0$ and $\nabla_{\btheta^t} \KL( p(\tau|\btheta^t); p(\tau| \btheta) ) = \bzero$, one can approximate  $\KL( p(\tau|\btheta^t); p(\tau| \btheta) )$ by its second order information and thus approximate \eqref{penalized KL} by
\begin{align}
\label{npg optimization}
\min_{\btheta} \la -\nabla V(\btheta^t), \btheta - \btheta^t \ra + \frac{1}{2\eta }( \btheta - \btheta^t)^\tran \mF(\btheta^t) ( \btheta - \btheta^t),
\end{align}
where $\mF(\btheta^t)=\E[\tau\sim p(\cdot | \btheta^t) ]{\nabla_{\btheta} \log p(\tau| \btheta^t) \left(\nabla_{\btheta} \log p(\tau| \btheta^t)\right)^\tran}$ is the Fisher information matrix (FIM) of $p(\tau|\btheta^t)$ and  $\mF(\btheta^t)^\dagger$ is the Moore-Penrose pseudoinverse of $\mF(\btheta^t)$. 
It can be easily seen that the optimal solution to \eqref{npg optimization} is given by
\begin{align*}
    \btheta^{t+1} = \btheta^t + \eta \mF(\btheta^t)^\dagger \nabla V(\btheta^t),\numberthis\label{eq: npgupdate}
\end{align*}
which yields the natural policy gradient update.

Given the definition of $p(\tau|\btheta)$  in \eqref{dist of tau}, the FIM can be further expressed as
\begin{align*}
    \mF(\btheta) &= \E[\tau\sim p(\cdot | \btheta)]{\left(\sum_{h=0}^{H-1}\nabla_{\btheta} \log \pi_{\btheta}(\va^h|\vs^h)\right) \left(\sum_{h=0}^{H-1}\nabla_{\btheta} \log \pi_{\btheta}(\va^h|\vs^h)\right)^\tran}\\
    &= \E[\tau\sim p(\cdot | \btheta)]{\sum_{h=0}^{H-1}\nabla_{\btheta} \log \pi_{\btheta}(\va^h|\vs^h) \left(\nabla_{\btheta} \log \pi_{\btheta}(\va^h|\vs^h)\right)^\tran}. \numberthis \label{fim 2}
\end{align*}
Here the second line has used the fact that the cross term is equal to $0$, which can be easily verified. When $H\rightarrow \infty$, the FIM may not be well determined. There are two typical ways to deal with this issue:
\begin{itemize}
    \item Averaged case. Let $\tau$ be a trajectory induced by $\pi_{\btheta}$ up to horizon $H$. The FIM in the average case is given by
    \begin{align}
    \label{eq limit}
        \mF(\btheta) 
        &= \lim_{H\rightarrow \infty} \frac{1}{H}\E[\tau\sim p(\cdot | \btheta)]{\sum_{h=0}^{H-1} \nabla_{\btheta} \log \pi_{\btheta}(\va^h|\vs^h) \left( \nabla_{\btheta} \log \pi_{\btheta}(\va^h|\vs^h)\right)^\tran},
    \end{align}
    It has been shown by \cite{bagnell2003covariant, peters2003reinforcement} that \eqref{eq limit} is equivalent to
    \begin{align*}
        \mF(\btheta) = \E[\vs\sim d^{\pi_\theta}, \va\sim \pi_{\btheta}(\cdot | \vs)]{\nabla_{\btheta} \log \pi_{\btheta}(\va|\vs) \left( \nabla_{\btheta} \log \pi_{\btheta}(\va|\vs)\right)^\tran},
    \end{align*}
    where $d^{\pi_{\btheta}}$ is the stationary distribution of state.
    
    \item Discounted case. On the other hand, one can consider infinite horizon but introduce a discounted factor $\gamma\in [0, 1)$. In this situation, the FIM is given by
    \begin{align}
    \label{eq discounted FIM}
        \mF(\btheta) = \E[\tau\sim p(\cdot | \btheta)]{ \sum_{h=0}^{+\infty}\gamma^h \nabla_{\btheta} \log \pi_{\btheta}(\va^h|\vs^h) \left( \nabla_{\btheta} \log \pi_{\btheta}(\va^h|\vs^h)\right)^\tran }.
    \end{align}
    Moreover, letting $\tau$ be the $H$-horizon trajectory induced by $\pi_{\btheta}$, where $H$ obeys the geometric distribution with parameter $1-\gamma$, then $F(\btheta)$ in \eqref{eq discounted FIM} is indeed the FIM associated with the random-length trajectory $\tau$. That is \cite{bagnell2003covariant, peters2003reinforcement},
    \begin{align*}
        \mF(\btheta) &=\E[H\sim \Geo(1-\gamma)]{\E[\tau\sim p(\cdot |\btheta)]{\sum_{h=0}^{H-1} \nabla_{\btheta} \log \pi_{\btheta}(\va^h|\vs^h) \left( \nabla_{\btheta} \log \pi_{\btheta}(\va^h|\vs^h)\right)^\tran \Bigg| H}}\notag\\
        &=\E[\vs\sim d_{\rho}^{\pi_{\btheta}}, \va\sim \pi_{\btheta}(\cdot | \vs)]{\nabla_{\btheta} \log \pi_{\btheta}(\va|\vs) \left( \nabla_{\btheta} \log \pi_{\btheta}(\va|\vs)\right)^\tran },
    \end{align*}
    where $d_{\rho}^{\pi_{\btheta}}(\vs)= \E[\vs^0\sim \rho]{(1-\gamma)\sum_{h=0}^{\infty} \gamma^h P(\vs^h = \vs | \vs^0,\pi_{\btheta}) }$ is the discounted state visitation distribution under the initial  distribution $\rho$. Such formulation has been widely used in the literature  \cite{kakade2001natural, bhatnagar2007incremental,agarwal2021theory}. 
\end{itemize}
With a slight abuse of notation, we will use the following definition of FIM in this paper:
\begin{align}
\label{def fim}
    \mF(\btheta) = \E[\tau\sim p(\cdot | \btheta)]{\frac{1}{H}\sum_{h=0}^{H-1}\nabla_{\btheta} \log \pi_{\btheta}(\va^h|\vs^h) \left(\nabla_{\btheta} \log \pi_{\btheta}(\va^h|\vs^h)\right)^\tran},
\end{align}
which agrees with \eqref{fim 2} up to a scale.


In contrast to PG, NPG  can be approximately viewed as a second-order method since the FIM serves as a structured preconditioner based on the underlying structure of the parameterized policy space \cite{amari1996neural, martens2020new}. Such a preconditioner can adaptively adjust the update direction to improve the convergence rate. In Algorithm~\ref{alg: moentum pg}, we have extended the NPG update~\eqref{eq: npgupdate} to the decentralized multi-agent setting~\eqref{eq: marl3}. For the $i$-th agent at the $t$-th iteration, we have
\begin{align*}
	\btheta_i^{t+1} = \sum_{j\in\calN(i)}W_{ij} \left( \btheta_j^t + \eta \mH_j^t \nabla V_j(\btheta^t_j) \right),
\end{align*}
where $\mH_j^t\in\mathbb{R}^{d\times d}$ denotes the Moore-Penrose pseudoinverse of $\mF_j(\btheta_j^t)$.
Namely, each agent searches along the (preconditioned) natural gradient direction of its own before the consensus. 

In addition, the following lemma establishes that the FIM in the collaborative RL setting mentioned in Section~\ref{section: marl} is indeed a block diagonal matrix for each agent due to the product structure of the joint policy, see \eqref{eq: pi1}.
\begin{lemma}
    \label{lemma: FIM}
    In collaborative RL, let $\tau_i$ be the $H$-horizon trajectory induced by the policy $\pi_{\btheta_i}$ for agent $i$. The FIM of the $i$-th agent  $\mF_i(\btheta_i)\in\mathbb{R}^{d\times d}$ is given by
\begin{align*}
\mF_i(\btheta_i) =  \diag\left (\mF_i(\btheta_{[1]}),\cdots, \mF_i(\btheta_{[n]}) \right ),
\end{align*}
where 
\begin{align*}\numberthis\label{eq: lemmafim}
    \mF_i(\btheta_{[j]})=  \frac{1}{H} E_{\tau_i\sim p(\cdot|\btheta)} \left\{\sum_{h=0}^{H-1}  \nabla_{\btheta_{[j]}} \log \pi_{\btheta_{[j]}}(\va^h_j | \vs^h) \left( \nabla_{\btheta_{[j]}} \log \pi_{\btheta_{[j]}}(\va^h_{j} | \vs^h)\right)^T \right\} \in \mathbb{R}^{d_j\times d_j},
\end{align*}
for $j=1,\cdots,n$.
\end{lemma}

\subsection{Theoretical result}
Before stating the main convergence result in Theorem~\ref{main result}, we first introduce some standard assumptions.

\begin{assumption}
\label{assumption 6} The weight matrix $\mW\in\R^{n\times n}$ associated with the communication graph $\calG$  is doubly stochastic, i.e., $\mW\bone_n = \bone_n$ and $\bone^\tran_n \mW = \bone_n^\tran$.
\end{assumption}
Under Assumption \ref{assumption 6}, one can show that \cite{qu2017harnessing} 
\begin{align}
\label{def rho}
    \rho := \opnorm{\mW - \frac{1}{n}\bone_n \bone_n^\tran} \in[0,1).
\end{align}


\begin{assumption}
\label{assumption 0}
	The objective function  $V(\btheta)$ is upper bounded, i.e., $V^\star := \sup_{\btheta\in\R^d} V(\btheta)<+\infty$.
\end{assumption}

\begin{assumption}
\label{assumption 1}
	Let $\pi_{\btheta}(\va | \vs)$ be the policy parameterized by $\btheta \in\R^{d}$. There are constants $G$ and $M$ such that the gradient and Hessian of the log-density of the policy function obey that 
	\begin{align*}
		\twonorm{\nabla_{\btheta} \log \pi_{\btheta}(\va | \vs)}^2 \leq G \text{ and } \opnorm{\nabla^2_{\btheta} \log \pi_{\btheta}(\va | \vs)} \leq M
	\end{align*}
for any $\va\in\calA$ and $\vs\in\calS$.
\end{assumption}

\begin{assumption}
\label{assumption 5}
  The variance of $\omega(\tau| \widetilde{\btheta}, \btheta)$, the importance sampling weight defined in \eqref{eq important sampling}, is bounded, 
	\begin{align*}
		\Var(\omega(\tau| \widetilde{\btheta}, \btheta) ) \leq W,
	\end{align*}
for any $\btheta, \widetilde{\btheta}\in\R^d$ and $\tau\sim p_i(\cdot| \btheta )$.
\end{assumption}

The following two auxiliary lemmas can be obtained from  Assumptions \ref{assumption 1} and \ref{assumption 5}.
\begin{lemma}[Proposition 5.2 in \cite{xu2020improved}]
	\label{lemma 1}
	Under Assumption \ref{assumption 1}, one has the following facts:
	\begin{itemize}
		\item The objective function $V(\btheta)$ is $L$-smooth with $L = H(M+HG)/(1-\gamma)$. 
		\item Let $\vg_i(\tau; \btheta)$ be the gradient estimator defined in \eqref{eq: grad2}. Then for all $\btheta, \widetilde{\btheta}\in\R^d$, one has
		\begin{align*}
			\twonorm{\vg_i(\tau; \btheta) - \vg_i(\tau; \widetilde{\btheta})} \leq L_g \twonorm{\btheta - \widetilde{\btheta}}
		\end{align*}
		and $\twonorm{\vg_i(\tau; \btheta)} \leq C_g$ for all $i\in[n]$, where $L_g = HM(1+|b| )/(1-\gamma)$ and $C_g = HG^{1/2}(1+|b|)/(1-\gamma)$. 
	\end{itemize}
\end{lemma}
\begin{lemma}[Lemma 3 in \cite{jiang2021mdpgt}]
Under Assumption \ref{assumption 1} and \ref{assumption 5}, one has 
\begin{align}
	\label{variance of importance weight}
	\Var(\omega(\tau| \widetilde{\btheta}, \btheta) )  =\E[\tau\sim p(\cdot | \btheta)]{\left(\omega(\tau| \widetilde{\btheta}, \btheta) - 1\right)^2}= C_{\omega}^2 \twonorm{\widetilde{\btheta} - \btheta}^2,
\end{align}
where $C_{\omega }^2 = H(2HG + M)(W+1)$.
\end{lemma}

The next assumption states that the variance of the gradient estimator is bounded, which is commonly used in stochastic optimization. 
\begin{assumption}
\label{assumption 3}
	There exists a constant $\nu_i$ such that 
	\begin{align}
		\label{boundedVariance}
		\Var( g_i(\tau; \btheta) ) = \E{\twonorm{\vg_i(\tau; \btheta) - \nabla V_i(\btheta)}^2} &\leq \nu_i^2
	\end{align}
for all policy $\pi_{\btheta}$, where $\tau\sim p(\cdot | \btheta)$. Define $\bar{\nu}^2 = \frac{1}{n}\sum_{i=1}^{n}\nu_i^2$.
\end{assumption}
\begin{assumption}
\label{assumption 4}
	Let $\mF(\btheta)\in\R^{d\times d}$ be the FIM defined in \eqref{def fim}.
There exists a constant $\mu_F >0$ such that $\mF(\btheta) \succeq \mu_F \mI_d$ for all $\btheta \in\R^d$.
\end{assumption}
\begin{remark}
The positive definiteness on FIM is fairly standard in the analysis of single-agent NPG algorithm \cite{liu2020improved} and is often made in both convex and nonconvex optimizations for establishing the convergence of preconditioned algorithms \cite{broyden1970convergence, fletcher1970new, goldfarb1970family, shanno1970conditioning}. Moreover, this assumption can be always satisfied if we use $\mF+\epsilon \mI_d$ for $\epsilon >0$ instead of $\mF$ as the precondition matrix.  
\end{remark}

Let $\mH_i$ be the inverse FIM $\mF(\btheta_i)$ for the $i$-th agent with policy parameter $\btheta_i\in\R^d$. Then Assumptions \ref{assumption 1} and \ref{assumption 4}  imply that 
\begin{align}
\label{eig of H}
\frac{1}{G }  \mI_d \preccurlyeq \mH_i^t  \preccurlyeq \frac{1}{\mu_F }\mI_d,
\end{align}
where the lower bound holds since the fact $\opnorm{\mF(\btheta)} \leq G$. Moreover, this fact can be proved as follows: 
\begin{align*}
	\opnorm{\mF(\btheta)} &=\opnorm{\E[\tau\sim p(\cdot | \btheta)]{\frac{1}{H}\sum_{h=0}^{H-1}\nabla_{\btheta} \log \pi_{\btheta}(\va^h|\vs^h) \left(\nabla_{\btheta} \log \pi_{\btheta}(\va^h|\vs^h)\right)^\tran}}\\
	&\leq \frac{1}{H}\sum_{h=0}^{H-1} \E[\tau\sim p(\cdot | \btheta)]{\opnorm{\nabla_{\btheta} \log \pi_{\btheta}(\va^h|\vs^h) \left(\nabla_{\btheta} \log \pi_{\btheta}(\va^h|\vs^h)\right)^\tran}}\\
	&\leq G,
\end{align*}
where the last line is due to Assumption \ref{assumption 1}.






We are in position to present the main result of this paper.
\begin{theorem}
\label{main result}
Let $\btheta_{\out}\in\R^d$ be the output of Algorithm \ref{alg: moentum pg}. Suppose that  
\begin{align*}
	 0 < \eta  <\frac{\mu_F(1-\rho^2)^3}{ \kappa_{F}\sqrt{1632000 (L^2 + \Phi^2) }}  
\end{align*}
and choose $\beta$ such that $\frac{1632000(L^2 + \Phi^2) \kappa_F^2 \eta^2 }{n\mu_F^2(1-\rho^2)^6} \leq \beta <\frac{1}{n}$, where $\Phi^2 = L_g^2 + C_g^2 C_{\omega}^2$ and $\kappa_F = G/\mu_F$.
Then under Assumptions \ref{assumption 6}-\ref{assumption 4},  one has
\begin{align}
	\E{\twonorm{\nabla V(\btheta_{\out}) }^2} &\leq \frac{8G^2 \Delta}{T\eta\mu_F} + \frac{76\bar{\nu}^2 \kappa_F^2}{nT\beta B }   + \frac{152\beta \bar{\nu}^2 \kappa_F^2 }{n}   \notag\\
	&\qquad +     \frac{44\rho^2\bar{\nu}^2 \kappa_F^2}{TB(1-\rho^2)}  + \frac{352\beta^2 \bar{\nu}^2\kappa_F^2 }{(1-\rho^2)^2}  + \frac{352\beta \bar{\nu}^2\kappa_F^2 }{TB(1-\rho^2)^2} +  \frac{704  \bar{\nu}^2 \kappa_F^2 \beta^3}{(1-\rho^2)^2}  + \frac{44\rho^2\kappa_F^2}{nT(1-\rho^2)} \twonorm{\tnablaV(\btheta^0)}^2.
\end{align}
\end{theorem}

\begin{remark}
If we choose $\eta$ and $\beta$ according to Theorem \ref{main result}, then the mean squared stationary gap $\E{\twonorm{\nabla V(\btheta_{\out}) }^2} $ converges to a steady-state error as $T\rightarrow \infty$ at a rate of $\calO(1/T)$, i.e., 
\begin{align*}
\E{\twonorm{\nabla V(\btheta_{\out}) }^2} \rightarrow \calO\left(\frac{\beta \bar{\nu}^2}{n} + \frac{\beta^2 \bar{\nu}^2}{(1-\rho^2)^2} \right),\quad T\rightarrow \infty.
\end{align*}
It can be seen that the steady-state error will decrease as the number of agents increases. Moreover, the second term in the steady-state error indicates that the impact of communication graph $\calG$ through $\rho$ can be reduced with small $\beta$.
\end{remark}

\begin{corollary}
	\label{corollary}
	Choose step size $\eta$, momentum parameter $\beta$, and batch size $B$ in initialization such that 
	\begin{align*}
		\eta = \frac{\mu_F n^{2/3}}{\kappa_{F} \sqrt{L^2 + \Phi^2}T^{1/3}}, \beta = \frac{n^{1/3}}{T^{2/3}} \text{ and } B = \left\lceil \frac{T^{1/3}}{n^{2/3}} \right\rceil.
	\end{align*}
	Then for all $T>  \frac{1632000^{3/2}n^2}{(1-\rho^2)^9} $, one has 
	\begin{align*}
		\E{\twonorm{\nabla V(\btheta_{\out})}^2}  & \leq  \frac{8\Delta   \kappa_{F}^3\sqrt{L^2 + \Phi}  + 228\bar{\nu}^2\kappa_F^2}{(nT)^{2/3}} \\
		&\qquad +\frac{44\rho^2\kappa_{F}^2 \twonorm{\tnablaV(\btheta^0)}^2}{ (1-\rho^2)} \cdot \frac{1}{nT}+ \frac{396\kappa_{F}^2\bar{\nu}^2}{ (1-\rho^2)^2}\cdot \frac{n^{2/3}}{T^{4/3}}+ \frac{1056 \bar{\nu}^2  \kappa_{F}^2 }{ (1-\rho^2)^2} \cdot \frac{n}{T^2}.
	\end{align*}

\end{corollary}
\begin{remark}
Corollary \ref{corollary} implies that 
\begin{align*} 
    \E{\twonorm{\nabla V(\btheta_{\out})}^2}  = \calO((nT)^{-2/3})
\end{align*}
when $T$ is large enough. Thus one can achieve $\varepsilon$-stationary, i.e., $\E{\twonorm{\nabla V(\btheta_{\out})}^2} \lesssim \varepsilon^2$, in $\calO( n^{-1}\varepsilon^{-3} )$ iteration complexity, which shows that MDNPG also enjoys the linear speedup convergence rate. 
\end{remark}

\subsection{Related work}
In this section, we discuss the recent progress that is mostly related to our work, especially those gradient based methods in reinforcement learning and decentralized stochastic optimization. 
\paragraph{Single-agent (natural) policy gradient.} Inspired by stochastic optimization, there has been extensive research in designing variance reduction methods for policy gradient estimator \cite{papini2018stochastic,xu2019sample,xu2020improved,shen2019hessian,huang2020momentum}. 
For instance, Papini et al. \cite{papini2018stochastic} show that SVRPG achieves an $\epsilon$-stationary point given $\calO(\epsilon^{-4})$ trajectories. Xu et al. \cite{xu2020improved} improve this sample complexity to $\calO(\epsilon^{-10/3})$. Moreover, SRVR-PG \cite{xu2019sample} and HAPG \cite{shen2019hessian} can obtain an $\epsilon$-stationary point provided $\calO(\epsilon^{-3})$ trajectories, both of which are nearly optimal in the sample complexity \cite{arjevani2022lower}. 
However, these methods require large batches or double-loop updates. Recently, \cite{huang2020momentum} incorporates the momentum-based variance reduction technique in policy gradient methods and achieve the sample complexity of $\calO(\epsilon^{-3})$ with a single trajectory at each iteration. The global convergence of policy gradients with variance reduction has also been studied in \cite{ding2021global,liu2020improved}. 

As already mentioned, NPG  \cite{kakade2001natural} and its generalizations, such as TRPO \cite{schulman2015trust} and PPO \cite{schulman2017proximal}, are widely used in RL. There has been a lot of interest in understanding the theoretical performance of this  class of methods, see \cite{agarwal2021theory,cen2021fast,lan2022policy} and the references therein. Different variance reduction techniques have also been utilized in NPG. For example, SRVR-NPG is proposed in \cite{liu2020improved} which reaches a sample complexity of $\calO(\epsilon^{-3})$. In addition, two variance reduced mirror ascent methods, named VRMPO and VR-BGPO, are developed in \cite{yang2022policy} and  \cite{huang2021bregman}, 
respectively.  These methods reduces to NPG with variance reduction if a special mirror mapping is used.

\paragraph{Multi-agent policy gradient.}  For collaborative RL problem, many decentralized policy gradient algorithms have been developed. Lu et al. \cite{lu2021decentralized} 
study a decentralized policy gradient method in safe MARL and show that an $\epsilon$-stationary point can be achieved from $\calO(\epsilon^{-4})$ iterations. Zhao et al. \cite{zhao2021distributed} 
study the convergence of decentralized policy gradient with variance reduction and gradient tracking in collaborative RL and establish the sample complexity of $\calO(\epsilon^{-3})$. However, the  method in  \cite{zhao2021distributed} requires very large batch gradients to obtain this optimal complexity. In contrast, Jiang et al. \cite{jiang2021mdpgt} 
adopt the momentum-based variance reduction technique for decentralized policy gradient which also achieves the optimal sample complexity but only uses a single trajectory in each iteration. For the MTRL problem, various policy gradient methods have been developed and studied. In \cite{espeholt2018impala, hessel2019multi}, 
a distributed framework is used to solve the learning problem. However, in these works, each agent collects local data, which are then shared to a centralized coordination. In a subsequent work, a decentralized policy gradient method is proposed in \cite{zeng2021decentralized}. However, the proposed method only adopts the vanilla gradient ascent without gradient tracking and variance reduction, thus resulting in a sample complexity of $\calO(\epsilon^{-4})$. In addition, decentralized optimization methods have also been studied in the framework of policy evaluation \cite{qu2019value, doan2019finite, lin2021decentralized}. 

\paragraph{Decentralized optimization.} In general, decentralized online optimization has been extensively studied for non-convex problems. There are many algorithms developed toward this line of research, including DSGD \cite{lian2017can}, EXTRA \cite{shi2015extra}, and Exact Diffusion \cite{yuan2018exact}. To achieve the lower oracle complexity, various kinds of variance reduced techniques have been utilized. The D-GET proposed in \cite{sun2020improving} is built upon gradient tracking and the SARAH gradient estimator and achieves an oracle complexity of $\calO(\epsilon^{-3})$. The same oracle complexity is also obtained by D-SPIDER-SFO \cite{pan2020d} which uses SPIDER in the variance reduction step. Recently, built on the hybrid variance reduction scheme introduced in \cite{cutkosky2019momentum, tran2022hybrid}, the GT-HSGD is developed in \cite{xin2021hybrid} which can achieve an  $\epsilon$-approximate ﬁrst-order stationary point within $\calO(n^{-1}\epsilon^{-3})$ samples for each node. It is worth noting that these approaches are designed for oblivious objective functions, where the randomness of data samples is independent of the optimized parameters. 

\section{Numerical experiments} \label{section experiments}
In this section, we empirically compare MDNPG with other state-of-the-art algorithms in several typical RL environments\footnote{Codes for reproducing the computational results in this section are available at \url{https://github.com/fccc0417/mdnpg}.}.
In our implementations, we use the sample version of \eqref{def fim} with one random trajectory to approximately compute the FIM for the single-agent experiments as well as the experiments about multi-task GridWorld:
\begin{align*}
    \mF(\btheta) \approx  \frac{1}{H}\sum_{h=0}^{H-1}\nabla_{\btheta} \log \pi_{\btheta}(\va^h|\vs^h) \left(\nabla_{\btheta} \log \pi_{\btheta}(\va^h|\vs^h)\right)^\tran.
\end{align*}
 Regrading the experiments for the collaborative RL setting on cooperative navigation, the FIM $\mF_i(\btheta_{[j]})$ is computed via  \eqref{eq: lemmafim} of Lemma~\ref{lemma: FIM}, 
 \begin{align*}
    	\mF_i(\btheta_{[j]}) \approx  \frac{1}{H} \sum_{h=0}^{H-1} \nabla_{\btheta_{[j]}} \log \pi_{\btheta_{[j]}}(\va^h_j | \vs^h) \left( \nabla_{\btheta_{[j]}} \log \pi_{\btheta_{[j]}}(\va^h_{j} | \vs^h)\right)^T.
 \end{align*}

\subsection{Single-agent experiments}\label{singleagent}
When $n=1$, MDNPG reduces to the single-agent NPG with momentum-based variance reduction, which can also be viewed as BGPO~\cite{huang2021bregman} with special mirror mappings. In this subsection, we compare the single-agent version of MDNPG with the momentum-based policy gradient~\cite{huang2020momentum}, PPO~\cite{schulman2017proximal}, and SRVR-NPG {(NPG with variance reduction via SRVR)}~\cite{liu2020improved} over two single-agent environments: GridWorld and MountainCar.

In a $10\times 10$ GridWorld, an agent at a random initial position try to reach the grid labeled as ``goal'' and at the same time avoid grids labeled as ``obstacle''  in a minimum number of steps. Five obstacle grids are set up in our experiments. The agent can select one of four discrete actions (up, down, right, and left) to move to another grid, and the state is simply the location of the agent. The received reward is $-0.1\times\mbox{(distance to the goal)} \pm 10$, up  to whether the goal is reached or the agent falls into an obstacle. The other environment, called MountainCar, is a continuous control task from OpenAI Gym~\cite{brockman2016openai}, in which the goal of an agent  is to reach the top of the hill. A detailed description of the environment is provided in \cite{brockman2016openai}.
 
In our implementation, a one-hidden-layer (of size 128) neural network with ReLU activation function is used to parameterize the policy. For the GirdWorld task, the parameterized policy can be obtained via a softmax layer. For the MountainCar task, the outputs of the network are $(\mu_{\btheta}(s), \sigma_{\btheta}(s))$,  the mean  and standard deviation  of a Gaussian distribution. Moreover, we also use a value network (one hidden layer of size 128, with ReLU as the activation function) for the estimation of value functions. All the parameters in the algorithms are finely tuned for the pursuit of better performance. 

\begin{figure}[ht!]
 \centering
\begin{subfigure}[b]{0.49\textwidth}
  \centering
  \includegraphics[width=0.95\linewidth]{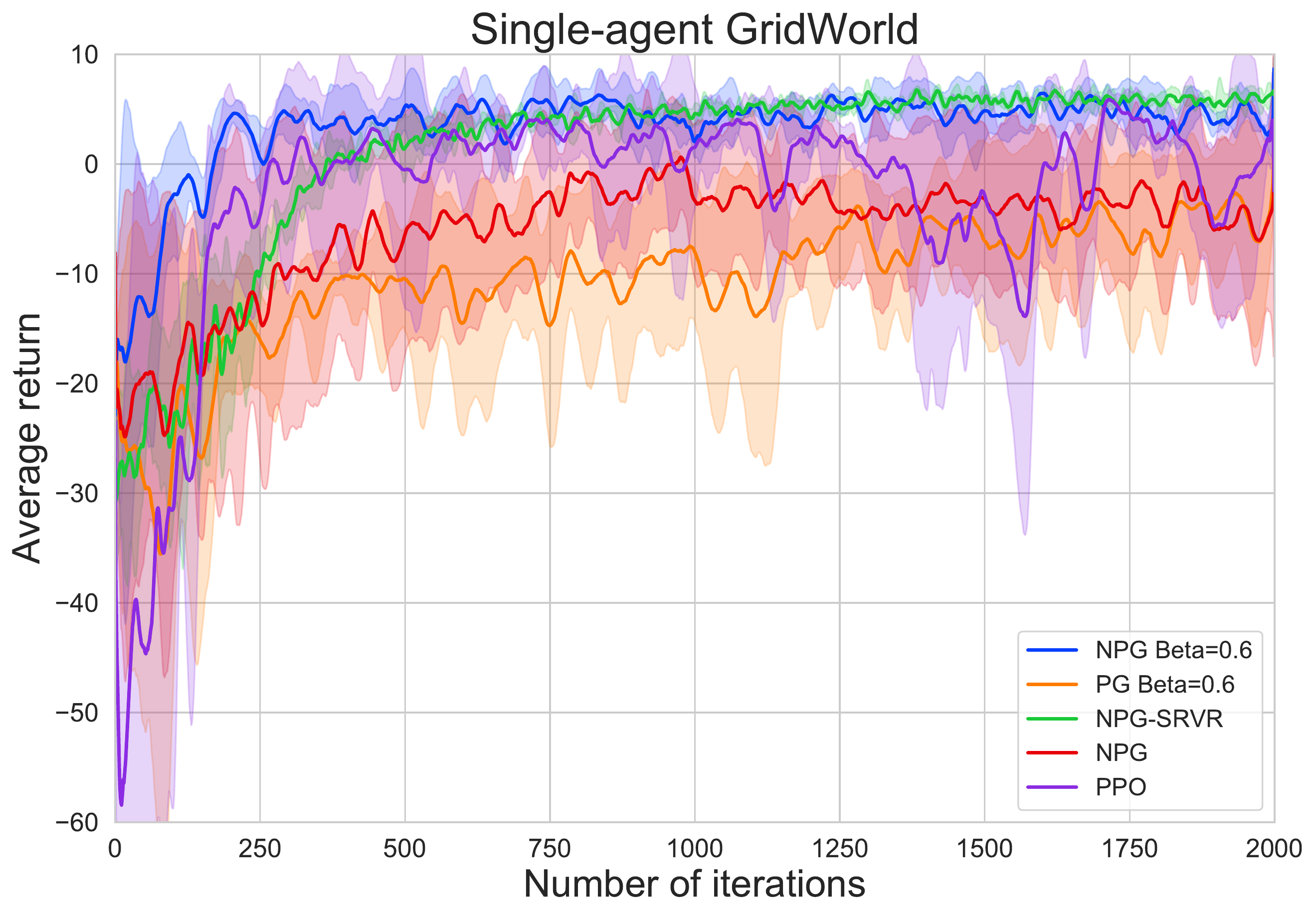}
  \caption{}
  \label{fig: gridworld}
\end{subfigure}
\begin{subfigure}[b]{0.49\textwidth}
  \centering
  \includegraphics[width=0.95\linewidth]{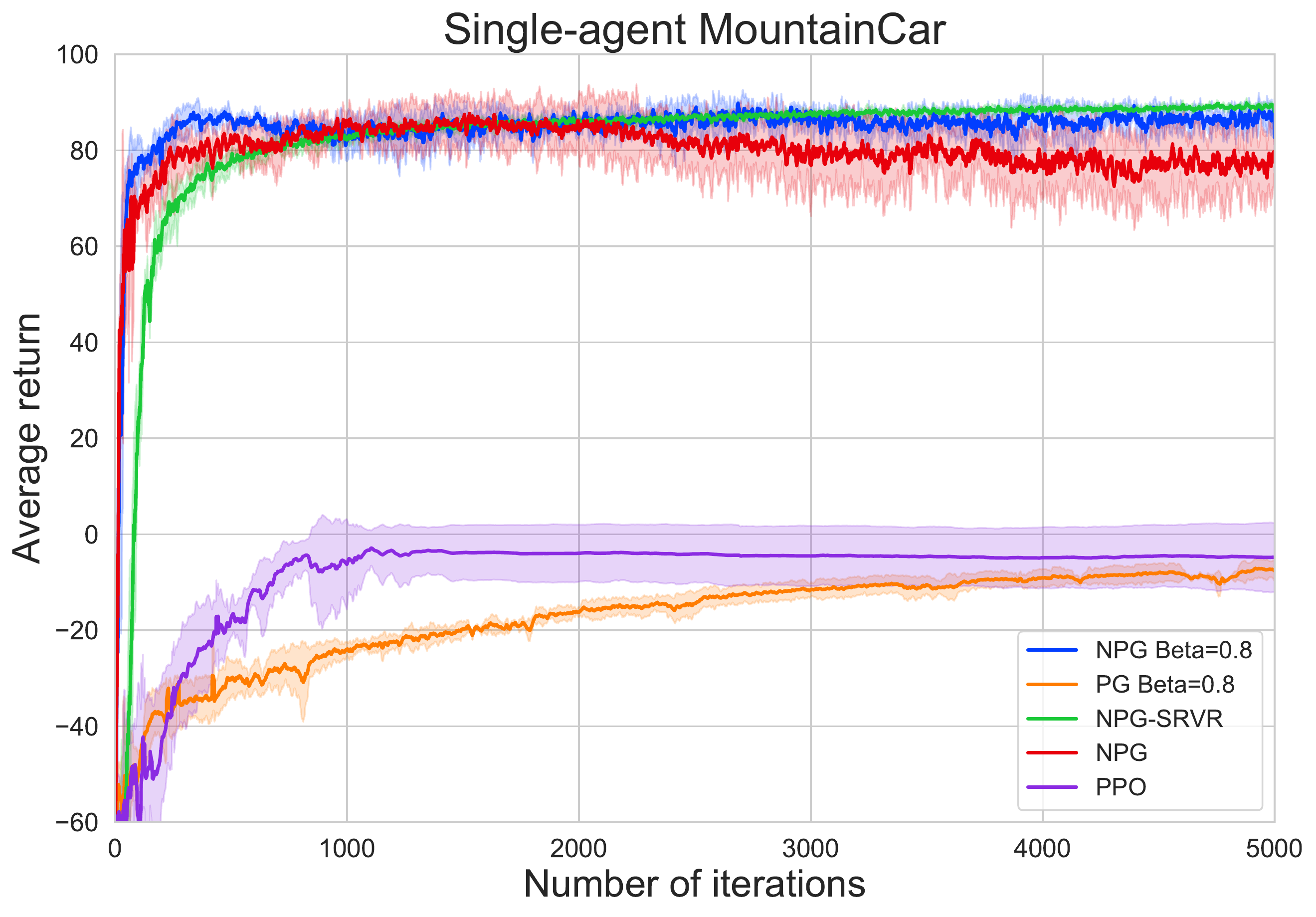}
  \caption{}
  \label{fig: mountaincar}
\end{subfigure}
\caption{Average return and standard deviation over five random instances against the number of iterations. The parameter $\beta$ in GridWorld (a)  and MountainCar (b) is set to $0.6$ and $0.8$ for the momentum-based NPG and PG methods. }
\label{fig: singleagent1}
\end{figure}

The plots of average return and standard deviation over five random instances against the number of iterations are presented in Figure~\ref{fig: singleagent1}. It can be observed that overall the momentum-based NPG method displays better convergence and stability than the other algorithms for both tasks. Note that even though  SRVR-NPG is competitive with  the momentum-based NPG method, the former one costs significantly more time and memory due to its double-loop nature for variance reduction.

\subsection{Multi-agent experiments}\label{sec: multiagent}

\subsubsection{Cooperative navigation}For the collaborative RL setting in Section~\ref{section: marl}, we compare  MDNPG with other state-of-the-art algorithms such as MDPGT~\cite{jiang2021mdpgt} and value propagation~\cite{qu2019value} on a simulated cooperative navigation environment introduced by \cite{lowe2017multi}.  As a benchmark multi-agent environment, it has been modified in several previous works such as \cite{zhang2018fully,qu2019value, jiang2021mdpgt} to be compatible with the collaborative RL setting. In the $n$-agent cooperative navigation, each agent at a randomly initialized location needs to find its specific landmark and avoid collisions with other agents in a rectangle region of size $2\times2$. Agents can move up, down, right, left, or keep still at each step. The globally observed state consists of the positions of all agents as well as their landmarks. The received reward of each agents is  $-\text{(distance to the landmark)} - \sum \one_{\{\text{if colliding with an agent}\}}$. 

More precisely, there are $5$ agents in  our experiments. The policy-based methods, MDNPG as well as MDPGT, utilize a policy network and a value network, both of which have two hidden layers with 64 and 128 units and use ReLU as the activation function. Additionally, the value propagation method utilizes another auxiliary network  to approximate the dual function. Since MDNPG and MDPGT are both on-policy algorithms, the on-policy  version of value propagation is implemented here. 


\begin{figure}[ht!]
  \centering
\begin{subfigure}[b]{.27\textwidth}
  \centering
  \includegraphics[width=.9\linewidth]{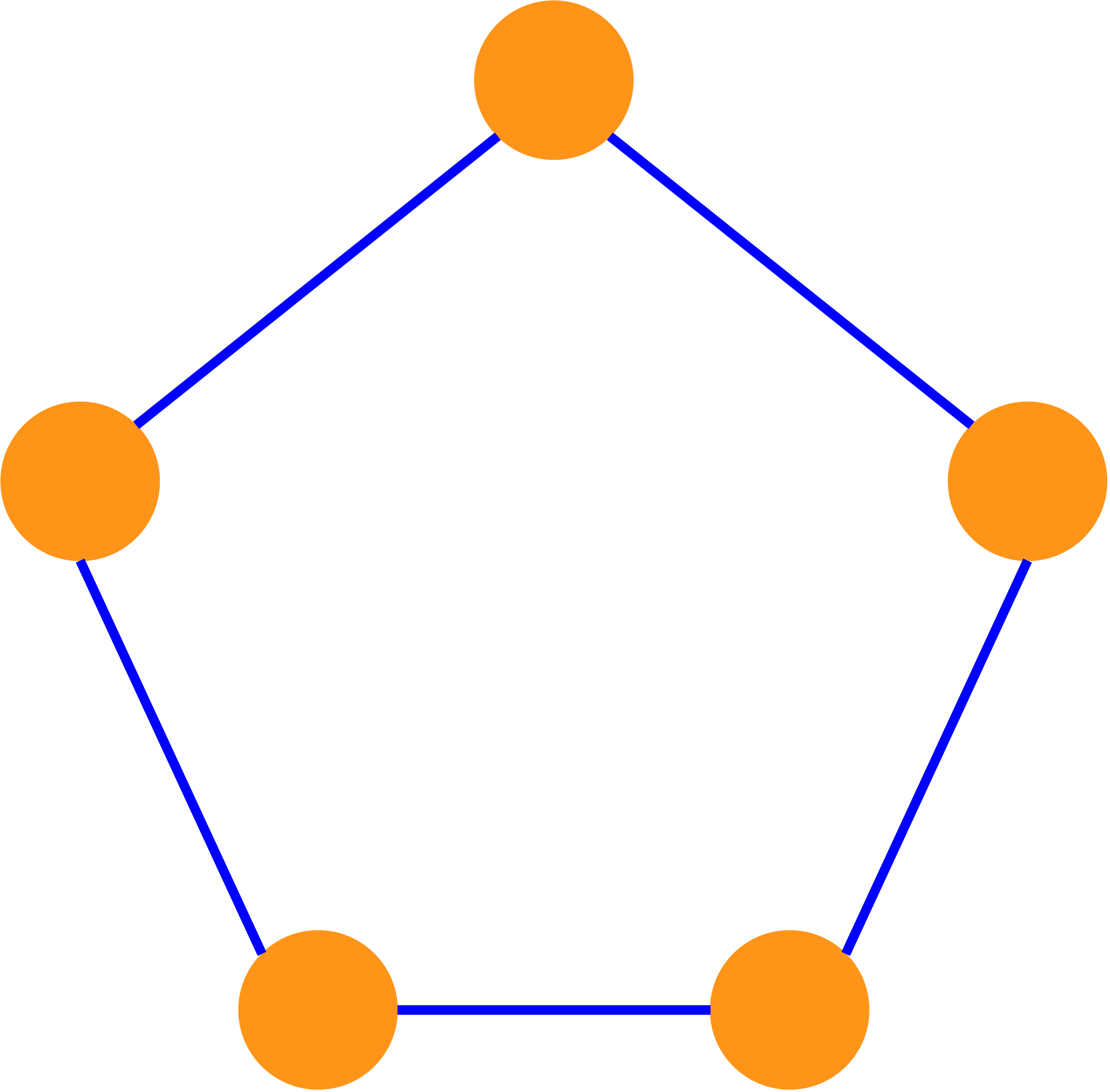}
  \caption{Ring}
\end{subfigure}
\begin{subfigure}[b]{.27\textwidth}
  \centering
  \includegraphics[width=.9\linewidth]{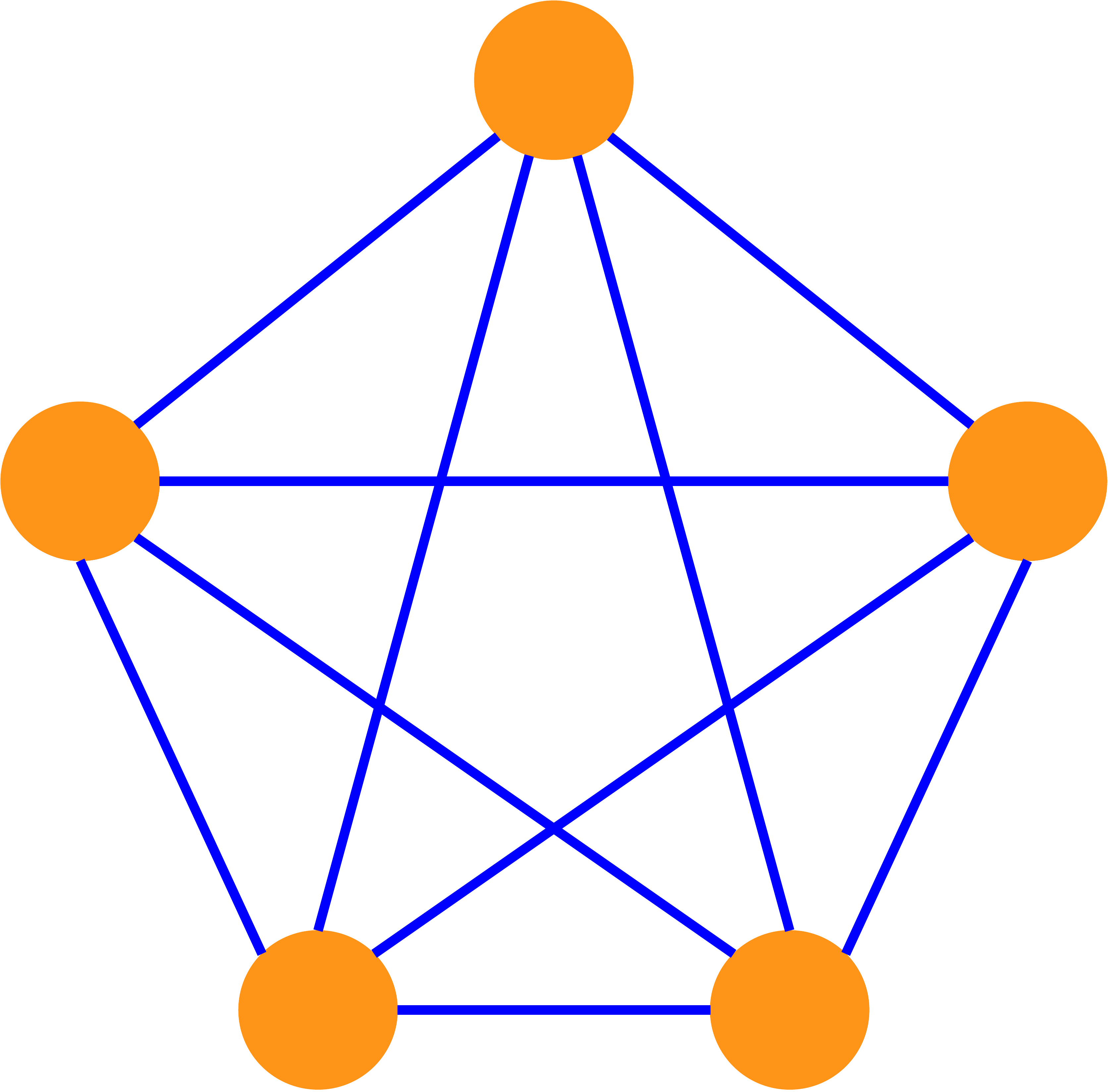}
  \caption{Fully-connected }
\end{subfigure}
\begin{subfigure}[b]{.27\textwidth}
  \centering
  \includegraphics[width=.9\linewidth]{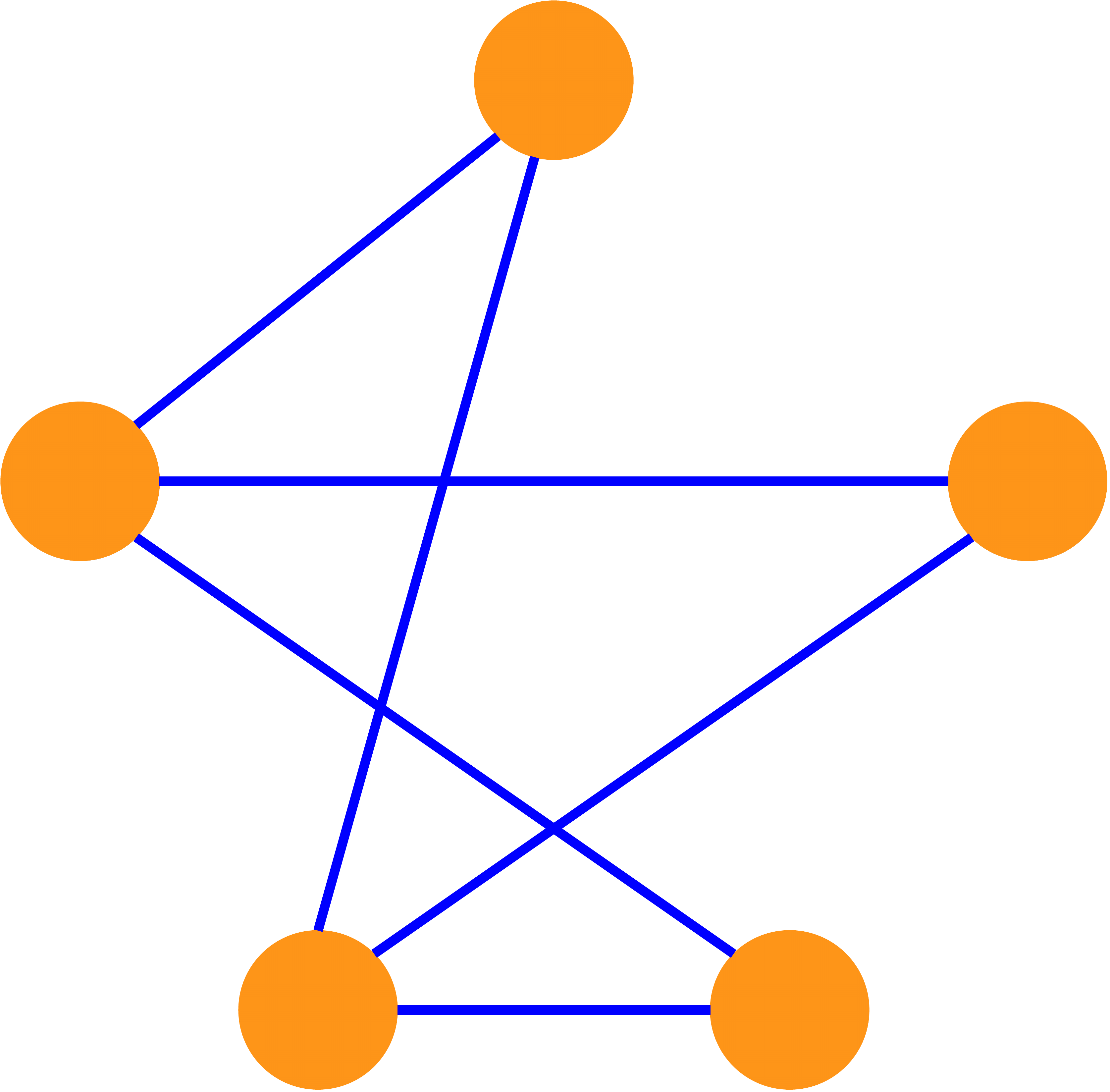}
  \caption{Bi-partite}
\end{subfigure}
\caption{Three network topologies. }
\label{fig: topology}
\end{figure}

\begin{figure}[ht!]
\centering
    \begin{subfigure}[]{0.49\textwidth}
    \centering
         \includegraphics[width=.95\linewidth]{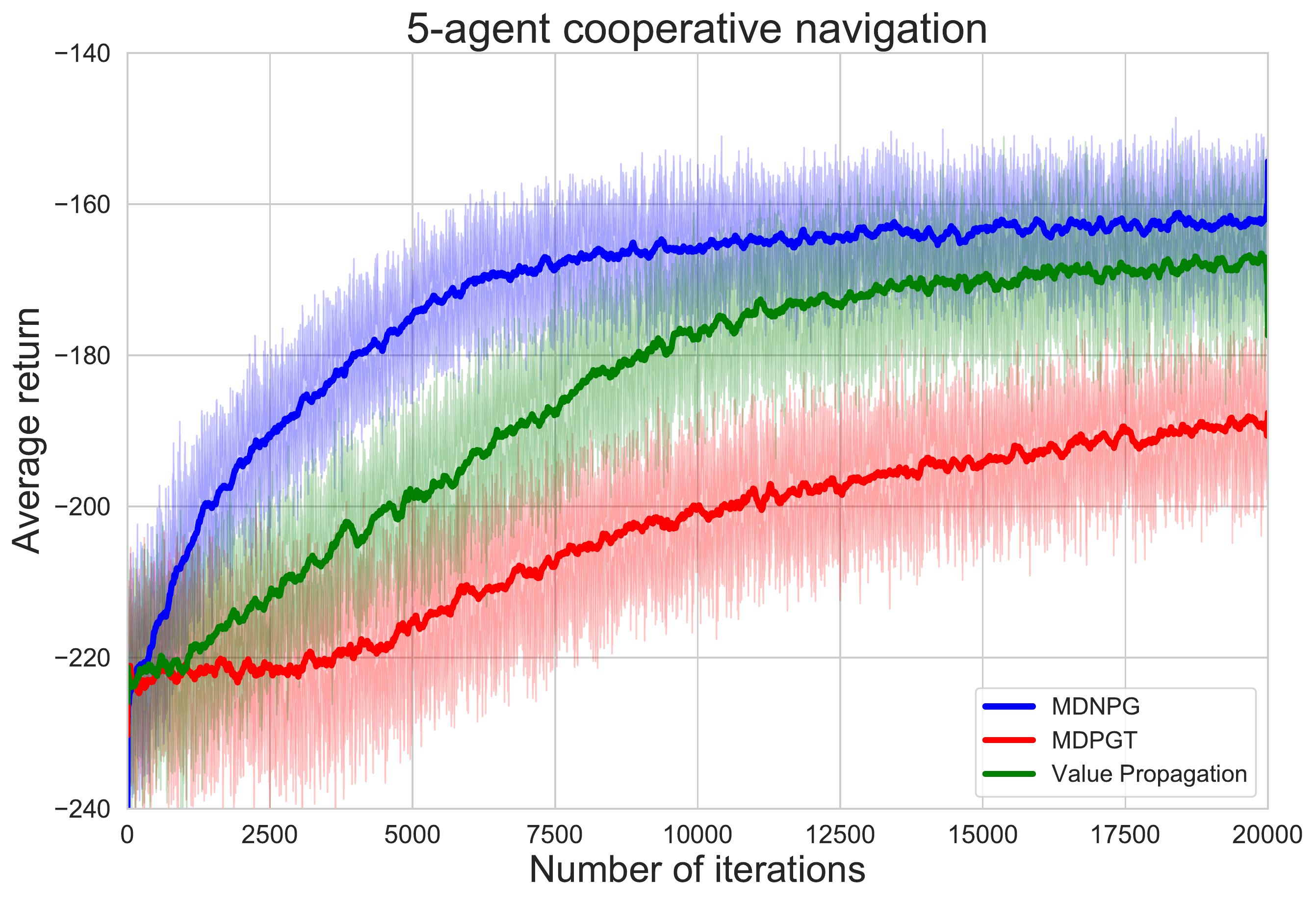}
         \caption[]{Ring} 
         \label{fig: navigation_ring}
     \end{subfigure}
     \hfill
     \begin{subfigure}[]{0.49\textwidth}
     \centering
         \includegraphics[width=.95\linewidth]{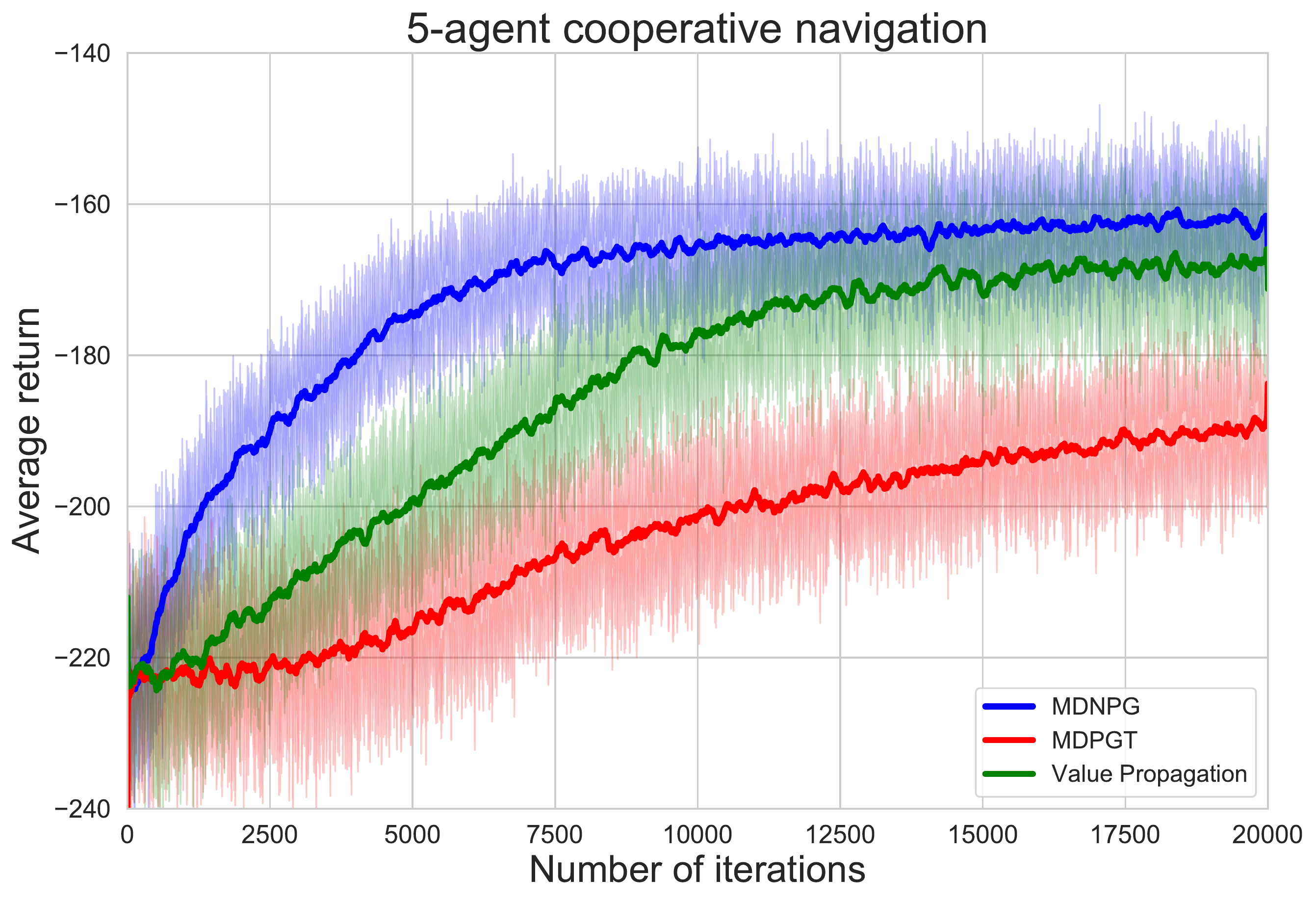}
         \caption[]{Fully-connected} 
         \label{fig: navigation_fc}
     \end{subfigure}
     
     \vskip\baselineskip
     \begin{subfigure}[]{0.49\textwidth}
         \centering
         \includegraphics[width=.95\linewidth]{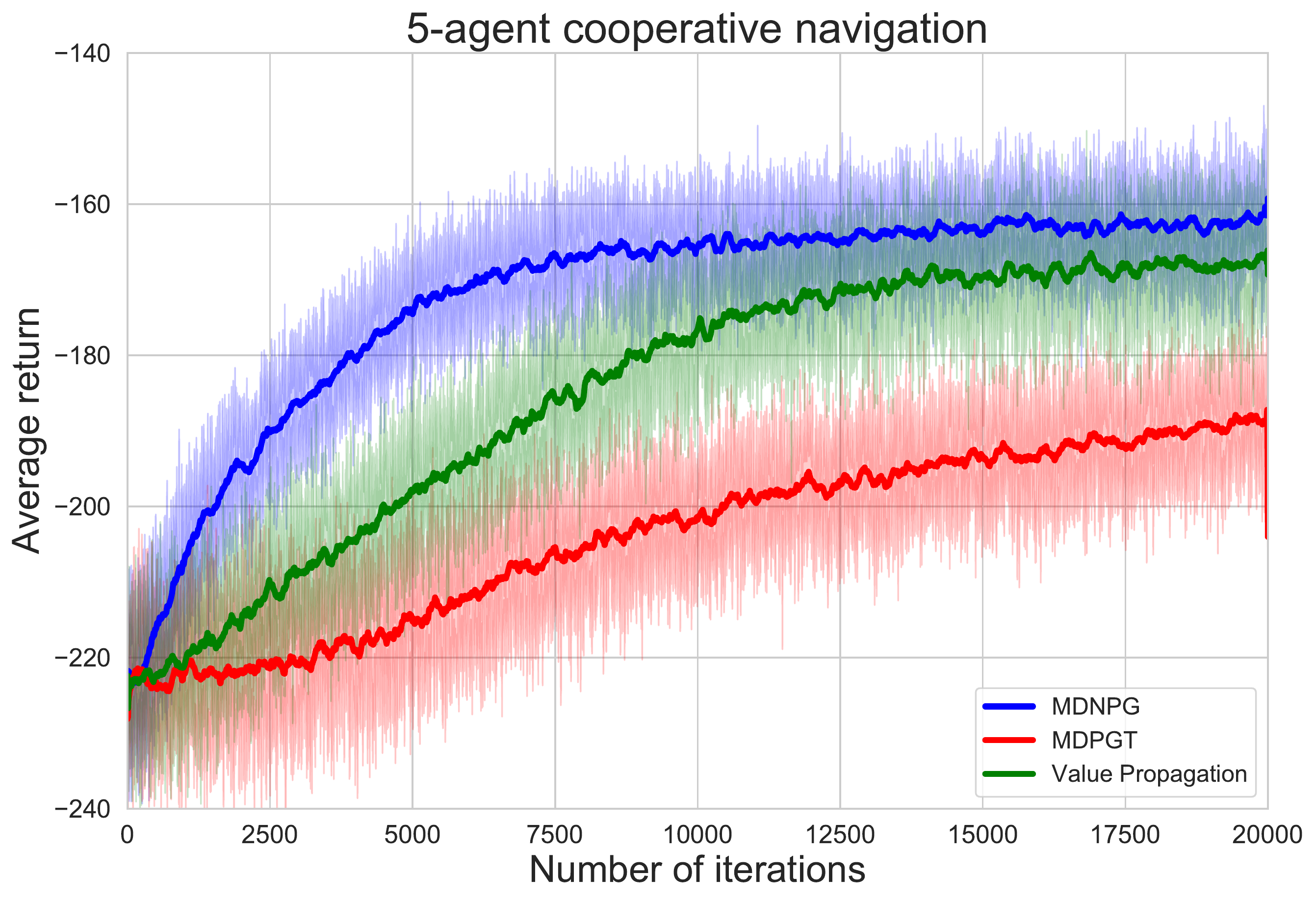}
         \caption[]{Bi-partite} 
         \label{fig: navigation_bi}
     \end{subfigure}
     \hfill
    \begin{subfigure}[]{0.49\textwidth}
    \centering
         \includegraphics[width=.95\linewidth]{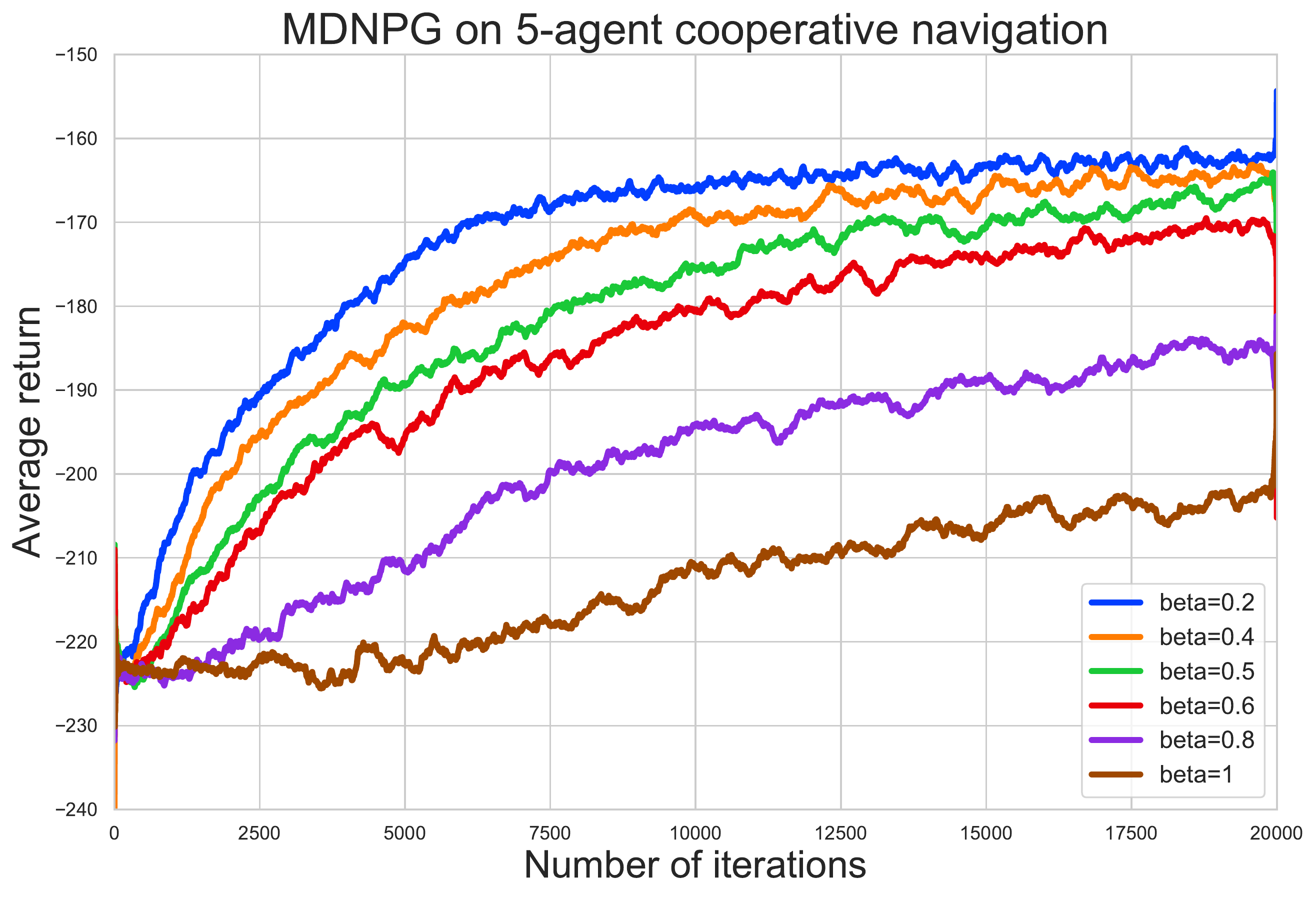}
         \caption[]{Influence of $\beta$ for the ring topology } 
         \label{fig: navigation_beta}
     \end{subfigure}
    \caption{Empirical results on cooperative navigation. In (a)--(c), plots of average return and standard deviation against  number of iterations over five random instances, where $\beta=0.2$ for MDNPG and MDPGT. In (d),  influence of $\beta$ on the performance MDNPG for the ring topology.} 
    \label{fig: navigation}
\end{figure}

In order to demonstrate the influence of the communication network on the  algorithms' performance, we follow the work in  \cite{jiang2021mdpgt} and  test three network topologies (see Figure~\ref{fig: topology}): ring, fully-connected, and bi-partite. The empirical results are displayed in Figure~\ref{fig: navigation}. It is evident from Figure~\ref{fig: navigation_ring}--\ref{fig: navigation_bi} that the performance of MDNPG is superior to the other two test methods in all the three network topologies. We have also tested the influence of the momentum parameter $\beta$ on the performance of MDNPG. Since similar trend has been observed for different topologies, only the results for the ring topology is presented in Figure~\ref{fig: navigation_beta}.

\subsubsection{Multi-task GridWorld}
For the MTRL setting in Section~\ref{section: marl},  experiments have been conducted on a multi-task GridWorld problem, whose setup is overall similar to the single-agent case in Section~\ref{singleagent} but with multiple individual environments. Each agent has its own environment but uses the same policy. By doing so, it is expected to obtain a policy with better generalization. 

\begin{figure}[ht!]
\centering
     \begin{subfigure}[]{0.49\textwidth}
     \centering
        \includegraphics[width=.95\linewidth]{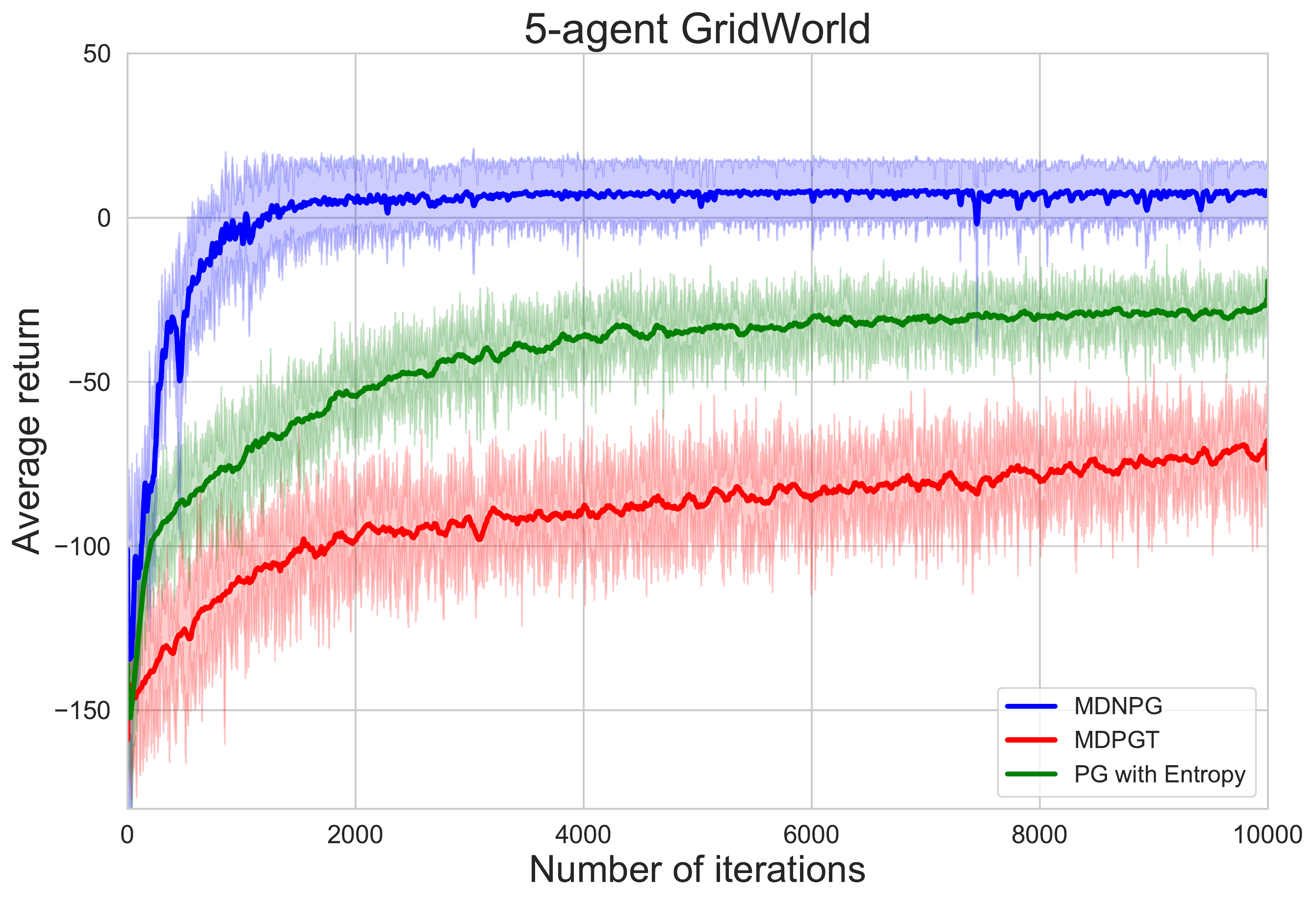}
        \caption[]{Ring} 
        \label{fig: task_ring}
     \end{subfigure}
     \hfill
     \begin{subfigure}[]{0.49\textwidth}
     \centering
        \includegraphics[width=.95\linewidth]{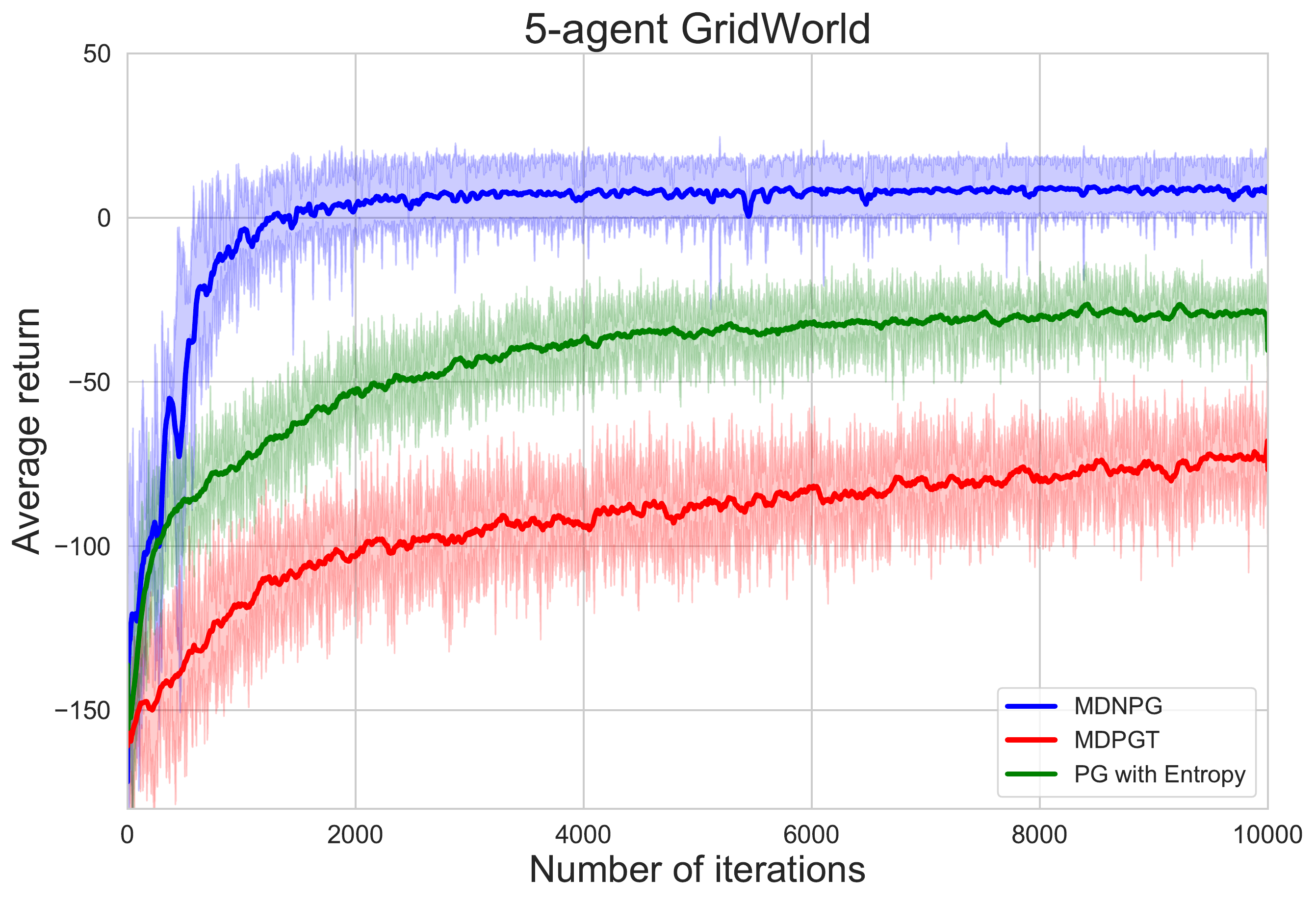}
        \caption[]{Fully-connected} 
        \label{fig: task_fc}
     \end{subfigure}

     \vskip\baselineskip
     \begin{subfigure}[]{0.49\textwidth}
        \centering
        \includegraphics[width=.95\linewidth]{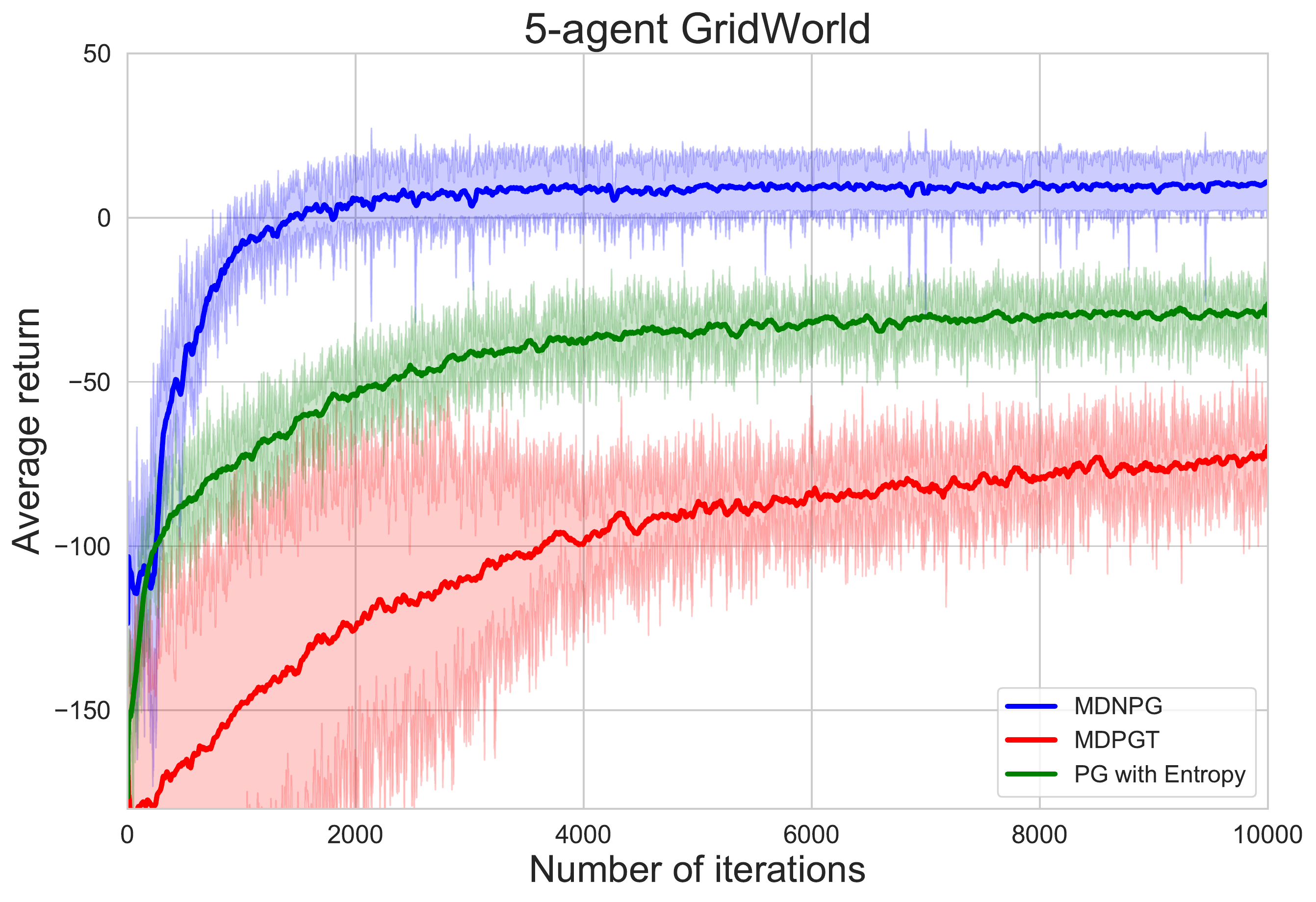}
        \caption[]{Bi-partite} 
        \label{fig: task_bi}
     \end{subfigure}
     \hfill
    \begin{subfigure}[]{0.49\textwidth}
    \centering
        \includegraphics[width=.95\linewidth]{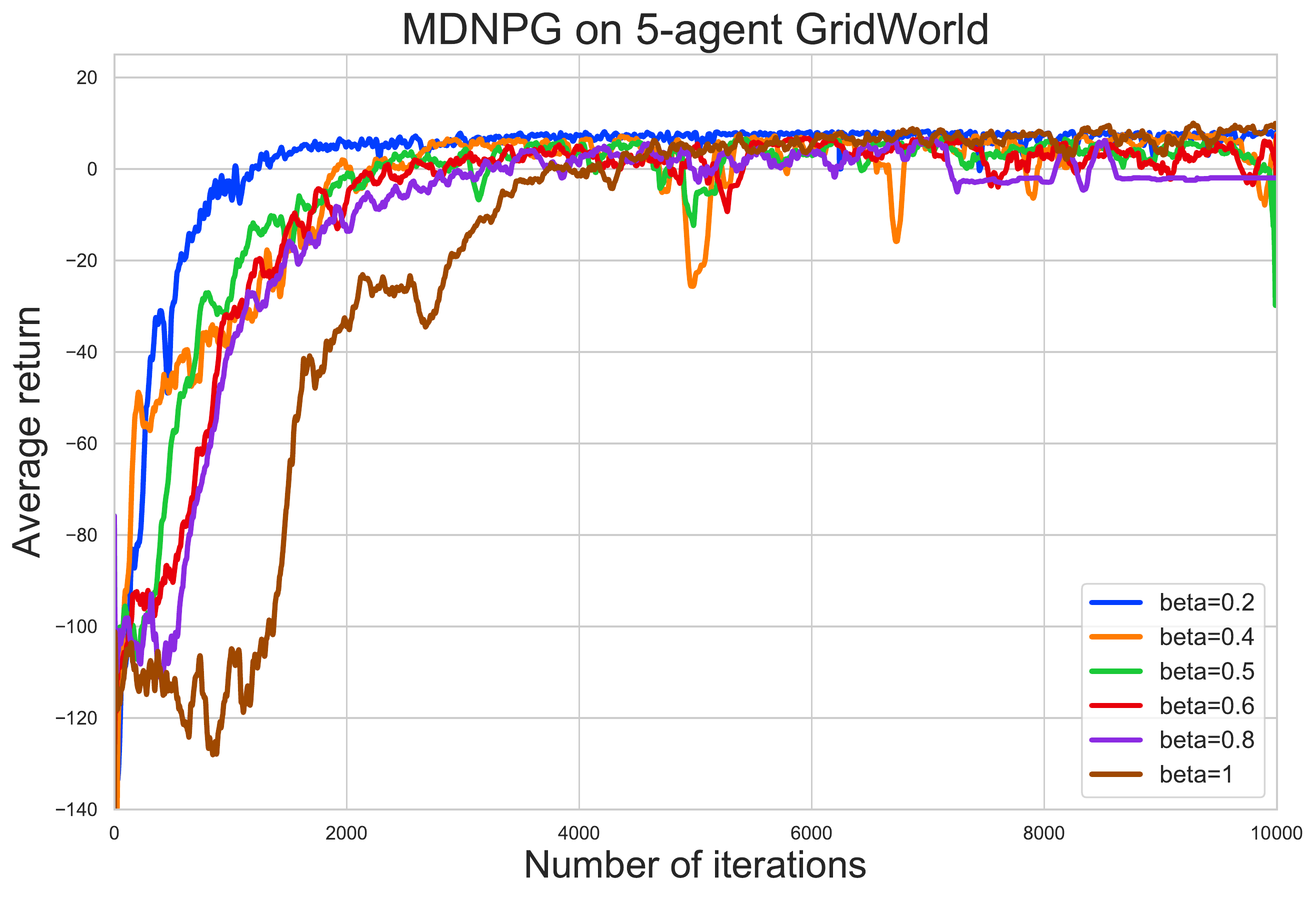}
        \caption[]{Influence of $\beta$ for the ring topology } 
        \label{fig: task_beta}
    \end{subfigure}
    \caption{Empirical results on multi-task GridWorld. In (a)--(c), plots of average return and standard deviation against number of iterations over five random instances where $\beta=0.2$ for MDNPG and MDPGT. In (d), influence of $\beta$ on the performance MDNPG for the ring topology.} 
    \label{fig: task}
\end{figure}

Five different yet similar environments are considered in our experiments and we  compare the proposed MDNPG with MDPGT~\cite{jiang2021mdpgt} and PG with entropy regularization~\cite{zeng2021decentralized}.  As with the single-agent case, one-hidden-layer (of size $128$) policy network and value network  with ReLU have been utilized, and all of the hyperparameters are properly tuned. Again, three network topologies have been tested and the empirical results are presented in Figure~\ref{fig: task_ring}--\ref{fig: task_bi}, which clearly shows that MDNPG outperforms the other two test methods. The influence of $\beta$ on the performance of MDNPG for the ring topology is presented in Figure~\ref{fig: task_beta}.


To evaluate the generalization effect of the  learned policies in the multi-task experiments, we  compare them with the policies learned by training each agent separately (that is, by solving $\max_{\btheta}V_i(\btheta)$ for each $i$ instead of solving $\max_{\btheta}\sum_{i=1}^n V_i(\btheta)$ as in \eqref{eq: marl1}). 
Table~\ref{tab: task} contains the average returns over 100 random trajectories computed from the policies trained in the multi-task environment (for the ring topology) as well as in each single environment. 
Note that even though the single-agent versions of MDNPG, MDPGT, and PG with entropy regularization are all tested for training each agent separately, only results for the single-agent MDNPG are presented due to its superior performance. It is clear that the shared policy learned by MDNPG overall generalizes better than the other two methods, and it is competitive with (or better than) the policy learned by training the individual agents in 4 out of 5 environments.

\begin{table}[ht!]
\centering
\caption{Average returns over 100 random trajectories based on different learned policies. Agent $i$ means that the policy is learned by training the $i$-th agent in the $i$-th grid using the single-agent MDNPG (i.e., the momentum-based NPG). In contrast, MDNPG, MDPGT, and PG with entropy regularization learn a shared policy by training the multi-task environment. }\label{tab: task}
\begin{tabular}{ c c c c c c c} 
\toprule[2pt]
 & Grid 1 & Grid 2 & Grid 3 & Grid 4 & Grid 5 & Sum\\
\midrule[1pt]
Agent 1 & {\bf7.64 } & -13.24&  -81.66& -136.32&  -14.1&  -237.68\\
Agent 2 & -10.95 &  {\bf 5.93 }& -81.09 &-11.94& -15.99&  -114.03\\
Agent 3  & -10.92& -17.97 &  {\bf 6.76}& -97.5&  -67.04 & -186.68\\
Agent 4   &-10.92 &-126.05 & -14.24  &  {\bf 9.25} & -14.04 & -156.01\\
Agent 5   & -11.02& -18.01 & -3.95 &-57.9  & {\bf  2.5} &  -88.37\\
\midrule[1pt]
MDNPG   &{\bf 7.94} &  -17.93  & {\bf 6.74} & {\bf 9.25} & {\bf 4.42} & 10.43\\
\midrule[1pt]
MDPGT   & -7.5 & -18.82 &  -7.04 & -3.37 & -20.25 &  -56.98\\
\midrule[1pt]
PG with entropy & 0.51& {\bf-11.82} &  -0.22 &  2.89 &-19.46 & -28.1\\
\bottomrule[2pt]
\end{tabular}
\end{table}

\section{Proof of main results} \label{section proofs}
We first introduce some convenient notations. Letting $\btheta_i\in\R^d$ be the local variable for $i$-th agent, we define the aggregated variable $\btheta\in\R^{nd}$ by 
\begin{align*}
    \btheta = \begin{bmatrix}
    \btheta_1^\tran &\cdots &\btheta_n^\tran
    \end{bmatrix}^\tran. 
\end{align*}
Let $\mH^t = \diag(\mH_1^t, \cdots, \mH_n^t)\in\R^{nd\times nd}$ be the block diagonal matrix and define $\vd_i^t = \mH_i^t \vy_i^{t+1}\in\R^d$. We apply the same aggregation rules to obtain other concatenated variables $\vy, \vv, \vd\in\R^{nd}$.  Using these notations, the key steps in Algorithm \ref{alg: moentum pg} can be rewritten in a more compact form:
\begin{align}
\label{update compact form}
	\vy^{t+1} = \pW(\vy^{t} + \vv^t - \vv^{t-1}) \text{ and } \btheta^{t+1} = \pW(\btheta^t + \eta \vd^{t}).
\end{align}
Let $\bbtheta\in\R^d$ be the average of $\bbtheta_i$ over all the agents, i.e, 
\begin{align*}
\bbtheta = \frac{1}{n}\sum_{i=1}^n \btheta_i = \frac{1}{n}(\bone_n^\tran \otimes \mI_d) \btheta.
\end{align*}
Similarly, $\by, \bv, \bd\in\R^d$ also denote the averages of related variables. By the update described in \eqref{update compact form}, it is straightforward to obtain that
\begin{align*}
\by^{t+1} = \bv^t \text{ and }
\bbtheta^{t+1} = \bbtheta^t + \eta \bd^t.
\end{align*}
Moreover, we define the aggregated gradient and averaged gradient by 
\begin{align*}
\tnablaV(\btheta) = \begin{bmatrix}
\nabla V_1(\btheta_1)^\tran &\cdots & \nabla V_n(\btheta_n)^\tran
\end{bmatrix}^\tran \in\R^{nd}, \text{ and } \bnablaV(\btheta) = \frac{1}{n}\sum_{i=1}^n \nabla V_i(\btheta_i)\in\R^d.
\end{align*}
 Throughout this work, we will frequently use the following relationship 
\begin{align}
\label{key eq}
	\twonorm{ \left((\mW - n^{-1} \mJ) \otimes \mI_d \right)\va } = \twonorm{ \left((\mW - n^{-1} \mJ) \otimes \mI_d \right)(\va - \bone_n \otimes \bar{\va} )} 
\end{align}
for any $\va=\begin{bmatrix}
	\va_1^\tran &\cdots &\va_n^\tran
\end{bmatrix}^\tran \in \R^{nd}$, where $\bar{\va} = n^{-1}\sum_{i=1}^{n}\va_i\in\R^{d}$ and $\mJ=\bone_n\bone_n^\tran\in\mathbb{R}^{n\times n}$. 

The following key lemma establishes the descent property of  MDNPG. This lemma may be of  broader interest in analyzing preconditioned stochastic first order methods in decentralized non-convex optimization.
\begin{lemma}
\label{lemma ascent}
Let $\{\btheta_i^t\}$ be generated by Algorithm \ref{alg: moentum pg} and $\Delta = V^\star - V(\btheta^0)$. Suppose $0 < \eta \leq \frac{\mu_F}{8L}$. Under Assumption \ref{assumption 0} and Assumption \ref{assumption 1}, one has
\begin{align*}
	\frac{1}{n}\sum_{t=0}^{T}\sum_{i=1}^{n} \twonorm{ \nabla V(\btheta_i^t)  }^2   	&\leq  \frac{8G^2 \Delta}{\eta\mu_F} - \frac{G^2}{n}\sum_{t=0}^{T}\twonorm{\vd^t}^2  + \frac{76G^2 }{\mu_F^2}  \sum_{t=0}^{T}\twonorm{\bnablaf (\btheta^t) - \bv^t}^2  \\
	&\qquad + \frac{10G^2 }{n\mu_F^2} \sum_{t=0}^{T} \twonorm{\vy^{t+1} - \bone_n\otimes \by^{t+1}}^2   + \frac{82G^2L^2  }{n\mu_F^2 }\sum_{t=0}^{T} \twonorm{\btheta^t - \bone_n \otimes \bbtheta^t}^2.
\end{align*}
\end{lemma}
Next, we will look for conditions on the  step size $\eta$ and weight factor $\beta$ such that 
\begin{align*}
&- \frac{G^2}{nT}\sum_{t=0}^{T}\twonorm{\vd^t}^2  + \frac{76G^2 }{T\mu_F^2}  \sum_{t=0}^{T}\twonorm{\bnablaf (\btheta^t) - \bv^t}^2 \\
&\qquad + \frac{10G^2 }{nT\mu_F^2} \sum_{t=0}^{T} \twonorm{\vy^{t+1} - \bone_n\otimes \by^{t+1}}^2   + \frac{82G^2L^2  }{nT\mu_F^2 }\sum_{t=0}^{T} \twonorm{\btheta^t - \bone_n \otimes \bbtheta^t}^2 = \calO\left(\eta, \beta, \frac{1}{B}, \frac{1}{T}\right).
\end{align*}
Assuming this holds, then the application of Lemma \ref{lemma ascent} will yield 
\begin{align*}
\frac{1}{n}\sum_{t=0}^{T}\sum_{i=1}^{n} \twonorm{ \nabla V(\btheta_i^t)  }^2   	&\leq  \frac{8G^2 \Delta}{\eta\mu_F} +\calO\left(\eta, \beta, \frac{1}{B}, \frac{1}{T}\right), 
\end{align*}
which implies the convergence of Algorithm \ref{alg: moentum pg}.  
To achieve this goal, we need several lemmas whose proofs are either deferred to Section~\ref{proof} or already given in the literature.

\begin{lemma}
	\label{errorX}
	Under Assumptions \ref{assumption 6}, \ref{assumption 1} and \ref{assumption 4}, for all $t\geq 0 $, one has 
	\begin{align}
\label{eqErrorX}
	\twonorm{\btheta^{t+1} - \bone_n \otimes \bbtheta^{t+1}}^2 &\leq \frac{1+\rho^2}{2} \twonorm{  \btheta^t - \bone_n \otimes \bbtheta^t  }^2 + \frac{4\eta^2}{\mu_F^2(1-\rho^2)}\twonorm{\vy^{t+1}- \bone_n \otimes \by^{t+1}}^2  + \frac{G^2\eta^2}{\mu_F^2(1-\rho^2)}  \twonorm{\vd^t}^2.
	\end{align}
Moreover, one has
\begin{align}
	\label{eqErrorX2}
\twonorm{\btheta^{t+1} - \bone_n \otimes \bbtheta^{t+1}}^2 &\leq 2\rho^2\twonorm{  \btheta^t - \bone_n \otimes \bbtheta^t  }^2 + \frac{4\eta^2\rho^2 }{\mu_F^2} \twonorm{\vy^{t+1}- \bone_n \otimes \by^{t+1}}^2+  \frac{G^2\eta^2\rho^2}{\mu_F^2} \twonorm{\vd^t}^2.
\end{align}

\end{lemma}


\begin{lemma} 
\label{lemma tracking}
Let $\{\vy^t\}$ be generated by Algorithm \ref{alg: moentum pg}. Under Assumptions~\ref{assumption 6} and \ref{assumption 1}-\ref{assumption 4}, we have
\begin{align}
\label{eq Y1}
	\E{\twonorm{\vy^1 - \bone_n \otimes \by^1}^2} \leq \frac{n\rho^2\bar{\nu}^2}{B} + \rho^2\sum_{i=1}^{n} \twonorm{\nabla V_i(\bbtheta^0) }^2.
\end{align}
Furthermore, if $\eta \leq \frac{\mu_F(1-\rho^2)}{24\sqrt{2}\Phi}$ and $0\leq\beta \leq 1$, then for all $t\geq 1$, one has
\begin{align}
	 \label{Error Y}
	\E{\twonorm{\vy^{t+1} - \bone_n \otimes \by^{t+1}}^2} &\leq \frac{3+\rho^2}{4}  \E{\twonorm{ \vy^t -\bone_n \otimes \by^t}^2 }+ \frac{8n\beta^2 \bar{\nu}^2 }{1-\rho^2} \notag \\
	&\qquad + \frac{8\beta^2 }{1-\rho^2}  \E{\twonorm{ \tnablaV(\btheta^{t-1}) - \vv^{t-1}}^2} \notag\\
	&\qquad +\frac{216\Phi^2}{1-\rho^2}  \E{ \twonorm{  \btheta^{t-1}- \bone_n \otimes \bbtheta^{t-1}}^2 } \notag\\
	&\qquad + \frac{144\Phi^2G^2\eta^2}{\mu_F^2 (1-\rho^2)}  \E{\twonorm{\vd^{t-1}}^2}.
\end{align}    
\end{lemma}

\begin{lemma}[Lemma 6 in \cite{xin2021hybrid}]
	\label{lemma6Xin}
	Let $\{a_t\}, \{b_t\}$ and $\{c_t\}$ be nonnegative sequences and $d>0$ be some constant such that $a_t\leq \gamma a_{t-1} + \gamma b_{t-1} + c_{t} + d$ for $t\geq 1$, where $\gamma \in (0,1)$. Then for $T\geq 1$, we have
	\begin{align}
		\label{lemma6a}
		\sum_{t=0}^{T} a_t \leq \frac{1}{1-\gamma} a_0 + \frac{1}{1-\gamma} \sum_{t=0}^{T-1}b_t + \frac{1}{1-\gamma} \sum_{t=1}^{T} c_t + \frac{dT}{1-\gamma}.
	\end{align}
	Moreover, if $a_{t+1}\leq \gamma a_{t} + b_{t-1} + d$ for $t\geq 1$, then for $T\geq 2$, one has 
	\begin{align}
		\label{lemma6b}
		\sum_{t=1}^{T}a_t \leq \frac{1}{1-\gamma} a_1 + \frac{1}{1-\gamma} \sum_{t=0}^{T-2}b_t + \frac{dT}{1-\gamma}.
	\end{align}
\end{lemma}

\begin{lemma}
\label{accumulated Error Y}
Let 
\begin{align*}
	A_1 = \frac{4n\rho^2\bar{\nu}^2}{B(1-\rho^2)}  + \frac{32n T\beta^2 \bar{\nu}^2 }{(1-\rho^2)^2}  + \frac{32n\beta \bar{\nu}^2 }{B(1-\rho^2)^2}   +  \frac{64nT  \bar{\nu}^2 \beta^3}{(1-\rho^2)^2}  + \frac{4\rho^2}{1-\rho^2} \twonorm{\tnablaV(\btheta^0)}^2.
\end{align*}
Suppose $\eta \leq \frac{\mu_F(1-\rho^2)}{24\sqrt{2}\Phi}$ and $\beta <1$. 
Then one has
\begin{align}
	\label{error accumulated Y}
	\sum_{t= 1}^{T} \E{\twonorm{ \vy^t - \bone_n \otimes \by^t}^2 } \leq A_1 + \frac{1632\Phi^2 }{(1-\rho^2)^2}\sum_{t=0}^{T}  \E{\twonorm{  \btheta^{t}- \bone_n \otimes \bbtheta^{t}}^2 }  +  \frac{960\Phi^2\kappa_F^2\eta^2 }{(1-\rho^2)^2}    \sum_{t=0}^{T-1} \E{\twonorm{\vd^t}^2} .
\end{align}
\end{lemma}

\begin{lemma}
	\label{accumulatedErrorX}
	Suppose \begin{align*}
		0 < \eta <\frac{\mu_F(1-\rho^2)^3}{ \kappa_{F}\sqrt{1632000(L^2 + \Phi^2) }},
	\end{align*} and $\beta < 1$. Then for any $T\geq 1$, one has
	\begin{align*}
		\sum_{t=0}^{T} \E{\twonorm{\btheta^{t} - \bone_n \otimes \bbtheta^{t}}^2} \leq\frac{16A_1\eta^2}{\mu_F^2(1-\rho^2)^2} + \frac{10G^2\eta^2}{\mu_F^2(1-\rho^2)^3} \sum_{t=0}^{T}\E{\twonorm{\vd^t}^2}.
	\end{align*}
\end{lemma}
 
\begin{lemma}[Lemma 8 in \cite{jiang2021mdpgt}]
	\label{lemma jiang2021}
	Let $\vv^t$ and $\btheta^t$ be generated by Algorithm \ref{alg: moentum pg} and let $\Phi^2=L_g^2+C_g^2C_\omega^2$. Then under Assumption \ref{assumption 1}, \ref{assumption 3} and \ref{assumption 5}, for any $t\geq 1$, one has
	\begin{align}
		\label{eq v1}
		\sum_{t=0}^{T} \E{\twonorm{\bv^t - \bnablaV(\btheta^t)}^2} &\leq \frac{\bar{\nu}^2}{n\beta B} + \frac{2\beta \bar{\nu}^2 T}{n}+ \frac{12\Phi^2\eta^2}{n\beta} \sum_{t=0}^{T-1} \E{\twonorm{\bd^t}^2} + \frac{24\Phi^2}{\beta n^2} \sum_{t=0}^{T} \E{\twonorm{\vx^t - \bone_n \otimes \bx^t}^2}
	\end{align}
	and
	\begin{align}
	\label{eq v2}
	\sum_{t=0}^{T}\E{\twonorm{\vv^t - \widetilde{\nabla V}(\btheta^t)}^2} &\leq \frac{n\bar{\nu}^2}{\beta B} + 2n \beta T \bar{\nu}^2 + \frac{12n  \eta^2\Phi^2}{\beta}\sum_{t=0}^{T-1} \E{\twonorm{\bd^t}^2}+ \frac{24\Phi^2}{\beta} \sum_{t=0}^{T} \E{\twonorm{\vx^t - \bone_n \otimes \bx^t}^2}.
	\end{align}
\end{lemma}

\subsection{Proof of Theorem \ref{main result}}
Lemma \ref{lemma ascent} implies that 
\begin{align*}
	&\frac{1}{n}\sum_{t=0}^{T}\sum_{i=1}^{n} \E{\twonorm{ \nabla V(\btheta_i^t)  }^2}   	\\
	&\leq  \frac{8G^2 \Delta}{\eta\mu_F} - \frac{G^2}{n}\sum_{t=0}^{T} \E{\twonorm{\vd^t}^2}  + \frac{10G^2 }{n\mu_F^2} \sum_{t=0}^{T} \E{\twonorm{\vy^{t+1} - \bone_n\otimes \by^{t+1}}^2 }  + \frac{82G^2L^2  }{n\mu_F^2 }\sum_{t=0}^{T} \E{\twonorm{\btheta^t - \bone_n \otimes \bbtheta^t}^2}\\
	&\qquad + \frac{76G^2 }{\mu_F^2}  \sum_{t=0}^{T}\E{ \twonorm{\bnablaf (\btheta^t) - \bv^t}^2} \\
	&\stackrel{(a)}{\leq } \frac{8G^2 \Delta}{\eta\mu_F} - \frac{G^2}{n}\sum_{t=0}^{T} \E{\twonorm{\vd^t}^2}+ \frac{10G^2 }{n\mu_F^2} \sum_{t=0}^{T} \E{\twonorm{\vy^{t+1} - \bone_n\otimes \by^{t+1}}^2}   + \frac{82G^2L^2  }{n\mu_F^2 }\sum_{t=0}^{T} \E{\twonorm{\btheta^t - \bone_n \otimes \bbtheta^t}^2}\\
	&\qquad + \frac{76G^2 }{\mu_F^2}   \bigg( \frac{\bar{\nu}^2}{n\beta B} + \frac{2\beta \bar{\nu}^2 T}{n}+ \frac{12\Phi^2\eta^2}{n\beta} \sum_{t=0}^{T-1} \E{\twonorm{\bd^t}^2} + \frac{24\Phi^2}{\beta n^2} \sum_{t=0}^{T} \E{\twonorm{\vx^t - \bone_n \otimes \bx^t}^2} \bigg)\\
	&=\frac{8G^2 \Delta}{\eta\mu_F} + \frac{76\bar{\nu}^2 G^2}{n\beta B \mu_F^2}   + \frac{152\beta T\bar{\nu}^2 G^2 }{n \mu_F^2}  -  \frac{G^2}{n}  \E{\sum_{t=0}^{T}\twonorm{\vd^t}^2}  +\frac{912\Phi^2\eta^2G^2}{n\beta \mu_F^2} \sum_{t=0}^{T-1} \E{\twonorm{\bd^t}^2} \\
	&\qquad  + \left( \frac{1824\Phi^2G^2}{\beta n^2 \mu_F^2}   +  \frac{82G^2L^2  }{n\mu_F^2 } \right)\sum_{t=0}^{T} \E{ \twonorm{\btheta^t - \bone_n \otimes \bbtheta^t}^2 }  + \frac{10G^2 }{n\mu_F^2} \sum_{t=0}^{T} \E{\twonorm{\vy^{t+1} - \bone_n\otimes \by^{t+1}}^2 }\\
	&\stackrel{(b)}{\leq} \frac{8G^2 \Delta}{\eta\mu_F} + \frac{76\bar{\nu}^2 G^2}{n\beta B \mu_F^2}   + \frac{152\beta T\bar{\nu}^2 G^2 }{n \mu_F^2}  -  \frac{G^2}{n}  \sum_{t=0}^{T} \E{\twonorm{\vd^t}^2} +\frac{912\Phi^2\eta^2G^2}{n^2\beta \mu_F^2} \sum_{t=0}^{T-1} \E{\twonorm{\vd^t}^2} \\
	&\qquad  + \left( \frac{1824\Phi^2G^2}{\beta n^2 \mu_F^2}   +  \frac{82G^2L^2  }{n\mu_F^2 } \right)\sum_{t=0}^{T} \E{ \twonorm{\btheta^t - \bone_n \otimes \bbtheta^t}^2 }  + \frac{10G^2 }{n\mu_F^2} \sum_{t=0}^{T} \E{\twonorm{\vy^{t+1} - \bone_n\otimes \by^{t+1}}^2 }, \numberthis \label{eq tmp theorem a}
\end{align*}
where step (a) follows from \eqref{eq v1} and step (b) is due to $\twonorm{\bd}^2 \leq \frac{1}{n}\twonorm{\vd}^2$.  Since 
\begin{align*}
	0 < \eta < \frac{\mu_F(1-\rho^2)^3}{ \kappa_{F}\sqrt{1632000(L^2 + \Phi^2) }} \text{ and }\frac{1632000(L^2 + \Phi^2) \kappa_F^2 \eta^2 }{n\mu_F^2(1-\rho^2)^6} \leq \beta <\frac{1}{n},
\end{align*}
it implies that 
\begin{align*}
	\frac{912\Phi^2\eta^2G^2}{n^2\beta \mu_F^2} &\leq \frac{912\Phi^2\eta^2G^2}{n^2  \mu_F^2} \cdot \frac{n\mu_F^2(1-\rho^2)^6 }{1632000 (L^2 + \Phi^2) \kappa_F^2 \eta^2} \leq \frac{G^2}{2n},\\
	\frac{1824\Phi^2G^2}{\beta n^2 \mu_F^2}   +  \frac{82G^2L^2  }{n\mu_F^2 } &\leq \frac{1824(L^2 + \Phi^2) G^2}{ n \mu_F^2} \left( 1+ \frac{1}{\beta n}\right) \\
	&\leq \frac{1824(L^2 + \Phi^2) G^2}{ n \mu_F^2}  \cdot \frac{2}{\beta n} \\
	&\leq \frac{1824(L^2 + \Phi^2) G^2}{ n \mu_F^2}  \cdot  \frac{2(1-\rho^2)^6 \mu_F^2 }{1632000(L^2 + \Phi^2) \kappa_F^2 \eta^2}\\
	&\leq \frac{(1-\rho^2)^6 G^2}{100n\kappa_F^2 \eta^2}.
\end{align*}
Plugging these inequalities into \eqref{eq tmp theorem a} yields that 
\begin{align*}
	\frac{1}{n}\sum_{t=0}^{T}\sum_{i=1}^{n} \E{ \twonorm{ \nabla V(\btheta_i^t)  }^2}   	 &\leq \frac{8G^2 \Delta}{\eta\mu_F} + \frac{76\bar{\nu}^2 G^2}{n\beta B \mu_F^2}   + \frac{152\beta T\bar{\nu}^2 G^2 }{n \mu_F^2}  -  \frac{G^2}{2n}  \sum_{t=0}^{T} \E{\twonorm{\vd^t}^2 } \\
	 &\qquad +\frac{(1-\rho^2)^6G^2}{100n\kappa_F^2 \eta^2}\sum_{t=0}^{T} \E{ \twonorm{\btheta^t - \bone_n \otimes \bbtheta^t}^2 } 
	  + \frac{10G^2 }{n\mu_F^2} \sum_{t=0}^{T} \E{\twonorm{\vy^{t+1} - \bone_n\otimes \by^{t+1}}^2 }\\
	&\stackrel{(a)}{\leq}\frac{8G^2 \Delta}{\eta\mu_F} + \frac{76\bar{\nu}^2 G^2}{n\beta B \mu_F^2}   + \frac{152\beta T\bar{\nu}^2 G^2 }{n \mu_F^2}  -  \frac{G^2}{2n}  \E{\sum_{t=0}^{T}\twonorm{\vd^t}^2} \\ 
	&\qquad +\frac{(1-\rho^2)^6 G^2}{100n\kappa_F^2 \eta^2}\sum_{t=0}^{T} \E{ \twonorm{\btheta^t - \bone_n \otimes \bbtheta^t}^2 } \\
	&\qquad + \frac{10G^2 }{n\mu_F^2}  \bigg( A_1 + \frac{1632\Phi^2 }{(1-\rho^2)^2}\sum_{t=0}^{T}  \E{\twonorm{  \btheta^{t}- \bone_n \otimes \bbtheta^{t}}^2 }  +  \frac{960\Phi^2\kappa_F^2\eta^2 }{(1-\rho^2)^2}    \sum_{t=0}^{T-1} \E{\twonorm{\vd^t}^2}  \bigg)\\
	&= \frac{8G^2 \Delta}{\eta\mu_F} + \frac{76\bar{\nu}^2 G^2}{n\beta B \mu_F^2}   + \frac{152\beta T\bar{\nu}^2 G^2 }{n \mu_F^2}    + \frac{10A_1G^2 }{n\mu_F^2}  - \left( \frac{G^2}{2n} - \frac{10G^2 }{n\mu_F^2}  \cdot \frac{960 \Phi^2\kappa_F^2\eta^2 }{(1-\rho^2)^2}\right) \sum_{t=0}^{T}\E{\twonorm{\vd^t}^2}  \\
	&\qquad + \left( \frac{(1-\rho^2)^6 G^2}{100n\kappa_F^2 \eta^2} + \frac{16320G^2\Phi^2 }{n\mu_F^2(1-\rho^2)^2}\right)\sum_{t=0}^{T} \E{ \twonorm{\btheta^t - \bone_n \otimes \bbtheta^t}^2 }, \numberthis \label{eq theorem b}
\end{align*}
where step (a) follows from \eqref{error accumulated Y}. 
Using the conditions for $\eta$ and $\beta$ again gives that
\begin{align*}
	\frac{10G^2 }{n\mu_F^2} \cdot  \frac{960\Phi^2\kappa_F^2\eta^2 }{(1-\rho^2)^2}  &\leq  \frac{10G^2 }{n\mu_F^2} \cdot  \frac{960\Phi^2\kappa_F^2  }{(1-\rho^2)^2} \cdot \frac{\mu_F^2(1-\rho^2)^6}{\kappa_F^2 \cdot 1632000 (L^2 + \Phi^2) }\leq \frac{G^2}{4n},\\
	\frac{(1-\rho^2)^6G^2}{100n\kappa_F^2 \eta^2} + \frac{16320G^2\Phi^2 }{n\mu_F^2 (1-\rho^2)^2}&= \frac{(1-\rho^2)^6G^2}{100n\kappa_F^2 \eta^2} + \frac{16320\kappa_F^2\Phi^2 }{n(1-\rho^2)^2}\\ 
	&\leq \frac{(1-\rho^2)^6G^2}{100n\kappa_F^2 \eta^2} +\frac{16320\kappa_{F}^2(L^2+\Phi^2) }{n (1-\rho^2)^2}\\
	&\leq \frac{(1-\rho^2)^6G^2}{100n\kappa_F^2 \eta^2} +\frac{16320}{n(1-\rho^2)^2}\cdot \frac{\mu_F^2(1-\rho^2)^6}{1632000\eta^2} \\
	&=\frac{(1-\rho^2)^6G^2}{100n\kappa_F^2 \eta^2} +\frac{G^2(1-\rho^2)^4}{100n\kappa_F^2\eta^2}\\
    &\leq\frac{(1-\rho^2)^4 G^2}{50n\kappa_F^2 \eta^2} .
\end{align*}
Substituting these inequalities into \eqref{eq theorem b} leads to
\begin{align*}
	\frac{1}{n}\sum_{t=0}^{T}\sum_{i=1}^{n} \E{\twonorm{ \nabla V(\btheta_i^t)  }^2}   	 
	&\leq \frac{8G^2 \Delta}{\eta\mu_F} + \frac{76\bar{\nu}^2 G^2}{n\beta B \mu_F^2}   + \frac{152\beta T\bar{\nu}^2 G^2 }{n \mu_F^2}    + \frac{10A_1G^2 }{n\mu_F^2}  -  \frac{G^2}{4n}  \sum_{t=0}^{T}\E{\twonorm{\vd^t}^2 } \\
	&\qquad + \frac{(1-\rho^2)^4 G^2}{50n\kappa_F^2 \eta^2} \sum_{t=0}^{T} \E{ \twonorm{\btheta^t - \bone_n \otimes \bbtheta^t}^2 }\\
	&\leq \frac{8G^2 \Delta}{\eta\mu_F} + \frac{76\bar{\nu}^2 G^2}{n\beta B \mu_F^2}   + \frac{152\beta T\bar{\nu}^2 G^2 }{n \mu_F^2}    + \frac{10A_1G^2 }{n\mu_F^2}  -  \frac{G^2}{4n}  \sum_{t=0}^{T}\E{\twonorm{\vd^t}^2}  \\
	&\qquad + \frac{(1-\rho^2)^4 G^2}{50n\kappa_F^2 \eta^2} \left( \frac{16A_1\eta^2}{\mu_F^2(1-\rho^2)^2} + \frac{10G^2\eta^2}{\mu_F^2(1-\rho^2)^3} \sum_{t=0}^{T}\E{\twonorm{\vd^t}^2}\right)\\
	&\leq \frac{8G^2 \Delta}{\eta\mu_F} + \frac{76\bar{\nu}^2 G^2}{n\beta B \mu_F^2}   + \frac{152\beta T\bar{\nu}^2 G^2 }{n \mu_F^2}    + \frac{10A_1G^2 }{n\mu_F^2}   + \frac{16A_1G^2 }{50n\mu_F^2}   \\
	&\qquad - \left( \frac{G^2}{4n}  - \frac{G^2}{5n}  \right)\sum_{t=0}^{T}\E{\twonorm{\vd^t}^2} \\
	&\leq  \frac{8G^2 \Delta}{\eta\mu_F} + \frac{76\bar{\nu}^2 G^2}{n\beta B \mu_F^2}   + \frac{152\beta T\bar{\nu}^2 G^2 }{n \mu_F^2}    + \frac{11A_1G^2 }{n\mu_F^2}, 
\end{align*}
where the second inequality follows from Lemma \ref{accumulatedErrorX}. Since $\btheta_{\out}$ is sampled uniformly from $\{\btheta_i^t\}_{i=1,\ldots, n; t = 0,\ldots, T}$, we have
\begin{align*}
	\E{\twonorm{\nabla V(\btheta_{\out}) }^2} &= \frac{1}{n(T+1)} \sum_{t=0}^{T} \sum_{i=1}^{n}\E{ \twonorm{ \nabla V(\btheta_i^t) }^2} \\
	&\leq   \frac{1}{nT} \sum_{t=0}^{T} \sum_{i=1}^{n}\E{ \twonorm{ \nabla V(\btheta_i^t) }^2}\\
	&\leq \frac{8G^2 \Delta}{T\eta\mu_F} + \frac{76\bar{\nu}^2 G^2}{nT\beta B \mu_F^2}   + \frac{152\beta \bar{\nu}^2 G^2 }{n \mu_F^2}    + \frac{11A_1G^2 }{nT\mu_F^2}\\
	&=\frac{8G^2 \Delta}{T\eta\mu_F} + \frac{76\bar{\nu}^2 \kappa_F^2}{nT\beta B }   + \frac{152\beta \bar{\nu}^2 \kappa_F^2 }{n}   \\
	&\qquad + \frac{11G^2 }{nT\mu_F^2} \left( \frac{4n\rho^2\bar{\nu}^2}{B(1-\rho^2)}  + \frac{32n T\beta^2 \bar{\nu}^2 }{(1-\rho^2)^2}  + \frac{32n\beta \bar{\nu}^2 }{B(1-\rho^2)^2}   +  \frac{64nT  \bar{\nu}^2 \beta^3}{(1-\rho^2)^2}  + \frac{4\rho^2}{1-\rho^2} \twonorm{\tnablaV(\btheta^0)}^2  \right)\\
	&=\frac{8G^2 \Delta}{T\eta\mu_F} + \frac{76\bar{\nu}^2 \kappa_F^2}{nT\beta B }   + \frac{152\beta \bar{\nu}^2 \kappa_F^2 }{n}   \\
	&\qquad +     \frac{44\rho^2\bar{\nu}^2 \kappa_F^2}{TB(1-\rho^2)}  + \frac{352\beta^2 \bar{\nu}^2\kappa_F^2 }{(1-\rho^2)^2}  + \frac{352\beta \bar{\nu}^2\kappa_F^2 }{TB(1-\rho^2)^2}   +  \frac{704  \bar{\nu}^2 \kappa_F^2 \beta^3}{(1-\rho^2)^2}  + \frac{44\rho^2\kappa_F^2}{nT(1-\rho^2)} \twonorm{\tnablaV(\btheta^0)}^2,
\end{align*}
which completes the proof of the main result.

\subsection{Proof of Corollary \ref{corollary}}
Since
\begin{align*}
	\eta = \frac{\mu_F n^{2/3}}{\kappa_{F} \sqrt{L^2 + \Phi}T^{1/3}}, \beta = \frac{n^{1/3}}{T^{2/3}} \text{ and } B = \left\lceil \frac{T^{1/3}}{n^{2/3}} \right\rceil,
\end{align*}
we have
\begin{align*}
	\frac{ 8\Delta G^2}{ \eta T\mu_F} &= \frac{ 8\Delta G^2}{ T\mu_F} \cdot \frac{\kappa_{F} \sqrt{L^2 + \Phi}T^{1/3} }{\mu_F n^{2/3}} =    \frac{8\Delta   \kappa_{F}^3\sqrt{L^2 + \Phi} }{(nT)^{2/3}},\\
	\frac{76\bar{\nu}^2\kappa_{F}^2}{n\beta TB  } &\leq  76\bar{\nu}^2\kappa_{F}^2\cdot \frac{1}{nT}\cdot  \frac{T^{2/3}} {n^{1/3}}\cdot \frac{n^{2/3}} {T^{1/3}} = \frac{76\bar{\nu}^2\kappa_{F}^2}{(nT)^{2/3}},\\
	 \frac{152\beta \bar{\nu}^2 \kappa_{F}^2 }{n } &=\frac{152\bar{\nu}^2\kappa_{F}^2}{(nT)^{2/3}},\\
	  \frac{44\rho^2\bar{\nu}^2\kappa_{F}^2}{BT(1-\rho^2)}  &\leq \frac{44\rho^2\bar{\nu}^2 \kappa_{F}^2}{(1-\rho^2)^2}\cdot \frac{n^{2/3}}{T^{4/3}},\\
	  \frac{352 \beta^2 \bar{\nu}^2 \kappa_{F}^2}{ (1-\rho^2)^2}  &\leq \frac{352 \bar{\nu}^2 \kappa_{F}^2 }{(1-\rho^2)^2}  \cdot \frac{n^{2/3}}{T^{4/3}},\\
	  \frac{352\beta \bar{\nu}^2\kappa_{F}^2 }{BT (1-\rho^2)^2}  &\leq\frac{352 \bar{\nu}^2 \kappa_{F}^2}{ (1-\rho^2)^2} \cdot \frac{1}{T}\cdot  \frac{n^{1/3}}{T^{2/3}} \cdot  \frac{n^{2/3}}{T^{1/3}} = \frac{352 \bar{\nu}^2 \kappa_{F}^2}{ (1-\rho^2)^2}\cdot \frac{n}{T^2},\\
	   \frac{704 \bar{\nu}^2 \kappa_{F}^2\beta^3}{ (1-\rho^2)^2} &=\frac{704 \bar{\nu}^2 \kappa_{F}^2  }{ (1-\rho^2)^2} \cdot \frac{n}{T^2}.
\end{align*}
Thus it can be seen that 
\begin{align*}
	\E{\twonorm{\nabla V(\btheta_{\out})}^2} & \leq  \frac{8\Delta   \kappa_{F}^3\sqrt{L^2 + \Phi}  + 228\bar{\nu}^2\kappa_F^2}{(nT)^{2/3}}+\frac{44\rho^2\kappa_{F}^2 \twonorm{\tnablaV(\btheta^0)}^2}{ (1-\rho^2)} \cdot \frac{1}{nT}+ \frac{396\kappa_{F}^2\bar{\nu}^2}{ (1-\rho^2)^2}\cdot \frac{n^{2/3}}{T^{4/3}}+ \frac{1056 \bar{\nu}^2  \kappa_{F}^2 }{ (1-\rho^2)^2} \cdot \frac{n}{T^2},
\end{align*}
which completes the proof of the corollary.

\section{Proofs}
\label{proof}

\subsection{Proof of Lemma~\ref{lemma: FIM}}
For simplicity, let $\btheta_i = [{\vvx^1}^\tran, \cdots, {\vvx^n}^\tran]^\tran\in\R^{d}$, where $\vvx^j\in\R^{d_j}$ and $d=\sum_{j=1}^n d_j$. 
From the definition~\eqref{eq: pi1} of the policy in collaborative RL, we have 
\begin{align*}
\nabla_{\btheta_i} \log \pi_{\btheta_i}(\va^h | \vs^h) &= \nabla_{\btheta_i} \sum_{j=1}^n\log \pi_{\vvx^j}(\va_j^h | \vs^h)\\&=\begin{bmatrix}\nabla_{\vvx^1} \log \pi_{\vvx^1}(\va_1^h | \vs^h) \\\vdots\\\nabla_{\vvx^n} \log \pi_{\vvx^n}(\va_n^h | \vs^h) \\\end{bmatrix}\in\R^{d\times 1}.
\end{align*}
Then the $(j,\ell)$-th block of $\mF_i(\btheta_i)$ is given by
\begin{align*}
    [\mF_i(\btheta_i)]_{j,\ell} = \E[\tau\sim p(\cdot | \btheta_i)]{\frac{1}{H}\sum_{h=0}^{H-1} \nabla_{\vvx^j} \log \pi_{\vvx^j}(\va_j^h|\vs^h) \left( \nabla_{\vvx^\ell} \log \pi_{\vvx^\ell}(\va_\ell^h|\vs^h) \right)^\tran }\in\R^{d_j\times d_\ell}.
\end{align*}
We will show that $[\mF_i(\btheta_i)]_{j,\ell}=\bzero$ for any $j\neq \ell$. To this end, for any $\alpha\in [d_j], \beta\in [d_\ell]$, one has
\begin{align*}
&\left[\E[\tau\sim p(\cdot | \btheta_i)]{\frac{1}{H}\sum_{h=0}^{H-1} \nabla_{\vvx^j} \log \pi_{\vvx^j}(\va_j^h|\vs^h) \left( \nabla_{\vvx^\ell} \log \pi_{\vvx^\ell}(\va_\ell^h|\vs^h) \right)^\tran }\right]_{\alpha, \beta}\\
=&\frac{1}{H}\sum_{h=0}^{H-1}\left[\mathbb{E}_{\tau\sim p(\cdot|\btheta_i)} \left\{ \nabla_{\vvx^j} \log \pi_{\vvx^j}(\va^h_j | \vs^h) \left( \nabla_{\vvx^{\ell}} \log \pi_{\vvx^\ell}(\va^h_{\ell} | \vs^h)\right)^T \right\}\right]_{\alpha, \beta}\\
=& \frac{1}{H}\sum_{h=0}^{H-1}\mathbb{E}_{\tau\sim p(\cdot|\btheta_i)} \left\{ \frac{\partial \log \pi_{\vvx^j}(\va^h_j|\vs^h)}{\partial \vvx_\alpha^j} \frac{\partial \log \pi_{\vvx^\ell}(\va^h_\ell|\vs^h)}{\partial \vvx_\beta^\ell}\right\}\\
=& \frac{1}{H}\sum_{h=0}^{H-1}\int p(\tau|\btheta_i) \frac{\partial \log \pi_{\vvx^j}(\va^h_j|\vs^h)}{\partial \vvx_\alpha^j} \frac{\partial \log \pi_{\vvx^\ell}(\va^h_\ell|\vs^h)}{\partial \vvx_\beta^\ell} d\tau\\
=&\frac{1}{H}\sum_{h=0}^{H-1}\int p(\tau^{-h}) \cdot  \pi_{\btheta_i}(\va^h|\vs^h)\frac{\partial \log \pi_{\vvx^j}(\va^h_j|\vs^h)}{\partial \vvx_\alpha^j} \frac{\partial \log \pi_{\vvx^\ell}(\va^h_\ell|\vs^h)}{\partial \vvx_\beta^\ell}d\tau\\
= &\frac{1}{H}\sum_{h=0}^{H-1}\int p(\tau^{-h}) \prod_{i=1}^n \pi_{\vvx^i}(\va^h_i|\vs^h)\frac{\partial \log \pi_{\vvx^j}(\va^h_j|\vs^h)}{\partial \vvx_\alpha^j} \frac{\partial \log \pi_{\vvx^\ell}(\va^h_\ell|\vs^h)}{\partial \vvx_\beta^\ell}d\tau\\
= &\frac{1}{H}\sum_{h=0}^{H-1}\int p(\tau^{-h}) \prod_{i\neq j,\ell}^n \pi_{\vvx^i}(\va^h_i|\vs^h)\cdot \pi_{\vvx^j}(\va^h_j|\vs^h)\frac{\partial \log \pi_{\vvx^j}(\va^h_j|\vs^h)}{\partial \vvx_\alpha^j} \cdot\pi_{\vvx^\ell}(\va^h_\ell|\vs^h)\frac{\partial \log \pi_{\vvx^\ell}(\va^h_\ell|\vs^h)}{\partial \vvx_\beta^\ell}d\tau\\
=& \frac{1}{H}\sum_{h=0}^{H-1} \int p(\tau^{-h}) \prod_{i\neq j,\ell}^n \pi_{\vvx^i}(\va^h_i|\vs^h)\cdot \frac{\partial  \pi_{\vvx^j}(\va^h_j|\vs^h)}{\partial \vvx_\alpha^j} \cdot \frac{\partial  \pi_{\vvx^\ell}(\va^h_\ell|\vs^h)}{\partial \vvx_\beta^\ell}d\tau\\
=&\frac{1}{H}\sum_{h=0}^{H-1} \frac{\partial^2}{\partial \vvx_{\alpha}^j \partial \vvx_{\beta}^\ell } \int p(\tau)d\tau \\
=&0,
\end{align*}
where $p(\tau^{-h}):= \rho(\vs^0) \prod_{h'\neq h} \pi_{\btheta}(\va^{h'}| \vs^{h'}) P(\vs^{h'+1} | \vs^{h'}, \va^{h'})\cdot P(\vs^{h+1}|\vs^h, \va^h)$. Thus we complete the proof.

\subsection{Proof of Lemma \ref{lemma ascent}}
\label{proof lemma ascent}
Since the objective function $V$ is $L$-smooth, one has 
\begin{align*}
	V(\bbtheta^{t+1}) &\geq V(\bbtheta^t) + \la \nabla V(\bbtheta^t), \bbtheta^{t+1} - \bbtheta^t\ra - \frac{L}{2}\twonorm{\bbtheta^{t+1} - \bbtheta^t}^2\\
	&=V(\bbtheta^t) + \eta \la \nabla V(\bbtheta^t), \bd^t\ra - \frac{L\eta^2}{2}\twonorm{\bd^t}^2, \numberthis\label{eq1}
\end{align*}
where the second line follows from $\bbtheta^{t+1} = \bbtheta^t + \eta \bd^t$. Moreover, for any $i\in [n]$, one has 
\begin{align*}
 \eta \mu_F \twonorm{\vd_i^t}^2 &\stackrel{(a)}{\leq}\eta \la {\mH_i^t }^{-1}\vd_i^t, \vd_i^t\ra \\
	&=	\eta \la \vy_i^{t+1}, \vd_i^t \ra  \\
	&= \eta 	\la \vy_i^{t+1} -\by^{t+1}, \vd_i^t \ra  + \eta \la \by^{t+1}, \vd_i^t \ra\\
	&\leq \eta \twonorm{ \vy_i^{t+1} -\by^{t+1} } \cdot \twonorm{\vd_i^t} + \eta \la \by^{t+1}, \vd_i^t \ra\\
	&\stackrel{(b)}{\leq} \frac{	\eta }{2\mu_F} \twonorm{\vy_i^{t+1} -\by^{t+1}}^2 + \frac{	\eta  \mu_F}{2} \twonorm{ \vd_i^t }^2 + \eta\la \by^{t+1}, \vd_i^t \ra\\
	&= \frac{	\eta }{2\mu_F} \twonorm{\vy_i^{t+1} -\by^{t+1}}^2 + \frac{	\eta  \mu_F}{2} \twonorm{ \vd_i^t }^2 + \eta\la \bv^{t}, \vd_i^t \ra,	\numberthis \label{tmp1}
\end{align*}
where step (a) is due to  \eqref{eig of H}, step (b) uses the elementary inequality that $x\cdot y \leq \frac{1}{2\alpha}x^2 + \frac{\alpha}{2}y^2$ with $\alpha=\mu_F$, and the last line follows from $\by^{t+1} = \bv^t$. Rearranging \eqref{tmp1} yields that 
\begin{align}
	0&\geq -\eta\la \bv^{t}, \vd_i^t \ra+ \frac{	\eta  \mu_F}{2} \twonorm{ \vd_i^t }^2-  \frac{	\eta }{2\mu_F} \twonorm{\vy_i^{t+1} -\by^{t+1}}^2 \label{tmp2}
\end{align}
holds for any fixed $i\in[n]$.  Taking an average over $i$ from $1$ to $n$ yields that
\begin{align*}
	0&\geq -\frac{\eta}{n}\sum_{i=1}^{n}\la \bv^{t}, \vd_i^t \ra+ \frac{	\eta  \mu_F}{2n}\sum_{i=1}^{n} \twonorm{ \vd_i^t }^2 -  \frac{	\eta }{2n\mu_F}\sum_{i=1}^{n} \twonorm{\vy_i^{t+1} -\by^{t+1}}^2\\
	&=-\eta \la \bv^t, \bd^t \ra +  \frac{	\eta  \mu_F}{2n} \twonorm{ \vd^t }^2-  \frac{	\eta }{2n\mu_F} \twonorm{\vy^{t+1} - \bone_n \otimes \by^{t+1}}^2.\numberthis \label{eq2}
\end{align*} 
Summing up \eqref{eq1} and  \eqref{eq2}, we obtain that 
\begin{align*}
	V(\bbtheta^{t+1}) &\geq V(\bbtheta^t) +\eta \la \nabla V(\bbtheta^t) - \bv^t, \bd^t\ra -\frac{\eta}{2n\mu_F} \twonorm{ \vy^{t+1} -\bone_n \otimes \by^{t+1} }^2 + \frac{\eta \mu_F}{2 n} \twonorm{\vd^t}^2- \frac{L\eta^2}{2}\twonorm{\bd^t}^2\\
	&\stackrel{(a)}{\geq} V(\bbtheta^t) -\frac{\eta }{2\gamma } \twonorm{\nabla V(\bbtheta^t) - \bv^t }^2- \frac{\gamma \eta}{2} \twonorm{\bd^t}^2-\frac{\eta}{2n\mu_F} \twonorm{ \vy^{t+1} -\bone_n \otimes \by^{t+1} }^2 + \frac{\eta \mu_F}{2 n} \twonorm{\vd^t}^2- \frac{L\eta^2}{2}\twonorm{\bd^t}^2\\
	&\stackrel{(b)}{\geq} V(\bbtheta^t) -\frac{\eta }{2\gamma} \twonorm{\nabla V(\bbtheta^t) - \bv^t }^2 -\frac{\eta}{2n\mu_F} \twonorm{ \vy^{t+1}-\bone_n \otimes \by^{t+1} }^2 - \frac{\gamma \eta+ L\eta^2}{2n} \twonorm{\vd^t}^2+ \frac{\eta \mu_F}{2 n} \twonorm{\vd^t}^2\\
	&=V(\bbtheta^t) -\frac{\eta }{2\gamma } \twonorm{\nabla V(\bbtheta^t) - \bv^t }^2 -\frac{\eta}{2n\mu_F} \twonorm{ \vy^{t+1} -\bone_n \otimes \by^{t+1} }^2 +\frac{\eta \mu_F - 2\gamma \eta-2 L\eta^2}{4n}\twonorm{\vd^t}^2+ \frac{\eta \mu_F}{4n} \twonorm{\vd^t}^2 \\
	&\stackrel{(c)}{\geq}V(\bbtheta^t) -\frac{4\eta }{\mu_F} \twonorm{\nabla V(\bbtheta^t) - \bv^t }^2 -\frac{\eta}{2n\mu_F} \twonorm{ \vy^{t+1} -\bone_n \otimes \by^{t+1} }^2 + \frac{\eta\mu_F}{8n}\twonorm{\vd^t}^2 + \frac{\eta \mu_F}{4 n} \twonorm{\vd^t}^2, \numberthis \label{eq3}
\end{align*}
where step (a) follows from the elementary inequality that $\la \va,\vb\ra \leq \frac{1}{2\gamma} \twonorm{\va}^2 + \frac{\gamma}{2} \twonorm{\vb}^2$ with $\gamma>0$ for any $\va $ and $\vb$, step (b) is due to $\twonorm{\bd^t}^2 \leq \frac{1}{n}\twonorm{\vd^t}^2$ and step (c) holds by choosing $\gamma = \frac{\mu_F}{8}$ and assuming $0 < \eta \leq \frac{\mu_F}{8L}$ (i.e, $\frac{\eta \mu_F - 2\gamma \eta-2 L\eta^2}{4n} =\frac{\eta\mu_F - \frac{1}{4}\eta \mu_F - 2L\eta^2}{4n}= \frac{\eta}{4n}\left( \frac{3\mu_F}{4} - 2L\eta\right)\geq \frac{\eta\mu_F}{8n}$).  Moreover,  the fact used in step (b) can be proved as follows:
\begin{align*}
	\twonorm{\bd^t}^2 =\twonorm{\frac{1}{n}\sum_{i=1}^{n}\vd_i^t}^2\leq\frac{1}{n}\sum_{i=1}^{n}\twonorm{\vd_i^t}^2	=\frac{1}{n}\twonorm{\vd^t}^2,
\end{align*}
where we have used the Jensen's inequality.

Notice that  
\begin{align*}
	\frac{1}{G^2} \twonorm{ \nabla V(\btheta_i^t)  }^2\leq\twonorm{\mH_i^t \nabla V(\btheta_i^t) }^2 \leq 2\twonorm{\mH_i^t \nabla V(\btheta_i^t) - \vd_i^t}^2 + 2\twonorm{\vd_i^t}^2
\end{align*}
holds for any $i\in [n]$. A direct computation yields that 
\begin{align*}
	\frac{\eta\mu_F }{4n}\twonorm{\vd^t}^2 &=	\frac{ \eta\mu_F }{4n} \sum_{i=1}^{n}\twonorm{\vd_i^t}^2 \\
	&\geq \frac{ \eta\mu_F }{4n} \sum_{i=1}^{n} \left( \frac{1}{2G^2} \twonorm{ \nabla V(\btheta_i^t)  }^2 - \twonorm{\mH_i^t \nabla V(\btheta_i^t) - \vd_i^t}^2  \right)\\
	&=\frac{\eta\mu_F  }{8nG^2}\sum_{i=1}^{n}\twonorm{ \nabla V(\btheta_i^t)  }^2  - \frac{\eta\mu_F  }{4n}\sum_{i=1}^{n}\twonorm{ \mH_i^t  \left(\vy_i^{t+1} - \nabla V(\btheta_i^t)\right)}^2\\
	&\stackrel{(a)}{\geq }\frac{\eta\mu_F  }{8nG^2}\sum_{i=1}^{n} \twonorm{ \nabla V(\btheta_i^t)  }^2  - \frac{\eta\mu_F }{4n}\sum_{i=1}^{n}\frac{1}{\mu_F^2} \twonorm{\vy_i^{t+1} - \nabla V(\btheta_i^t)}^2\\
	&=\frac{\eta\mu_F  }{8nG^2}\sum_{i=1}^{n} \twonorm{ \nabla V(\btheta_i^t)  }^2   - \frac{\eta }{4n\mu_F }\sum_{i=1}^{n} \twonorm{\vy_i^{t+1}-\by^{t+1}+ \by^{t+1} - \nabla V(\bbtheta^t) +\nabla V(\bbtheta^t)  - \nabla V(\btheta_i^t)}^2\\
	&\stackrel{(b)}{\geq }\frac{\eta\mu_F  }{8nG^2}\sum_{i=1}^{n} \twonorm{ \nabla V(\btheta_i^t)  }^2   - \frac{\eta}{4n\mu_F} \sum_{i=1}^{n} \left( 3\left( \twonorm{\vy_i^{t+1} -\by^{t+1}}^2 + \twonorm{\bv^t - \nabla V(\bbtheta^t) }^2+ L^2\twonorm{ \bbtheta^t -  \btheta_i^t}^2 \right) \right)\\
	&=\frac{\eta\mu_F  }{8nG^2}\sum_{i=1}^{n} \twonorm{ \nabla V(\btheta_i^t)  }^2  - \frac{3 \eta }{4n\mu_F}\twonorm{\vy^{t+1} - \bone_n \otimes \by^{t+1}}^2  -  \frac{3 \eta }{4\mu_F } \twonorm{\bv^{t} - \nabla V(\bbtheta^t) }^2 - \frac{3L^2 \eta }{4 n\mu_F } \twonorm{\btheta^t - \bone_n \otimes \bbtheta^t}^2, \numberthis\label{eq4}
\end{align*} 
where step (a) follows from \eqref{eig of H} and step (b) is due to the $L$-smoothness of $V$. Then plugging \eqref{eq4} into \eqref{eq3} yields that
\begin{align*}
	V(\bbtheta^{t+1}) &\geq V(\bbtheta^t) -\frac{4\eta }{\mu_F} \twonorm{\nabla V(\bbtheta^t) - \bv^t }^2 -\frac{\eta}{2n\mu_F} \twonorm{ \vy^{t+1} -\bone_n \otimes \by^{t+1} }^2 + \frac{\eta\mu_F}{8n}\twonorm{\vd^t}^2 + \frac{\eta\mu_F  }{8nG^2}\sum_{i=1}^{n} \twonorm{ \nabla V(\btheta_i^t)  }^2 \\
	&\qquad - \frac{3 \eta }{4n\mu_F}\twonorm{\vy^{t+1} - \bone_n\otimes \by^{t+1}}^2  -  \frac{3 \eta }{4\mu_F } \twonorm{\bv^{t} - \nabla V(\bbtheta^t) }^2 - \frac{3L^2 \eta }{4 n\mu_F } \twonorm{\btheta^t - \bone_n \otimes \bbtheta^t}^2\\
	&= V(\bbtheta^t) -\frac{19\eta }{4\mu_F} \twonorm{\nabla V(\bbtheta^t) - \bv^t }^2+ \frac{\eta\mu_F}{8n}\twonorm{\vd^t}^2 + \frac{\eta\mu_F  }{8nG^2}\sum_{i=1}^{n} \twonorm{ \nabla V(\btheta_i^t)  }^2  \\
	&\qquad - \frac{5 \eta }{4n\mu_F}\twonorm{\vy^{t+1} - \bone_n\otimes \by^{t+1}}^2   - \frac{3L^2 \eta }{4 n\mu_F } \twonorm{\btheta^t - \bone_n \otimes \bbtheta^t}^2. \numberthis \label{eq5}
\end{align*}
Furthermore, it can be seen that
\begin{align*}
	\twonorm{\nabla V(\bbtheta^t)  - \bv^t}^2 &=\twonorm{\nabla V(\bbtheta^t)  - \bnablaf (\btheta^t) + \bnablaf (\btheta^t) - \bv^t}^2 \\
	&\leq 2\twonorm{\nabla V(\bbtheta^t)  - \bnablaf (\btheta^t)}^2 + 2\twonorm{\bnablaf (\btheta^t) - \bv^t}^2 \\
	&\leq\frac{2L^2}{n}\twonorm{\btheta^t - \bone_n\otimes \bbtheta^t}^2 + 2\twonorm{\bnablaf (\btheta^t) - \bv^t}^2, \numberthis\label{eq6}
\end{align*}
where the last line follows from the fact that $\twonorm{\nabla V(\bbtheta^t)  - \bnablaf (\btheta^t)}^2 \leq \frac{L^2}{n}\twonorm{\btheta^t - \bone_n\otimes \bbtheta^t}^2$. Indeed, this fact can be proved as follows:
\begin{align*}
	\twonorm{\nabla V(\bbtheta^t)  - \bnablaf (\btheta^t)}^2 &=\twonorm{\frac{1}{n}\sum_{i=1}^{n} \left( \nabla V_i(\bbtheta^t) - \nabla V_i(\btheta_i^t)\right)}^2\\
	&\leq \frac{1}{n}\sum_{i=1}^{n} \twonorm{ \nabla V_i(\bbtheta^t) - \nabla V_i(\btheta_i^t) }^2\\
	&\leq \frac{L^2}{n}\sum_{i=1}^{n} \twonorm{\btheta_i^t - \bbtheta^t}^2\\
	&= \frac{L^2}{n}\twonorm{\btheta^t - \bone_n\otimes \bbtheta^t}^2,
\end{align*}	
where the second line is due to Jensen's inequality and the third line is due to $L$-smoothness of $V_i$. Thus, plugging \eqref{eq6} into \eqref{eq5} yields that
\begin{align*}
	V(\bbtheta^{t+1}) & \geq V(\bbtheta^t) -\frac{19\eta }{4\mu_F} \left(\frac{2L^2}{n}\twonorm{\btheta^t - \bone_n\otimes \bbtheta^t}^2 + 2\twonorm{\bnablaf (\btheta^t) - \bv^t}^2 \right) + \frac{\eta\mu_F}{8n}\twonorm{\vd^t}^2 + \frac{\eta\mu_F  }{8nG^2}\sum_{i=1}^{n} \twonorm{ \nabla V(\btheta_i^t)  }^2  \\
	&\qquad - \frac{5 \eta }{4n\mu_F}\twonorm{\vy^{t+1} - \bone_n\otimes \by^{t+1}}^2   - \frac{3L^2 \eta }{4 n\mu_F } \twonorm{\btheta^t - \bone_n \otimes \bbtheta^t}^2\\
	&=V(\bbtheta^t) -\frac{19\eta }{2\mu_F}  \twonorm{\bnablaf (\btheta^t) - \bv^t}^2  + \frac{\eta\mu_F}{8n}\twonorm{\vd^t}^2 +\frac{\eta\mu_F  }{8nG^2}\sum_{i=1}^{n} \twonorm{ \nabla V(\btheta_i^t)  }^2  \\
	&\qquad - \frac{5 \eta }{4n\mu_F}\twonorm{\vy^{t+1} - \bone_n\otimes \by^{t+1}}^2   - \frac{41L^2 \eta }{4 n\mu_F } \twonorm{\btheta^t - \bone_n \otimes \bbtheta^t}^2. \numberthis \label{tmp 3}
\end{align*}
Rearranging \eqref{tmp 3} yields that
\begin{align*}
	\frac{1 }{n}\sum_{i=1}^{n} \twonorm{ \nabla V(\btheta_i^t)  }^2 &\leq \frac{8G^2}{\eta\mu_F} \left(  V(\bbtheta^{t+1})  -V(\bbtheta^t) \right) - \frac{G^2}{n}\twonorm{\vd^t}^2  \\
	&\qquad + \frac{76G^2 }{\mu_F^2}  \twonorm{\bnablaf (\btheta^t) - \bv^t}^2  + \frac{10G^2 }{n\mu_F^2}\twonorm{\vy^{t+1} - \bone_n\otimes \by^{t+1}}^2   + \frac{82G^2L^2  }{n\mu_F^2 } \twonorm{\btheta^t - \bone_n \otimes \bbtheta^t}^2. \numberthis \label{tmp 4}
\end{align*}
Taking the telescoping sum of \eqref{tmp 4} over $t$ from $0$ to $T$ for any $T\geq 0$, one has
\begin{align*}
		\frac{1}{n}\sum_{t=0}^{T}\sum_{i=1}^{n} \twonorm{ \nabla V(\btheta_i^t)  }^2   & \leq  \frac{8G^2}{\eta\mu_F}(V(\bx^{T+1}) - V(\bx^0) ) - \frac{G^2}{n}\sum_{t=0}^{T}\twonorm{\vd^t}^2 + \frac{76G^2 }{\mu_F^2}  \sum_{t=0}^{T}\twonorm{\bnablaf (\btheta^t) - \bv^t}^2  \\
		&\qquad + \frac{10G^2 }{n\mu_F^2} \sum_{t=0}^{T} \twonorm{\vy^{t+1} - \bone_n\otimes \by^{t+1}}^2   + \frac{82G^2L^2  }{n\mu_F^2 }\sum_{t=0}^{T} \twonorm{\btheta^t - \bone_n \otimes \bbtheta^t}^2\\
		&\leq  \frac{8G^2}{\eta\mu_F}(V^\ast - V(\bx^0) ) - \frac{G^2}{n}\sum_{t=0}^{T}\twonorm{\vd^t}^2    + \frac{76G^2 }{\mu_F^2}  \sum_{t=0}^{T}\twonorm{\bnablaf (\btheta^t) - \bv^t}^2  \\
		&\qquad+ \frac{10G^2 }{n\mu_F^2} \sum_{t=0}^{T} \twonorm{\vy^{t+1} - \bone_n\otimes \by^{t+1}}^2   + \frac{82G^2L^2  }{n\mu_F^2 }\sum_{t=0}^{T} \twonorm{\btheta^t - \bone_n \otimes \bbtheta^t}^2,
\end{align*}
where the last line has used the Assumption \ref{assumption 0}.

\subsection{Proof of Lemma \ref{errorX}}
\label{proofErrorX}
A sample computation yields that 
\begin{align*}
	\bone_n \otimes \bbtheta^{t+1} = \bone_n \otimes \left( \frac{1}{n} \sum_{i=1}^{n} \btheta_i^{t+1} \right) = \frac{1}{n}\bone_n \otimes \left( \bone_n^\tran \otimes \mI_d  \right) \bbtheta^{t+1}= \frac{1}{n} \left(\mJ_n \otimes \mI_d \right)\btheta^{t+1}.
\end{align*} 
Thus by the update rule described in \eqref{update compact form}, it is straightforward to obtain that 
\begin{align*}
&\twonorm{\btheta^{t+1} - \bone_n \otimes \bbtheta^{t+1}}^2 \\
&= \twonorm{\pW(\btheta^t + \eta \vd^t)  -\frac{1}{n}(\mJ_n\otimes \mI_d)(\btheta^t + \eta \vd^t)}^2\\
	&= \twonorm{ \left( \left(\mW - \frac{1}{n}\mJ_n\right)\otimes \mI_d \right)\btheta^t + \eta\left( \left(\mW - \frac{1}{n}\mJ_n\right)\otimes \mI_d \right)   \vd^t}^2\\
	&\stackrel{(a)}{\leq }\left(1+\frac{1-\rho^2}{2\rho^2} \right) \twonorm{ \left( \left(\mW - \frac{1}{n}\mJ_n\right)\otimes \mI_d \right)\btheta^t }^2 + \eta^2 \left(1+\frac{2\rho^2}{1-\rho^2} \right) \twonorm{\left( \left(\mW - \frac{1}{n}\mJ_n\right)\otimes \mI_d \right)\vd^t}^2\\
	&\stackrel{(b)}{=}\left(1+\frac{1-\rho^2}{2\rho^2} \right) \twonorm{ \left( \left(\mW - \frac{1}{n}\mJ_n\right)\otimes \mI_d \right)(\btheta^t - \bone_n \otimes \bbtheta^t) }^2 \\
	&\qquad + \eta^2 \left(1+\frac{2\rho^2}{1-\rho^2} \right) \twonorm{\left( \left(\mW - \frac{1}{n}\mJ_n\right)\otimes \mI_d \right)(\vd^t-\bone_n \otimes \bd^t)}^2\\
	&\leq \left(1+\frac{1-\rho^2}{2\rho^2} \right)\cdot \opnorm{\mW - \frac{1}{n}\mJ_n}^2\cdot \twonorm{ \btheta^t - \bone_n \otimes \bbtheta^t}^2 + \eta^2 \left(1+\frac{2\rho^2}{1-\rho^2} \right)\cdot \opnorm{\mW - \frac{1}{n}\mJ_n}^2\cdot  \twonorm{\vd^t - \bone_n \otimes \bd^t}^2\\
	&=\left(1+\frac{1-\rho^2}{2\rho^2} \right)\rho^2\cdot \twonorm{ \btheta^t - \bone_n \otimes \bbtheta^t}^2 + \eta^2 \left(1+\frac{2\rho^2}{1-\rho^2} \right)\rho^2 \twonorm{\vd^t - \bone_n \otimes \bd^t}^2\\
	&=\frac{1+\rho^2}{2} \twonorm{  \btheta^t - \bone_n \otimes \bbtheta^t  }^2 + \frac{(1+\rho^2) \rho^2\eta^2}{1-\rho^2}  \twonorm{ \vd^t - \bone_n \otimes \bd^t}^2\\
	&\stackrel{(c)}{\leq } \frac{1+\rho^2}{2} \twonorm{  \btheta^t - \bone_n \otimes \bbtheta^t  }^2 + \frac{2\eta^2}{1-\rho^2}  \twonorm{ \vd^t - \bone_n \otimes \bd^t}^2\\
	&\stackrel{(d)}{\leq} \frac{1+\rho^2}{2} \twonorm{  \btheta^t - \bone_n \otimes \bbtheta^t  }^2 + \frac{2\eta^2}{1-\rho^2} \cdot \left( \frac{G^2}{2\mu_F^2} \twonorm{\vd^t}^2+ \frac{2}{\mu_F^2} \twonorm{\vy^{t+1}- \bone_n \otimes \by^{t+1}}^2 \right)\\
	&=\frac{1+\rho^2}{2} \twonorm{  \btheta^t - \bone_n \otimes \bbtheta^t  }^2 + \frac{4\eta^2}{\mu_F^2(1-\rho^2)}\twonorm{\vy^{t+1}- \bone_n \otimes \by^{t+1}}^2  + \frac{G^2\eta^2}{\mu_F^2(1-\rho^2)}  \twonorm{\vd^t}^2, 
\end{align*}
where step (a) follows from the element inequality that $\twonorm{\va+\vb}^2 \leq (1+\gamma)\twonorm{\va}^2 + (1+\gamma^{-1})\twonorm{\vb}^2$ with $\gamma = \frac{1-\rho^2}{2\rho^2}$, step (b) is due to the fact that $((\mW - n^{-1}\mJ_n ) \otimes \mI_d )(\bone_n \otimes \va) = (\mW\bone_n - n^{-1}\mJ_n \bone_n)\otimes \va = \bzero$ for any $\va\in\R^d$, step (c) is due to $(1+\rho^2)\rho^2\leq 2$, step (d) follows from the fact that 
\begin{align}
	\label{eq d}
	\twonorm{ \vd^t - \bone_n \otimes \bd^t}^2 &\leq  \frac{G^2}{2\mu_F^2} \twonorm{\vd^t}^2+ \frac{2}{\mu_F^2} \twonorm{\vy^{t+1}- \bone_n \otimes \by^{t+1}}^2.
\end{align}
Thus we complete the first part of this lemma. For the second part, using the same argument yields that 
\begin{align*}
\twonorm{\btheta^{t+1} - \bone_n \otimes \bbtheta^{t+1}}^2 &\leq 2\twonorm{ \left( \left(\mW - \frac{1}{n}\mJ_n\right)\otimes \mI_d \right)\btheta^t}^2 + 2\twonorm{ \eta\left( \left(\mW - \frac{1}{n}\mJ_n\right)\otimes \mI_d \right)   \vd^t}^2\\
&\leq 2\rho^2\twonorm{  \btheta^t - \bone_n \otimes \bbtheta^t  }^2 + 2\eta^2\rho^2 \twonorm{ \vd^t - \bone_n \otimes \bd^t}^2 \\
	&\leq 2\rho^2\twonorm{  \btheta^t - \bone_n \otimes \bbtheta^t  }^2 + \frac{4\eta^2\rho^2 }{\mu_F^2} \twonorm{\vy^{t+1}- \bone_n \otimes \by^{t+1}}^2+  \frac{G^2\eta^2\rho^2}{\mu_F^2} \twonorm{\vd^t}^2.
\end{align*}
It only remains to prove the fact \eqref{eq d} used in step (d). Firstly, we have
\begin{align*}
	&\twonorm{ \vd^t - \bone_n \otimes \bd^t} \\
	&= \twonorm{\left(\mI_{nd } - \pJ\right) \vd^t}\\
	&=\twonorm{\left(\mI_{nd } - \pJ\right) (\mH^t \vy^{t+1})}\\
	&=\twonorm{\left(\mI_{nd } - \pJ\right) \left(\mH^t \vy^{t+1} - \left( \frac{1}{2\mu_F} + \frac{1}{2G} \right)\vy^{t+1} + \left( \frac{1}{2\mu_F} + \frac{1}{2G} \right)\vy^{t+1} \right)}\\
	&\leq \twonorm{\left(\mI_{nd } - \pJ\right) \left(\mH^t \vy^{t+1} - \left( \frac{1}{2\mu_F} + \frac{1}{2G} \right) \vy^{t+1} \right)}  + \left( \frac{1}{2\mu_F} + \frac{1}{2G} \right)\twonorm{\left(\mI_{nd } - \pJ\right)\vy^{t+1}}\\
	&\leq \opnorm{  \mI_{nd } - \pJ  } \cdot \opnorm{\mH^t - \left( \frac{1}{2\mu_F} + \frac{1}{2G} \right)\mI_{nd }}\cdot \twonorm{  \vy^{t+1}}  + \left( \frac{1}{2\mu_F} + \frac{1}{2G} \right)\twonorm{\left(\mI_{nd } - \pJ\right)\vy^{t+1}}\\
	&=\opnorm{  \mI_{nd } - \pJ  }\cdot \opnorm{\mH^t - \left( \frac{1}{2\mu_F} + \frac{1}{2G} \right)\mI_{nd }}\cdot \twonorm{  \vy^{t+1} }  + \left( \frac{1}{2\mu_F} + \frac{1}{2G} \right)\twonorm{\vy^{t+1} - \bone_n \otimes \by^{t+1}}\\
	&\stackrel{(a)}{\leq}\frac{1}{2}\left(\frac{1}{\mu_F} -\frac{1}{G}\right)\twonorm{  \vy^{t+1} }  + \frac{1}{2}\left(\frac{1}{\mu_F} + \frac{1}{G}\right)\twonorm{\vy^{t+1}- \bone_n \otimes \by^{t+1}},
\end{align*}
where the last line follows from that $\opnorm{  \mI_{nd } - \pJ  } = \opnorm{(\mI_d - n^{-1}\mI_d)\otimes \mI_n} \leq 1-\frac{1}{n}\leq 1$ and 
\begin{align*}
\opnorm{\mH^t - \left( \frac{1}{2\mu_F} + \frac{1}{2G} \right) \mI_{nd}}&\leq \frac{1}{2}\left(\frac{1}{\mu_F} -\frac{1}{G}\right).
\end{align*} 
Then a direct computation yields that 
\begin{align*}
	\twonorm{ \vd^t - \bone_n \otimes \bd^t} ^2 &\leq \frac{1}{2}\left(\frac{1}{\mu_F} -\frac{1}{G}\right)^2\twonorm{  \vy^{t+1} }^2  + \frac{1}{2}\left(\frac{1}{\mu_F} + \frac{1}{G}\right)^2\twonorm{\vy^{t+1}- \bone_n \otimes \by^{t+1}}^2\\
	&\stackrel{(a)}{\leq }\frac{1}{2}\left(\frac{1}{\mu_F} -\frac{1}{G}\right)^2 G^2\twonorm{\vd^t}^2+ \frac{1}{2}\left(\frac{1}{\mu_F} + \frac{1}{G}\right)^2\twonorm{\vy^{t+1}- \bone_n \otimes \by^{t+1}}^2\\
	&=\frac{1}{2}\left(\frac{G}{\mu_F} -1\right)^2 \twonorm{\vd^t}^2+ \frac{1}{2}\left(\frac{1}{\mu_F} + \frac{1}{G}\right)^2\twonorm{\vy^{t+1}- \bone_n \otimes \by^{t+1}}^2\\
	&\stackrel{(b)}{\leq } \frac{G^2}{2\mu_F^2}  \twonorm{\vd^t}^2+ \frac{1}{2}\left(\frac{1}{\mu_F} + \frac{1}{G}\right)^2\twonorm{\vy^{t+1}- \bone_n \otimes \by^{t+1}}^2\\
	&\leq \frac{G^2}{2\mu_F^2} \twonorm{\vd^t}^2+ \frac{2}{\mu_F^2} \twonorm{\vy^{t+1}- \bone_n \otimes \by^{t+1}}^2,
\end{align*}
where step (a) is due to the fact $\twonorm{\vy^{t+1}}^2 =\sum_{i=1}^{n}\twonorm{\vy_i^{t+1}}^2 = \sum_{i=1}^{n}\twonorm{ {\mH_i^{t}}^{-1}\vd_i^t}^2\leq G^2\twonorm{\vd^t}^2$, and step (b) is due to $ \frac{G}{\mu_F}\geq 1$.

\subsection{Proof of Lemma \ref{lemma tracking}}
\label{proofTracking}

\subsubsection{Proof of \eqref{eq Y1} }
Recall the initialization in Algorithm \ref{alg: moentum pg} that $\vy_i^{0}= \bzero, \vv_i^{-1} = \bzero $, and $\vv_i^0 = \frac{1}{B}\sum_{b=1}^{B} \vg_i( \tau_{i,b}^0 | \btheta_i^0 )$. A direct computation yields that 
\begin{align*}
\E{\twonorm{ \vy^1 - \bone_n  \otimes \by^1}}^2 &= \E{ \twonorm{ \pW \vv^0 - \bone_n \otimes \bv^0 }^2}\\
&= \E{ \twonorm{ \pW \vv^0 - \frac{1}{n}\left(\mJ_n \otimes \mI_d\right) \vv^0}^2}\\
&\leq  \opnorm{\mW - \frac{1}{n}\mJ_n}^2\cdot \E{\twonorm{\vv^0  }^2}\\
&=\rho^2 \sum_{i=1}^{n} \E{\twonorm{\vv_i^0 - \nabla V_i(\bbtheta^0) + \nabla V_i(\bbtheta^0) }^2}\\
&= \rho^2\sum_{i=1}^{n} \E{\twonorm{\vv_i^0 - \nabla V_i(\bbtheta^0) }^2 }+ \rho^2\sum_{i=1}^{n} \twonorm{\nabla V_i(\bbtheta^0) }^2 +2 \rho^2 \sum_{i=1}^n \la\E{\vv_i^0} - \nabla V_i(\bbtheta^0), \nabla V_i(\bbtheta^0) \ra\\
&= \rho^2\sum_{i=1}^{n} \E{\twonorm{\vv_i^0 - \nabla V_i(\bbtheta^0) }^2 }+ \rho^2\sum_{i=1}^{n} \twonorm{\nabla V_i(\bbtheta^0) }^2,\numberthis \label{eq43a}
\end{align*}
where the last line follows from $\E{\vv_i^0}=\nabla V_i(\btheta^0_i) =\nabla V_i(\bbtheta^0)$. Moreover, for any $i\in[n]$, it can be seen that 
\begin{align*}
	&\E{\twonorm{\vv_i^0 - \nabla V_i(\bbtheta^0) }^2 } \\
	&= \E{\twonorm{ \frac{1}{B}\sum_{b=1}^{B} \vg_i( \tau_{i,b}^0 | \bbtheta^0 )- \nabla V_i(\bbtheta^0) }^2 }\\
	&= \E{\twonorm{ \frac{1}{B}\sum_{b=1}^{B}  \left( \vg_i( \tau_{i,b}^0 | \bbtheta^0 )- \nabla V_i(\bbtheta^0) \right)}^2 }\\
	&=\frac{1}{B^2} \sum_{b=1}^{B}\E{ \twonorm{ \vg_i( \tau_{i,b}^0 | \bbtheta^0 )- \nabla V_i(\bbtheta^0)}^2} + \frac{1}{B^2} \sum_{b\neq b'}\E{\la \vg_i( \tau_{i,b}^0 | \bbtheta^0 )- \nabla V_i(\bbtheta^0), \vg_i( \tau_{i,b'}^0 | \bbtheta^0 )- \nabla V_i(\bbtheta^0)\ra}\\
	&\stackrel{(a)}{=}\frac{1}{B^2} \sum_{b=1}^{B}\E{ \twonorm{ \vg_i( \tau_{i,b}^0 | \bbtheta^0 )- \nabla V_i(\bbtheta^0)}^2}\\
	&\stackrel{(b)}{\leq } \frac{\nu_i^2}{B}, \numberthis \label{eq 43b} 
\end{align*}
where step (a) is due to the fact that $\{\tau_{i,b}^0\}_{b=1}^B$ are independent trajectories and step (b) follows from Assumption \ref{assumption 3}.  Substituting \eqref{eq 43b} into \eqref{eq43a} yields that
\begin{align*}
	\E{\twonorm{\vy^1 - \bone_n \otimes \by^1}^2} &\leq \frac{\rho^2}{B}\sum_{i=1}^{n}\nu_i^2 + \rho^2\sum_{i=1}^{n} \twonorm{\nabla V_i(\bbtheta^0) }^2\\
	&=\frac{n\rho^2\bar{\nu}^2}{B} + \rho^2\sum_{i=1}^{n} \twonorm{\nabla V_i(\bbtheta^0) }^2,
\end{align*}
which completes the proof of \eqref{eq Y1}.

\subsubsection{Proof of \eqref{Error Y}}
Following the gradient tracking update in \eqref{update compact form}, we have
\begin{align*}
	&\E{\twonorm{\vy^{t+1} - \bone_n \otimes \by^{t+1}}^2 }\\
	&=\E{\twonorm{\pW(\vy^t + \vv^t - \vv^{t-1}) - \frac{1}{n}\left(\bone_n\bone_n^\tran \otimes\mI_d \right)\pW(\vy^t + \vv^t - \vv^{t-1}) }^2}\\
	&= \E{ \twonorm{  \left(\left( \mW - \frac{1}{n}\mJ_n \right) \otimes \mI_d \right)(\vy^t + \vv^t - \vv^{t-1})  }^2}\\
	&\leq \left( 1+ \frac{1-\rho^2}{2\rho^2}\right) \E{\twonorm{\left(\left(\mW - \frac{1}{n}\mJ_n \right)\otimes \mI_d \right) \vy^t }^2 }+ \left( 1+ \frac{2\rho^2}{1-\rho^2}\right)\E{\twonorm{\left(\left(\mW - \frac{1}{n}\mJ_n \right)\otimes \mI_d \right) (\vv^t - \vv^{t-1})}^2 }\\
	&= \frac{1+\rho^2}{2\rho^2}  \E{\twonorm{\left(\left(\mW - \frac{1}{n}\mJ_n \right)\otimes \mI_d \right) (\vy^t - \bone_n \otimes \by^t) }^2 }+   \frac{1+\rho^2}{1-\rho^2} \E{\twonorm{\left(\left(\mW - \frac{1}{n}\mJ_n \right)\otimes \mI_d \right) (\vv^t - \vv^{t-1})}^2 }\\
	&\leq \frac{1+\rho^2}{2} \E{\twonorm{ \vy^t -\bone_n \otimes \by^t}^2 }+ \frac{(1+\rho^2)\rho^2}{1-\rho^2} \E{\twonorm{  \vv^t - \vv^{t-1}}^2},\numberthis \label{eq Y2}
\end{align*}
where the third line is due to the element inequality that $\twonorm{\va+\vb}^2 \leq (1+c)\twonorm{\va}^2 + (1+c^{-1})\twonorm{\vb}^2$ with $c=\frac{1-\rho^2}{2\rho^2}$ for any $\va$ and $\vb$. Moreover, 
we have the following relationship:
\begin{align}
	\label{fact 1}
	\E{\twonorm{  \vv^t - \vv^{t-1}}^2} &\leq  \left(8(1-\beta)^2L_g^2 + 8(1-\beta)^2C_g^2 C_{\omega}^2 + 4\beta^2 L_g^2 \right)\E{\twonorm{ \btheta^t - \btheta^{t-1}}^2 } + 4n\beta^2 \bar{ \nu}^2 \notag \\
	&\qquad + 4\beta^2\sum_{i=1}^{n} \E{\twonorm{  \nabla V_i(\btheta_i^{t-1})  - \vv_i^{t-1 } }^2},
\end{align}
which has been shown in (61) of \cite{jiang2021mdpgt}. Substituting \eqref{fact 1} into \eqref{eq Y2} yields that
\begin{align*}
	&\E{\twonorm{\vy^{t+1} - \bone_n \otimes \by^{t+1}}^2 }\\
	&\leq\frac{1+\rho^2}{2} \E{\twonorm{ \vy^t -\bone_n \otimes \by^t}^2 }  + \frac{(1+\rho^2)\rho^2}{1-\rho^2} \left( 12\Phi^2\E{\twonorm{ \btheta^t - \btheta^{t-1}}^2 } + 4n\beta^2 \bar{ \nu}^2+ 4\beta^2\sum_{i=1}^{n} \E{\twonorm{  \nabla V_i(\btheta_i^{t-1})  - \vv_i^{t-1 } }^2}\right)\\
	&= \frac{1+\rho^2}{2} \E{ \twonorm{ \vy^t -\bone_n \otimes \by^t}^2 } +\frac{4n\beta^2 \bar{\nu}^2\rho^2(1+\rho^2)}{1-\rho^2}  \\
	&\qquad + \frac{12\Phi^2(1+\rho^2)\rho^2}{1-\rho^2}  \E{\twonorm{\btheta^t - \btheta^{t-1}}^2}+ \frac{4\beta^2\rho^2(1+\rho^2)}{1-\rho^2}\sum_{i=1}^{n} \E{ \twonorm{\nabla V_i(\btheta_i^{t-1} ) - \vv_i^{t-1}}^2}\\
	&\leq \frac{1+\rho^2}{2} \E{\twonorm{ \vy^t -\bone_n \otimes \by^t}^2 }+ \frac{8n\beta^2 \bar{\nu}^2 }{1-\rho^2} \\
	&\qquad  + \frac{24\Phi^2 }{1-\rho^2}  \E{\twonorm{\btheta^t - \btheta^{t-1}}^2} + \frac{8\beta^2 }{1-\rho^2} \sum_{i=1}^{n} \E{ \twonorm{\nabla V_i(\btheta_i^{t-1} ) - \vv_i^{t-1}}^2}, \numberthis \label{eq Y3}
\end{align*}   
where the first inequality follows from that $8(1-\beta)^2L_g^2 + 8(1-\beta)^2C_g^2 C_{\omega}^2 + 4\beta^2 L_g^2 \leq 12 (L_g^2 + C_g^2 C_{\omega}^2):=12\Phi^2$ for $0\leq \beta\leq 1$ and the last line is due to $\rho <1$. Furthermore, the term $\E{\twonorm{\btheta^t - \btheta^{t-1}}^2} $ can be bounded as follows:
\begin{align*}
	&\E{\twonorm{\btheta^t - \btheta^{t-1}}^2} \\
	&= \E{\twonorm{ \btheta^t -\bone_n \otimes \bbtheta^t + \bone_n \otimes \bbtheta^t - \bone_n \otimes \bbtheta^{t-1} +\bone_n \otimes \bbtheta^{t-1} - \btheta^{t-1} }^2}\\
	&\leq 3\E{\twonorm{\btheta^t - \bone_n \otimes \bx^t}^2} + 3\E{\twonorm{\btheta^{t-1} -\bone_n \otimes \bx^{t-1}}^2} + 3\E{\twonorm{\bone_n \otimes(\bx^t - \bx^{t-1})}^2} \\
	& = 3\E{\twonorm{\btheta^t - \bone_n \otimes \bx^t}^2} + 3\E{\twonorm{\btheta^{t-1} -\bone_n \otimes \bx^{t-1}}^2} + 3n\eta^2\E{\twonorm{\bd^{t-1}}^2}\\
	&\stackrel{(a)}{\leq} 3\left( 2\rho^2 \E{\twonorm{  \btheta^{t-1} - \bone_n \otimes \bbtheta^{t-1}  }^2 }+ \frac{4\eta^2\rho^2 }{\mu_F^2} \E{\twonorm{\vy^{t}- \bone_n \otimes \by^{t}}^2 }+  \frac{G^2 \eta^2\rho^2}{\mu_F^2} \E{\twonorm{\vd^{t-1}}^2}\right)\\
	&\qquad + 3\E{\twonorm{\btheta^{t-1} -\bone_n \otimes \bx^{t-1}}^2} + 3n\eta^2\E{\twonorm{\bd^{t-1}}^2}\\
	&\stackrel{(b)}{\leq} 3\left( 2\rho^2 \E{ \twonorm{  \btheta^{t-1} - \bone_n \otimes \bbtheta^{t-1}  }^2} + \frac{4\eta^2\rho^2 }{\mu_F^2} \E{ \twonorm{\vy^{t}- \bone_n \otimes \by^{t}}^2} +  \frac{G^2\eta^2\rho^2}{\mu_F^2} \E{\twonorm{\vd^{t-1}}^2}\right)\\
	&\qquad + 3\E{\twonorm{\btheta^{t-1} -\bone_n \otimes \bx^{t-1}}^2} + 3\eta^2\E{\twonorm{\vd^{t-1}}^2}\\
	& \stackrel{(c)}{\leq } 9 \E{\twonorm{  \btheta^{t-1}- \bone_n \otimes \bbtheta^{t-1}}^2} + \frac{12\eta^2\rho^2}{\mu_F^2} \E{\twonorm{\vy^{t} - \bone_n \otimes \by^{t}}^2}+ \frac{6G^2\eta^2}{\mu_F^2}\E{\twonorm{\vd^{t-1}}^2}, \numberthis \label{diff theta}
\end{align*}
where step (a) is due to \eqref{eqErrorX2}, step (b) follows from the fact that $\twonorm{\bd^{t-1}}^2 = \twonorm{\frac{1}{n}\sum_{i=1}^n \vd_i^{t-1}}^2 \leq \frac{1}{n}\sum_{i=1}^n \twonorm{\vd_i^{t-1}}^2 = \frac{1}{n} \twonorm{\vd^{t-1}}^2$, and step (c) holds since $\rho <1$ and $\mu_F \leq G$. Substituting \eqref{diff theta} into \eqref{eq Y3} yields that 
\begin{align*}
&\E{\twonorm{\vy^{t+1} - \bone_n \otimes \by^{t+1}}^2 } \\
&\leq \frac{1+\rho^2}{2} \E{\twonorm{ \vy^t -\bone_n \otimes \by^t}^2 }+ \frac{8n\beta^2 \bar{\nu}^2 }{1-\rho^2} + \frac{8\beta^2 }{1-\rho^2} \sum_{i=1}^{n} \E{ \twonorm{\nabla V_i(\btheta_i^{t-1} ) - \vv_i^{t-1}}^2}\\
	&\qquad  + \frac{24\Phi^2}{1-\rho^2} \left( 9 \E{\twonorm{  \btheta^{t-1}- \bone_n \otimes \bbtheta^{t-1}}^2} + \frac{12\eta^2\rho^2}{\mu_F^2} \E{\twonorm{\vy^{t} - \bone_n \otimes \by^{t}}^2}+ \frac{6G^2\eta^2}{\mu_F^2}\E{\twonorm{\vd^{t-1}}^2} \right)\\
	&=\left( \frac{1+\rho^2}{2} +\frac{24\Phi^2}{1-\rho^2} \cdot \frac{12\eta^2\rho^2}{\mu_F^2} \right) \E{\twonorm{ \vy^t -\bone_n \otimes \by^t}^2} + \frac{8n\beta^2 \bar{\nu}^2 }{1-\rho^2} + \frac{8\beta^2 }{1-\rho^2} \sum_{i=1}^{n}\E{  \twonorm{\nabla V_i(\btheta_i^{t-1}) - \vv_i^{t-1}}^2} \\
	&\qquad +\frac{216\Phi^2}{1-\rho^2}  \E{ \twonorm{  \btheta^{t-1}- \bone_n \otimes \bbtheta^{t-1}}^2} + \frac{144\Phi^2 G^2\eta^2}{\mu_F^2 (1-\rho^2)}  \E{\twonorm{\vd^{t-1}}^2} \\
	&\leq \frac{3+\rho^2}{4}  \E{\twonorm{ \vy^t -\bone_n \otimes \by^t}^2 }+ \frac{8n\beta^2 \bar{\nu}^2 }{1-\rho^2} + \frac{8\beta^2 }{1-\rho^2}  \sum_{i=1}^{n}\E{  \twonorm{\nabla V_i(\btheta_i^{t-1}) - \vv_i^{t-1}}^2}\\
	&\qquad +\frac{216\Phi^2}{1-\rho^2}  \E{ \twonorm{  \btheta^{t-1}- \bone_n \otimes \bbtheta^{t-1}}^2 } + \frac{144\Phi^2G^2\eta^2}{\mu_F^2 (1-\rho^2)}  \E{\twonorm{\vd^{t-1}}^2}\\
	&=\frac{3+\rho^2}{4}  \E{\twonorm{ \vy^t -\bone_n \otimes \by^t}^2 }+ \frac{8n\beta^2 \bar{\nu}^2 }{1-\rho^2} + \frac{8\beta^2 }{1-\rho^2}  \E{\twonorm{ \tnablaV(\btheta^{t-1}) - \vv^{t-1}}^2}\\
	&\qquad +\frac{216\Phi^2}{1-\rho^2}  \E{ \twonorm{  \btheta^{t-1}- \bone_n \otimes \bbtheta^{t-1}}^2 } + \frac{144\Phi^2G^2\eta^2}{\mu_F^2 (1-\rho^2)}  \E{\twonorm{\vd^{t-1}}^2}
\end{align*}
where the second inequality is due to $\eta \leq \frac{\mu_F(1-\rho^2)}{24\sqrt{2}\Phi}$, i.e., 
\begin{align*}
	\frac{1+\rho^2}{2} +\frac{24\Phi^2}{1-\rho^2} \cdot \frac{12\eta^2\rho^2}{\mu_F^2} \leq \frac{3+\rho^2}{4}.
\end{align*}
Now the proof is complete.

\subsection{Proof of Lemma \ref{accumulated Error Y}}
\label{proof lemma accumulated error Y}
Applying \eqref{lemma6b} to \eqref{Error Y} yields that 
	\begin{align*}
		&\sum_{t=1}^{T} \E{\twonorm{ \vy^t - \bone_n \otimes \by^t}^2 } \\
		&\leq \frac{4}{1-\rho^2} \E{ \twonorm{ \vy^1 - \bone_n \otimes \by^1}^2 } +  \frac{4T}{1-\rho^2} \cdot \frac{8n\beta^2 \bar{\nu}^2 }{1-\rho^2} \\
		&\qquad +\frac{4}{1-\rho^2} \sum_{t=0}^{T-2}\bigg( \frac{8\beta^2 }{1-\rho^2}  \E{\twonorm{ \tnablaV(\btheta^{t}) - \vv^{t}}^2}   +\frac{216\Phi^2}{1-\rho^2}  \E{ \twonorm{  \btheta^{t}- \bone_n \otimes \bbtheta^{t}}^2 }  \\
		&\qquad + \frac{144\Phi^2G^2\eta^2}{\mu_F^2 (1-\rho^2)}  \E{\twonorm{\vd^{t}}^2} \bigg) \\
		&\leq \frac{4}{1-\rho^2} \E{ \twonorm{ \vy^1 - \bone_n \otimes \by^1}^2 } +  \frac{32nT\beta^2 \bar{\nu}^2 }{(1-\rho^2)^2} + \frac{576\Phi^2G^2\eta^2 }{\mu_F^2(1-\rho^2)^2} \sum_{t=0}^{T} \E{\twonorm{\vd^{t}}^2} \\
		&\qquad + \frac{864\Phi^2 }{(1-\rho^2)^2} \sum_{t=0}^{T}  \E{\twonorm{  \btheta^{t}- \bone_n \otimes \bbtheta^{t}}^2 } \\
		&\qquad +  \frac{32\beta^2 }{(1-\rho^2)^2} \sum_{t=0}^{T}  \E{\twonorm{ \tnablaV(\btheta^{t}) - \vv^{t}}^2}   \\
		&\stackrel{(a)}{\leq}\frac{4}{1-\rho^2} \E{ \twonorm{ \vy^1 - \bone_n \otimes \by^1}^2 } +  \frac{32nT\beta^2 \bar{\nu}^2 }{(1-\rho^2)^2} + \frac{576\Phi^2G^2\eta^2 }{\mu_F^2(1-\rho^2)^2} \sum_{t=0}^{T} \E{\twonorm{\vd^{t}}^2} \\
		&\qquad + \frac{864\Phi^2 }{(1-\rho^2)^2} \sum_{t=0}^{T}  \E{\twonorm{  \btheta^{t}- \bone_n \otimes \bbtheta^{t}}^2 } \\
		&\qquad +  \frac{32\beta^2 }{(1-\rho^2)^2}\left(\frac{n\bar{\nu}^2}{\beta B} + 2n \beta T \bar{\nu}^2 + \frac{12n  \eta^2\Phi^2}{\beta}\sum_{t=0}^{T-1} \E{\twonorm{\bd^t}^2} + \frac{24\Phi^2}{\beta} \sum_{t=0}^{T} \E{\twonorm{\vx^t - \bone_n \otimes \bx^t}^2}\right)\\
		&\stackrel{(b)}{\leq } \frac{4}{1-\rho^2} \left(\frac{n\rho^2\bar{\nu}^2}{B} + \rho^2\twonorm{\tnablaV(\btheta^0)}^2 \right) +  \frac{32nT\beta^2 \bar{\nu}^2 }{(1-\rho^2)^2} + \frac{576\Phi^2G^2\eta^2 }{\mu_F^2(1-\rho^2)^2} \sum_{t=0}^{T} \E{\twonorm{\vd^{t}}^2} \\
		&\qquad + \frac{864\Phi^2 }{(1-\rho^2)^2} \sum_{t=0}^{T}  \E{\twonorm{  \btheta^{t}- \bone_n \otimes \bbtheta^{t}}^2 } \\
		&\qquad +  \frac{32\beta^2 }{(1-\rho^2)^2}\left(\frac{n\bar{\nu}^2}{\beta B} + 2n \beta T \bar{\nu}^2 + \frac{12n  \eta^2\Phi^2}{\beta}\sum_{t=0}^{T-1} \E{\twonorm{\bd^t}^2} + \frac{24\Phi^2}{\beta} \sum_{t=0}^{T} \E{\twonorm{\vx^t - \bone_n \otimes \bx^t}^2}\right)\\
		&= \frac{4}{1-\rho^2} \bigg(\frac{n\rho^2\bar{\nu}^2}{B} + \rho^2\twonorm{\tnablaV(\btheta^0)}^2\bigg) +  \frac{32n T\beta^2 \bar{\nu}^2 }{(1-\rho^2)^2} + \frac{32\beta^2 }{(1-\rho^2)^2} \bigg( \frac{n\bar{\nu}^2}{\beta B} + 2n \beta T \bar{\nu}^2  \bigg) \\
		&\qquad + \left( \frac{864\Phi^2 }{(1-\rho^2)^2}+\frac{32\beta^2 }{(1-\rho^2)^2}\cdot   \frac{24\Phi^2 }{\beta}\right) \sum_{t=0}^{T}  \E{\twonorm{  \btheta^{t}- \bone_n \otimes \bbtheta^{t}}^2 } \\
		&\qquad + \frac{576\Phi^2G^2\eta^2 }{\mu_F^2(1-\rho^2)^2} \sum_{t=0}^{T} \E{\twonorm{\vd^{t}}^2} +\frac{32\beta^2 }{(1-\rho^2)^2}\cdot \frac{12n  \eta^2\Phi^2}{\beta}\sum_{t=0}^{T-1} \E{\twonorm{\bd^t}^2}  \\
		&\stackrel{(c)}{\leq} \frac{4}{1-\rho^2} \bigg(\frac{n\rho^2\bar{\nu}^2}{B} + \rho^2\twonorm{\tnablaV(\btheta^0)}^2\bigg) +  \frac{32n T\beta^2 \bar{\nu}^2 }{(1-\rho^2)^2} + \frac{32\beta^2 }{(1-\rho^2)^2} \bigg( \frac{n\bar{\nu}^2}{\beta B} + 2n \beta T \bar{\nu}^2  \bigg) \\
		&\qquad + \left( \frac{864\Phi^2 }{(1-\rho^2)^2}+\frac{32\beta^2 }{(1-\rho^2)^2}\cdot   \frac{24\Phi^2 }{\beta}\right) \sum_{t=0}^{T}  \E{\twonorm{  \btheta^{t}- \bone_n \otimes \bbtheta^{t}}^2 } \\
		&\qquad + \left(\frac{576\Phi^2G^2\eta^2 }{\mu_F^2(1-\rho^2)^2}   + \frac{384\beta \eta^2\Phi^2}{(1-\rho^2)^2} \right)\sum_{t=0}^{T-1} \E{\twonorm{\vd^t}^2} \\
		&\leq  \frac{4}{1-\rho^2} \bigg(\frac{n\rho^2\bar{\nu}^2}{B} + \rho^2\twonorm{\tnablaV(\btheta^0)}^2\bigg) +  \frac{32n T\beta^2 \bar{\nu}^2 }{(1-\rho^2)^2} + \frac{32\beta^2 }{(1-\rho^2)^2} \bigg( \frac{n\bar{\nu}^2}{\beta B} + 2n \beta T \bar{\nu}^2  \bigg) \\
		&\qquad + \frac{1632\Phi^2 }{(1-\rho^2)^2}\sum_{t=0}^{T}  \E{\twonorm{  \btheta^{t}- \bone_n \otimes \bbtheta^{t}}^2 }  +  \frac{960\Phi^2G^2\eta^2 }{\mu_F^2(1-\rho^2)^2}    \sum_{t=0}^{T-1} \E{\twonorm{\vd^t}^2} 
	\end{align*}
	where step (a) is due to \eqref{eq v2}, step (b) follows from \eqref{eq Y1}, step (c) holds since $\twonorm{\bd^t}^2 \leq \frac{1}{n}\twonorm{\vd^t}^2$, and the last inequality is due to $\beta < 1$, i.e.,
	\begin{align*}
		\frac{864\Phi^2 }{(1-\rho^2)^2}+\frac{32\beta^2 }{(1-\rho^2)^2}\cdot   \frac{24\Phi^2 }{\beta} &=	\frac{864\Phi^2 }{(1-\rho^2)^2}+\frac{768\Phi^2 \beta  }{(1-\rho^2)^2} \leq 	\frac{1632\Phi^2 }{(1-\rho^2)^2}.
	\end{align*}
 Thus we complete the proof.
 
\subsection{Proof of Lemma \ref{accumulatedErrorX}}
\label{proof lemma accumulated X}
Due to \eqref{eqErrorX}, it can be seen that 
\begin{align*}
	&\E{\twonorm{\btheta^{t} - \bone_n \otimes \bbtheta^{t}}^2} \\
	&\leq \frac{1+\rho^2}{2} \E{\twonorm{  \btheta^{t-1} - \bone_n \otimes \bbtheta^{t-1} }^2} + \frac{4\eta^2}{\mu_F^2(1-\rho^2)} \E{\twonorm{\vy^{t}- \bone_n \otimes \by^{t}}^2 } +   \frac{G^2\eta^2}{\mu_F^2(1-\rho^2)}  \E{\twonorm{\vd^{t-1}}^2}\\
	&=\frac{1+\rho^2}{2} \E{\twonorm{  \btheta^{t-1} - \bone_n \otimes \bbtheta^{t-1} }^2 }+   \frac{1+\rho^2}{2}  \cdot \frac{2}{1+\rho^2}\cdot \frac{\eta^2\kappa_{F} ^2}{1-\rho^2} \E{ \twonorm{\vd^{t-1}}^2} + \frac{4\eta^2}{\mu_F^2(1-\rho^2)}\E{ \twonorm{\vy^{t}- \bone_n \otimes \by^{t}}^2  }\\
	&\leq \frac{1+\rho^2}{2}\E{ \twonorm{  \btheta^{t-1} - \bone_n \otimes \bbtheta^{t-1} }^2} +   \frac{1+\rho^2}{2}  \cdot \frac{2}{1-\rho^2}\cdot \frac{\eta^2\kappa_{F} ^2}{1-\rho^2} \E{ \twonorm{\vd^{t-1}}^2 } + \frac{4\eta^2}{\mu_F^2(1-\rho^2)} \E{\twonorm{\vy^{t}- \bone_n \otimes \by^{t}}^2 }\\
	&= \frac{1+\rho^2}{2} \E{\twonorm{  \btheta^{t-1} - \bone_n \otimes \bbtheta^{t-1} }^2 } +   \frac{1+\rho^2}{2}  \cdot  \frac{2\eta^2\kappa_{F} ^2}{(1-\rho^2)^2}  \E{\twonorm{\vd^{t-1}}^2 }+ \frac{4\eta^2}{\mu_F^2(1-\rho^2)} \E{\twonorm{\vy^{t}- \bone_n \otimes \by^{t}}^2} \numberthis \label{eq tmp1}
\end{align*}
for any $t\geq 0$. Applying \eqref{lemma6a} to \eqref{eq tmp1} leads to that 
\begin{align*}
	&\sum_{t=0}^{T} \E{\twonorm{\btheta^{t} - \bone_n \otimes \bbtheta^{t}}^2}  \\
	&\leq   \frac{4G^2\eta^2}{\mu_F^2(1-\rho^2)^3}  \sum_{t=0}^{T}\E{\twonorm{\vd^t}^2} + \frac{8\eta^2}{\mu_F^2(1-\rho^2)^2}    \sum_{t=1}^{T}\twonorm{\vy^{t} - \bone_n \otimes \by^{t}}^2\\
	&\leq \frac{4G^2\eta^2}{\mu_F^2(1-\rho^2)^3}  \sum_{t=0}^{T}\E{\twonorm{\vd^t}^2} \\
	&\qquad + \frac{8\eta^2}{\mu_F^2(1-\rho^2)^2} \left( A_1 + \frac{1632\Phi^2 }{(1-\rho^2)^2}\sum_{t=0}^{T}  \E{\twonorm{  \btheta^{t}- \bone_n \otimes \bbtheta^{t}}^2 }  +  \frac{960\Phi^2G^2\eta^2 }{\mu_F^2(1-\rho^2)^2}    \sum_{t=0}^{T-1} \E{\twonorm{\vd^t}^2}  \right),
\end{align*}
where the first inequality has used the fact that $\btheta^0_i = \bbtheta^0$ for all $i\in [n]$ and the second inequality follows from \eqref{error accumulated Y}. Since
\begin{align*}
	0 < \eta <\frac{\mu_F(1-\rho^2)^3}{ \kappa_{F}\sqrt{1632000(L^2 + \Phi^2) }},
\end{align*}
it can be seen that 
\begin{align*}
\frac{8\eta^2}{\mu_F^2(1-\rho^2)^2} \cdot\frac{1632\Phi^2 }{(1-\rho^2)^2} &\leq   \frac{8\cdot 1632\Phi^2 }{\mu_F^2(1-\rho^2)^4} \cdot \frac{\mu_F^2(1-\rho^2)^6}{ 1632000\kappa_F^2(L^2 +\Phi^2)}  \leq \frac{1}{2},\\
\frac{4G^2\eta^2}{\mu_F^2(1-\rho^2)^3} + \frac{8\eta^2}{\mu_F^2(1-\rho^2)^2}\cdot  \frac{960\Phi^2G^2\eta^2 }{\mu_F^2(1-\rho^2)^2} &= \frac{4G^2\eta^2}{\mu_F^2(1-\rho^2)^3} +  \frac{G^2 \eta^2}{\mu_F^2 (1-\rho^2)^4}\cdot \frac{8\cdot 960\Phi^2 \eta^2}{\mu_F^2}\\
&\leq \frac{4G^2\eta^2}{\mu_F^2(1-\rho^2)^3} +  \frac{G^2 \eta^2}{\mu_F^2 (1-\rho^2)^4}\cdot \frac{8\cdot  960\Phi^2 }{\mu_F^2}\cdot \frac{\mu_F^2(1-\rho^2)^6}{1632000 \kappa_{F}^2(L^2 + \Phi^2) } \\
&\leq \frac{5G^2\eta^2}{\mu_F^2(1-\rho^2)^3} .
\end{align*}
Thus we have
\begin{align*}
	\sum_{t=0}^{T} \E{\twonorm{\btheta^{t} - \bone_n \otimes \bbtheta^{t}}^2} & \leq  \frac{16A_1\eta^2}{\mu_F^2(1-\rho^2)^2}  + 2\left( \frac{4G^2\eta^2}{\mu_F^2(1-\rho^2)^3} + \frac{8\eta^2}{\mu_F^2(1-\rho^2)^2}\cdot  \frac{960\Phi^2G^2\eta^2 }{\mu_F^2(1-\rho^2)^2}\right) \sum_{t=0}^{T}\E{\twonorm{\vd^t}^2} \\
	&\leq \frac{16A_1\eta^2}{\mu_F^2(1-\rho^2)^2} + \frac{10G^2\eta^2}{\mu_F^2(1-\rho^2)^3} \sum_{t=0}^{T}\E{\twonorm{\vd^t}^2},
\end{align*}
which completes the proof.

\section{Conclusions} \label{section conclusions}
In this work, we propose a novel decentralized algorithm named MDNPG for MARL. We have established the sample complexity for local convergence of MDNPG, which achieves the best available rate. The key ingredient to our development is a new stochastic ascent inequality for non-convex objectives, which could be of independent interest. Numerical results have demonstrated the efficiency of the proposed method. 

There are several interesting directions for future research. Firstly, it is natural to study the global convergence of MDNPG and extend our framework to the class of entropy-regularized natural policy gradient methods in MARL. 
Secondly,  the  Fisher information matrix in this paper is empirically estimated by sample averaging, which may incur large variance. Thus, we may also consider variance reduction for the estimation of the precondition matrix. Lastly, though  importance sampling is widely used to address the varying data distribution issue when developing variance reduced policy gradient methods, there are also a few recent works  \cite{shen2019hessian,salehkaleybar2022adaptive} which instead use  a hessian-based technique in the single-agent setting. Therefore, it is also interesting to investigate whether  importance sampling  can be removed in the multi-agent setting when developing decentralized (natural) policy gradient methods.

\bibliographystyle{IEEEtran}
\bibliography{refs}
\end{document}